

\magnification = 1200
 at 12truept

\input psfig
\mathsurround = 2pt
\abovedisplayskip=6pt
\belowdisplayskip=6pt

\def \QED {\rlap{$\sqcup$}$\sqcap$\smallskip}

\def\R{{\bf R}}
\def\C{{\bf C}}
\def\D{{\bf D}}
\def\H{{\bf H}}
\let\oeqno\eqno
\def\eqno(#1){\oeqno\hbox{(#1)}}

\def\ref{\hangindent=1pc \hangafter=1 \noindent}
\def\mod{{\rm mod\,}}

\def\cl{{\rm cl}\,}
\def\inter{{\rm int}\,}
\def\orb{{\rm orb}\,}
\def\rel{\;{\rm rel}\;}
\def\dist{{\rm dist}\,}
\def\diam{{\rm diam}\,}

\def\inter{{\rm int}\,}
\def\width{{\rm width}\,}

\def\bd{{\partial}}
\def\capac{{\rm cap}\,}
\def\area{{\rm area}\,}
\def\id{{\rm id}\,}

\def\AA{{\cal A}}

\def\CC{{\cal C}}
\def\DD{{\cal D}}

\def\HH{{\cal H}}
\def\II{{\cal I}}
\def\JJ{{\cal J}}

\def\RR{{\cal R}}
\def\QQ{{\cal Q}}
\def\MM{{\cal M}}
\def\SS{{\cal S}}
\def\YY{{\cal Y}}
\def\VV{{\cal V}}
\def\UU{{\cal U}}
\def\WW{{\cal W}}

\input psfig

\font\eightrm=cmr8

\centerline{\bf Dynamics of quadratic polynomials.}\smallskip
\centerline{\bf I. Combinatorics and geometry of the Yoccoz puzzle.}
\medskip

\centerline{February 5, 1995}
\medskip

\centerline{Mikhail Lyubich
\footnote{*}{\eightrm Supported in part by Sloan Research Fellowship
and NSF grants DMS-8920768 and DMS-9022140.}
}\smallskip

\vskip 40pt
{\narrower
{
Abstract}.
\eightrm
This work studies combinatorics and geometry of the Yoccoz puzzle for
quadratic polynomials. It is proven that the moduli of the ``principal
nest'' of annuli grow at linear rate. As a corollary we obtain complex
a priori bounds and local connectivity of the Julia set for many
infinitely renormalizable quadratics.\par}

\vskip 40pt

\centerline{\S 1. Introduction}\bigskip

This is the beginning of a series of notes on dynamics 
 of quadratic polynomials. 
The forthcoming notes will be concerned with parameter geometry, 
 applications to the
complex and real Rigidity Problems, and measurable dynamics. 
This  first part is a detailed version of \S\S 2,3 of
[L4].

The goal of this part is to study combinatorics and
 geometry of the  Yoccoz puzzle,
with an application to the problem of local connectivity of Julia sets.
The puzzle was introduced by Branner and Hubbard [BH]
for cubic maps with one escaping critical point 
and by Yoccoz for quadratics (see [H], [M2]).
 The main geometric result of these works is the divergence property
of moduli of a certain nest of annuli. This implies that
the corresponding domains (``puzzle pieces") shrink to points, which
yields local connectivity of the Julia set. The corresponding
 result in the parameter plane yields rigidity.

The geometric result of Branner-Hubbard and Yoccoz does not contain
information on the {\it rate} at which the pieces shrink to points. In this
work we tackle this problem. We consider a smaller nest $V^0\supset V^1\supset\ldots$
called {\it principal}, and prove that the annuli moduli of this 
nest grow at linear rate over a certain combinatorially specified
 subsequence of levels (Theorem III).
This implies that  the first return maps
 $g_n: V^n\rightarrow V^{n-1}$ 
are becoming purely quadratic at
exponential rate. We expect this fact to have many applications.

Theorem III also yields {\it a priori bounds} for infinitely 
renormalizable quadratics of ``sufficiently big type"  (Theorems IV and IV$'$).
For  real quadratics of ``bounded type" the problem of complex  bounds
was resolved by Sullivan (see [S2], [MvS]). Our result for real quadratics
is a kind of complement to
Sullivan's, though there is still
one special combinatorial situation,  ``saddle-node cascades", which  is
not covered by either of these results.
In the forthcoming notes we
will give an appropriate extension of Sullivan's  bounds which
pick saddle-node cascades [LY]: Local connectivity 
of the Julia sets for all real quadratics follows.\footnote{**}{\eightrm This result
with a different proof has been also announced by Levin and  van Strien [LS].}
As to the non-real situation, to the best of our knowledge
 the problem of a priori bounds was not settled
before for any one quadratic of bounded type (compare Rees [R]).

An immediate consequence of our a priori bounds is that the puzzle pieces
of the ``full principal nest" shrink to points, which yields local
connectivity of the Julia set. 
This is a nice property, since  such a
Julia set has an explicit  topological model (see Douady [D2]).
Note that a large pool of locally connected examples was already known (Douady - Hubbard,
Yoccoz (see [DH1], [H]), Hu - Jiang [HJ], [J], Petersen [P]).
Unfortunately there  are also counter-examples:
Cremer and some infinitely renormalizable quadratics have non-locally
connected  Julia sets (see [L3], [M2]). 

 However, the main reason why  complex a priori bounds are important is that
they provide a key to the rigidity problem and the
 renormalization theory
(compare Sullivan [S2] and McMullen [McM]).
These applications will be the subject of the forthcoming notes.

Let us now describe the structure of the paper.

In the next section, \S 2,
 we overview the necessary preliminaries in holomorphic 
dynamics, particularly Douady-Hubbard renormalization and the Yoccoz puzzle.

In \S 3 we present our approach to combinatorics of the puzzle.
The main concepts involved are 
the {\it  principal nest}
 of puzzle pieces, {\it generalized renormalization} and 
{\it central cascades.}
As we indicated above, the principal nest $V^0\supset V^1\supset\ldots$
 contains the key combinatorial and geometric
information about the puzzle. We describe the combinatorics of this nest
by means of {\it generalized renormalizations}, 
that is, 
appropriately restricted first return maps considered up to rescaling.

It may happen that the quadratic-like map $g_n: V^n\rightarrow V^{n-1}$ has "almost
connected" Julia set. 
 This phenomenon often requires a special treatment.
Such a map  generates a subnest of
the principal nest called a {\it central cascade.}
The number of central cascades  in the principal nest 
 is called the {\it height} $\chi(f)$  of a map $f$. In other words,
$\chi(f)$ is the number of {\it different} quadratic-like maps among the $g_n$'s
(where "different" means: "with different Julia sets").
 It will play an important role for our discussion.

In \S 4 we study the initial geometry of the puzzle.
The main result of this section is the construction of an initial annulus
$V^0\backslash V^1$  
with definite modulus, provided the hybrid class of a map is
selected from a truncated secondary limb (Theorem I). 

In \S 5 we prove the main geometric result of the paper saying
 that the principal moduli  grow linearly with the number of central cascades.  
We introduce a new geometric parameter (worked out jointly with J. Kahn),
 {\it the asymmetric modulus},
and prove first that it is {\it monotonically non-decreasing }
when we go down along the 
principal nest (Theorem II). This
already provides us
with lower bounds for the principal moduli $\mu_m=\mod(V^{m-1}\backslash V^m)$
(which, by  the way, implies the  Branner-Hubbard-Yoccoz divergence property),
and upper bounds on the distortion.
We reach these results by means of a purely combinatorial analysis
plus the standard  Gr\"{o}tzsch inequality. 

However this analysis does not  always
yield the linear growth of moduli. In particular, it is not enough for
the basic example called the {\it Fibonacci map}. 
The crucial part of the proof is to gain this extra growth for
Fibonacci-like combinatorics.
The argument based on the {\it Definite Gr\"{o}tzsch inequality} involves
estimates of hyperbolic  distances between puzzle pieces and
analysis of their shapes. The key observation is that sufficiently
pinched pieces make
 a definite  extra contribution  to  the moduli growth (Theorem III).

At the end of this section we prove  a priori bounds for infinitely
renormalizable quadratics of sufficiently big type (Theorems  IV and
IV$'$). 
The meaning of this condition is that certain combinatorial parameters
of the renormalized maps $R^n f$ are sufficiently big (depending on
the truncated secondary limb to which the internal class of $R^n f$
belongs). The main such parameter is the above mentioned height,
but there are also several others.
These conditions together  mean roughly that the {\it periods}
of $R^n f$ are sufficiently big. The only 
difference is a possibility of 
long ``parabolic or Siegel cascades".

In the next short section, \S 6, we show that the Julia sets of the quadratics
under consideration are  locally connected. This follows from shrinking
of the puzzle-pieces.

In the last section, \S 7, we outline the content of the forthcoming notes.

Before finishing  let us draw the reader's attention to the work of
 Graczyk and Swiatek [GS] containing a related result on the 
growth of moduli (within a certain  nest of domains different from the 
principal nest) for {\it real} quadratics. 

\medskip\noindent {\bf Acknowledgement.} I would like to thank Jeremy Kahn
for a fruitful suggestion, and Curt McMullen, Mary Rees and
Mitsuhiro Shishikura for useful comments on the results. 
I also had  useful discussions with John Milnor on combinatorics of external
rays. 
Feedback from the Stony Brook dynamical group during my course in the fall
1994 was very helpful for cleaning up the exposition. 
Let me also thank  Scott Sutherland and Brian Yarrington for help with the 
computer pictures.
The results of this work were obtained
in June 1993 during the Warwick Workshop on
hyperbolic geometry. I am grateful to 
the organizers, particularly  David Epstein and Caroline Series, 
for that wonderful time.

\bigskip
\centerline{\bf \S 2. Preliminaries: Douady-Hubbard renormalization and 
   Yoccoz puzzle.}\bigskip

\noindent{\bf 2.1. General terminology and notations.}
 Given two subsets $V$ and $W$ of the complex plane, 
 we say that  $V$ is {\sl strictly
contained} in $W$, $V\subset\subset W$,  if
$\cl V\subset \inter W$. 

Let $B(z,\epsilon)$ denote the disk of radius $\epsilon$ centered at $z$.
Given two sets $A$ and $B$, 
let $\dist(A,B)=\inf\{ \dist(z,\zeta): z\in A, \zeta\in B\}$.

By $\orb z$ we denote the forward orbit $\{f^n z\}_{n=0}^{\infty}$ of $z$,
and by $\omega(z)$ its $\omega$-limit set. Let also orb$_n z=\{f^m z\}_{m=0}^n$.  
By a {\sl topological disk} we will mean a simply connected
region in ${\bf C}$. By an {\it annulus} we mean a doubly connected
region. A {\it horizontal curve} in an annulus $A$ is a preimage of 
a circle centered at 0 by the Riemann mapping $A\rightarrow \{z: 1<|z|<R\}.$ 

We assume that the reader is familiar with the basic theory of quasi-conformal
maps (see [A]).
Quasi-conformal and quasi-symmetric maps will be abbreviated as qc and qs
correspondingly.

\medskip\noindent{\bf 2.2. Polynomials.}
 By now there are many surveys and books on
holomorphic dynamics. The reader can consult   [Be], [GC], [M1]
 for general reference,
and [B], [DH1] for the quadratic case. Below we will state the main 
definitions and facts required  for  discussion. However we assume
that the reader is familiar with the notions of periodic points, and their
classification as attracting, neutral, parabolic and repelling.

 Let $f: \C\rightarrow \C$ be a polynomial
of degree $d\geq 2$. The {\sl basin of} $\infty$ is 
the set of points escaping to $\infty$:
$$D_f(\infty)\equiv D(\infty)=\{z\in \C: f^n z\to \infty\}.$$
Its complement is called {\sl the filled Julia set}:
     $$K(f)=\C\backslash D(\infty).$$ 
The {\sl Julia set} is the common boundary of $K(f)$ and $D(\infty)$:
$$J(f)=\partial K(f)=\partial D(\infty).$$ The Fatou set $F(f)$ is
defined as $\C\backslash J(f)$.
The Julia set (and the filled Julia set) is connected if and only if
non of the critical points  escape to $\infty$, that is, all of them belong
to $K(f)$.

Given a polynomial $f$, there is a conformal map 
(the {\it B\"{o}ttcher function})
 $$B_f: U_f\rightarrow \{ z: |z|>r_f\geq 1\}$$ of a neighborhood $U_f$ of 
infinity onto the exterior of a disk such that
$B_f(fz)=(B_f z)^d$ and $B_f(z)\sim z$ as $z\to\infty$.
There is an
explicit dynamical formula for this map:
$$B_f(z)=\lim_{n\to\infty} (f^n z)^{1/d^n} \eqno (2-1)$$
with an appropriate choice of the branch of the $d^n$th  root.

 If the Julia set 
$J(f)$ is connected then $\partial U_f$
contains  a critical point $b$ of $f$.   
Otherwise $B_f$ coincides with he Riemann mapping of the whole basin of infinity
$D(\infty)$ onto $\{z: |z|>1\}$ (in this case
$r_f=1$). 
 
The {\sl external rays}
$ R_{\theta}$  with angle $\theta$ and
 {\sl equipotentials} $E_h$ of level $h$ are defined as the
 $B_f$-preimages of the straight rays 
$\{ re^{i\theta}: r_f<r<\infty\}$
and the round circles   $\{ e^{i\theta}: 0\leq\theta\leq 2\pi\}$.
They form two orthogonal invariant foliations of $D(\infty)$.

\proclaim Theorem 2.1 (see [M1], \S 18, or [H]). If $a$ is a
 repelling periodic point of $f$, 
then there is at least one but at most finitely many
external rays landing at $a$. 

\noindent{\bf Sketch.}
Let us sketch the geometric construction of these rays (see the above references
for the details). Let us linearize $f^p$ near the periodic point $a$:
Let $\psi: (U,0)\rightarrow (V, a)$ be a local conformal map  normalized by
$\psi'(0)=1$ and such that
$$\psi (\lambda z)=f^p (\psi z) \eqno (2-2),$$
where $\lambda$ is the multiplier at $a$.
By means of this functional equation $\psi$ can be analytically
extended to the whole complex plane. Thus we obtain an entire function
$\psi : \C\rightarrow \C$ satisfying (2-2).

Let $\DD_i$ be the components of  $\psi^{-1} D(\infty)$. It is not hard
 to see that
$\psi: \DD_i\rightarrow D(\infty)$ is a universal covering map. Moreover,
the function $$h_i=\log B_f\circ \psi: \DD_i\rightarrow \H \eqno (2-3)$$
 gives the Riemann
map of $\DD_i$ onto the upper half-plane $\H$. 

The key point  is that there are only finitely many components 
$\DD_i$ which hence invariant under some iterate $f^p$. 
Then the map $h_i$ conjugates 
$z\mapsto \lambda^p z$ on $D_i$ to $z\mapsto d^p z$ on $\H$. 
 Take now the hyperbolic
geodesic $\gamma$  joining 0 and $\infty$ in $\H$ 
(the vertical ray) and the corresponding
hyperbolic geodesic $\hat R_i=h_i^{-1} \gamma $ in $\DD_i$. It follows that
 $\hat R_i$ is invariant with respect to $z\mapsto \lambda^p z$, 
and hence lands at 0.  Hence the external 
ray $R_i=\psi \hat R_i$ lands at $a$.\QED

Thus the external rays landing at $a$ are organized in several cycles.
The rotation number of these cycles is the same, and
 is called the {\sl combinatorial rotation number} $\rho(a)$ of $a$.
Let $\RR(a)$ denote the union of the external rays landing at $a$,
and 
$$\RR(\bar a)=\cup_{k=0}^{p-1} \RR(f^k a)$$
(where $p$ is the period of $a$). This configuration (with the external angles
marked on the  rays) is called the {\it
rays portrait} of the cycle $\bar a$. The class of isotopic portraites
is called the {\it abstract rays portrait.}

\medskip\noindent{\bf 2.3. Quadratic family.}
Let now $f\equiv P_c: z\mapsto z^2+c$ be the quadratic family. In this case
the rays portraites of periodic cycles have quite specific combinatorial
properties. The reader can consult [DH1], [At], [GM], [Sch], [M4] 
for the proofs of the results quated below.\medskip

{\bf Proposition 2.2 (see [M4]).} {\it Let  $\bar a=\{a_k\}_{k=0}^{p-1}$
 be a repelling periodic cycle such that there are at least two rays
landing at each point $a_k$. \smallskip

 (i) Let
 $S_1$ be the components of $\C\backslash \RR(\bar a)$ 
containing the critical value
$c$. Then $S_1$ is a sector bounded by two external rays.\smallskip

(ii) Let $S_0$ be the component of $\C\backslash f^{-1} \RR(\bar a)$ containing the
  critical point 0. Then $S_0$ is bounded by four external rays: two of them
  lands at a periodic point $a_k$, and two others land at the symmetric point
  $-a_k$.  \smallskip

(iii) The rays of $\RR(\bar a)$ form either one or two cycles under
  iterates of $f$.} \medskip

A particular situation of such kind is the following. Let 
$\bar b=\{b_k\}_{k=0}^{p-1}$ be an attracting cycle, $p>1$. Let $D_k$ be the 
components of its basin of attraction containing $b_k$. Then the 
boundaries of $D_k$ are Jordan curves, and the restrictions $f^p|D_k$
are topologically conjugate to the doubling map $z\mapsto z^2$ of the
unit circle. Hence there is a unique $f^p$-fixed point 
$a_k\in D_k$. Altogether these points form a repelling periodic cycle $\bar a$
(whose period may be smaller than $p$), with at least two rays landing at
each $a_k$. The  portrait $\RR(\bar a)$ will be also called the rays
portrait associated to the attracting cycle $\bar b$.

A case of special interest for what follows is the fixed points portraites.
There is always a fixed point called $\beta$ which is the landing point of
the invariant ray $\RR_0$. Moreover, this is the only ray landing at this
point, so that it is non-dividing: the set $K(f)\backslash \{\beta\}$ is connected.

If the second fixed point called $\alpha$ is also repelling, it is much
more interesting. Namely, this point is dividing, so that there are 
several external rays landing at it. These rays are cyclically permuted
by dynamics with some combinatorial rotation number $q/p$.

The 
{\sl Mandelbrot set} $M$ is defined as the set of $c\in \C$ for which
$J(P_c)$ is connected, that is, 0 does not escape to $\infty$ under
iterates of $P_c$. If $c\in \C\backslash M$, then $J(P_c)$ is a Cantor set. 

The Mandelbrot set itself is connected (see [DH1], [CG]). This is proven by
constructing explicitly the Riemann mapping 
$B_M :  \C\backslash  M\rightarrow \{z: |z|>1\}$. Namely, let $D_c(\infty)$ be the 
basin of $\infty$ of $P_c$, and  $B_c$ be the B\"{o}ttcher
function  (2-1) of $P_c$. Then 
$$B_M(c)=B_c(c).   \eqno (2-4)$$
The meaning of this formula is that the ``conformal position" of a 
parameter $c\in \C\backslash M$ coincides with the ``conformal position" of
the critical value $c$ in the basin $D_c(\infty)$. This relation is the key
to the similarity between dynamical and parameter planes.

Using the Riemann mapping $B_M$ we can define the {\it parameter
external rays} and {\it equipotentials} as the preimages of the
straight rays going to $\infty$ and round circles centered at 0. 
This gives us two orthogonal foliations in the complement of the
Mandelbrot set.

A quadratic polynomial $P_c$ with $c\in M$ 
is called {\sl hyperbolic} if it has an attracting cycle. The set of hyperbolic
parameter values is the union of some components of $\inter M$ called 
{\sl hyperbolic components}. Conjectually all components of $\inter M$ are
hyperbolic. This Conjecture would follow from the MLC Conjecture
(Douady \& Hubbard [DH1]).

Let $H\subset \inter M$ be hyperbolic component of the Mandelbrot set,
and let $\bar b(c)=\{b_k(c)\}_{k=0}^{p-1}$ be the corresponding 
attracting cycle.
On the boundary of $H$ the cycle $\bar b$ becomes neutral, and there is 
a single point $c_H\in \partial H$ where $f^p(b_0)=1$ [DH1]. This  point
is called the {\it root } of $H$. 

Given a $c\in H$, there is the rays portrait $\RR_c$ associated to the
corresponding  attracting periodic point. Let $\theta_1$ and $\theta_2$ 
be the external angles of the two rays bounding the sector $S_1$ of
Proposition 2.2. 

\proclaim Theorem 2.3 (see [DH1], [M4], [Sch]). 
The parameter rays $\RR^H_1$ and $\RR^H_2$
with angles $\theta_1$ and
$\theta_2$ land at the root $c_H$. There are no other rays landing at
$c_H$.
 
The region $W_H$ in the parameter plane bounded by the rays $\RR^H_i$ and 
containing $H$ is called the {\it wake} (originated from $H$).
 The part of the Mandelbrot set
contained in the wake together with the root  $c_H$ is called the
{\it limb } $L_H$ of the Mandelbrot set originated at $H$.
 The root of $H$
is also called the root of wake $W_H$ and limb $L_H$. 

Recall that for $c\in H$, $\bar a_c$ denotes the repelling cycle
associated to the basin of attraction of the attracting cycle $\bar b_c$.
 The dynamical
meaning of the wakes is reflected in the following statement.

\proclaim Proposition 2.4 (see [GM]). 
The repelling cycle $\bar a_c$  stays repelling
throughout the wake $W_H$.
 The corresponding rays portrait $\RR(a_c)$ preserves
its isotopic type throughout  this wake.

The {\sl main cardioid} of $M$ is defined as the set of points $c$ for 
which $P_c$ has a neutral fixed point $\alpha_c$, that is, 
$|P_c'(\alpha_c)|=1$. It encloses the hyperbolic component where
$P_c$ has an attracting fixed point. In the exterior of the 
 main cardioid both fixed points are repelling.

The limbs attached to the main cardioid are called {\it primary}.
Let $H$ be a hyperbolic component attached to the main cardioid.
The limbs attached to such a  component are called {\it secondary}.
We refer to a {\sl truncated limb} if we remove from it a neighborhood
of its root  (Figure~1).




We refer to a {\sl truncated limb} if we remove from it a neighborhood
of its root  (Figure~1).

\midinsert
\centerline{\psfig{figure=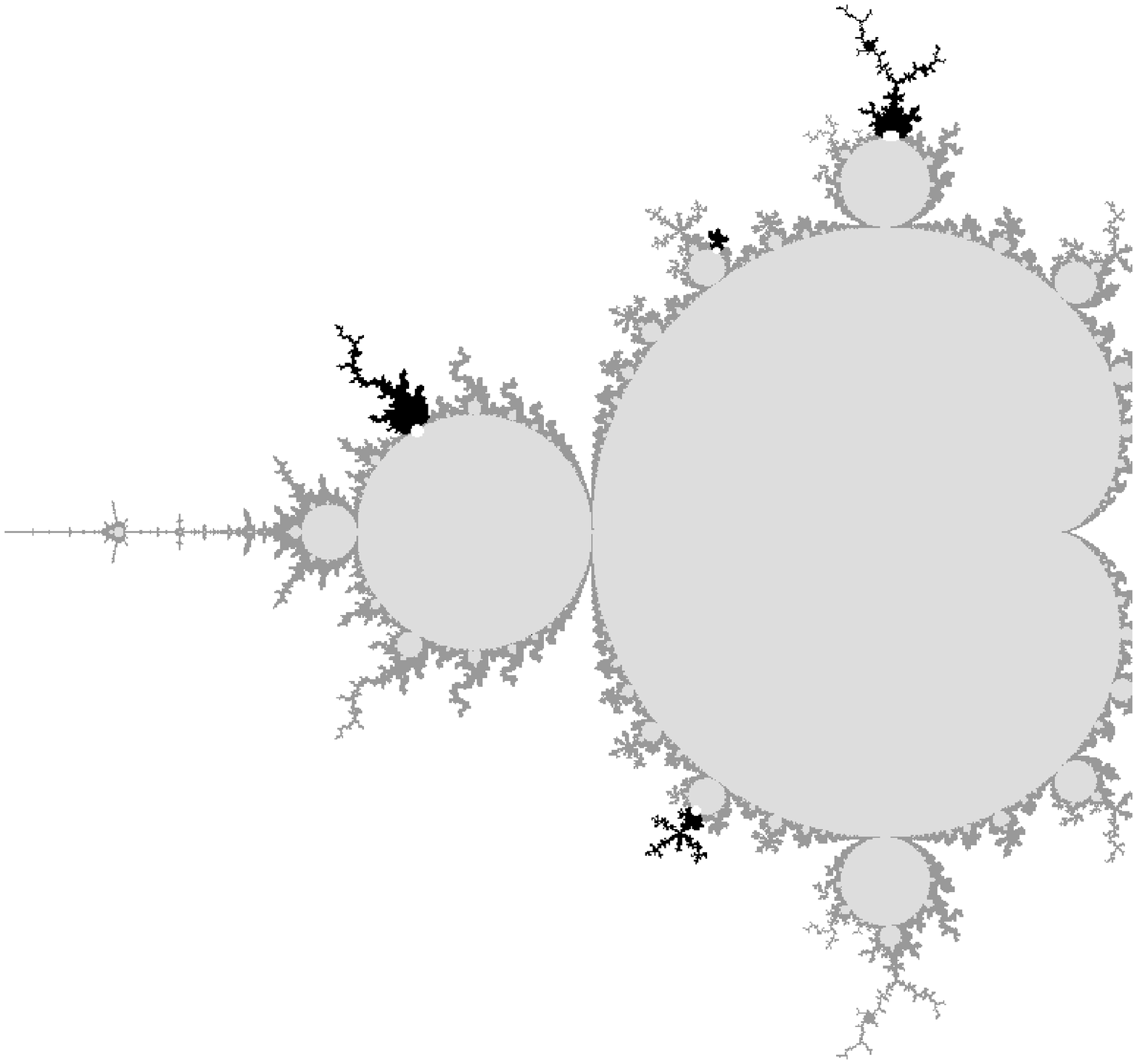,width=.8\hsize}}
\bigskip\centerline{Figure 1: Truncated secondary limbs of the Mandelbrot set.}
\medskip\endinsert 

\medskip\noindent {\bf 2.4. Douady-Hubbard polynomial-like maps.}
 The main reference
for the following material is [DH2]. 
Let $U'\subset\subset U$ be two topological disks (with a non-degenerate 
annulus $U\backslash U'$). A branched covering   $f: U'\rightarrow U$
is called 
a {\sl DH polynomial-like map} (we will sometimes skip ``DH" in case this
does not cause confusion with ``generalized" polynomial-like maps
defined below). Every polynomial with connected Julia set
can be viewed as a 
polynomial-like map after restricting it onto an appropriate neighborhood
of the filled Julia set. Polynomial-like maps of degree 2 are
called {\sl (DH) quadratic-like}. We will always normalize quadratic-like maps
so that the origin 0 will be its critical point.  

One can naturally define the filled Julia set
of $f$ as the set of non-escaping points:
 $$K(f)=\{z: f^n z\in U': n=0,1,\ldots\}.$$
The Julia set is defined as $J(f)=\partial K(f)$. 
These sets are connected if and only if non of the critical points is
escaping.


The choice of the domain $U'$ and range $U$ of a polynomial-like map is
not canonical. Two polynomial-like maps $f: U'\rightarrow U$ and
 $g: V'\rightarrow V$
are considered to be the same if $K(f)=K(g)$ and $f$ coincides with $g$
on a neighborhood of $K(f)$.

Given a polynomial-like map $f: U'\rightarrow U$, we can consider a
{\sl fundamental annulus} $A=U'\backslash U$. It is certainly not a canonical object
but rather depending on the choice of $U'$ and $U$. Let
 $$\mod(f)=\sup \mod A,$$
where $A$ runs over all fundamental annuli of $f$.

Two polynomial-like maps $f$ and $g$ are called
 {\sl topologically (quasi-conformally, conformally) conjugate} if 
 there is a homeomorphism (qc map, conformal isomorphism correspondingly)
 $h$ from a neighborhood of $K(f)$ to a neighborhood of $K(g)$ such that
 $h\circ f=g\circ h$.    

If there is a qc conjugacy $h$ between $f$ and $g$ with $\bar\partial h=0$
almost everywhere on the filled Julia set $K(f)$, then $f$ and $g$ are
called {\sl hybrid} or {\sl internally}  equivalent. A {\sl hybrid class}
$\HH(f)$ is the space of DH polynomial-like maps hybrid equivalent
to $f$ modulo conformal equivalence. According to
Sullivan  [S1], a hybrid class of polynomial-like maps
should be viewed as  as infinitely dimensional
Teichm\"{u}ller space. 
In contrast with the classical Teichm\"{u}ller theory
 this space has a preferred point:

\proclaim Straightening Theorem [DH2]. Any hybrid class $\HH(f)$
of DH polynomial-like maps
with connected Julia set contains a unique (up to affine conjugacy)
  polynomial. 

In particular,
any hybrid class of quadratic-like maps with connected Julia set
 contains a unique quadratic polynomial
$z\mapsto z^2+c$ with $c=c(f)\in M$. 
So the hybrid classes of quadratic-like maps are labeled by the points
of the Mandelbrot set. In what follows we will freely identify a  
quadratic hybrid class with its label $c\in M$.  

Sullivan supplied any hybrid class with the following
{\sl Teichm\"{u}ller metric} [S1]:
$$\dist_T(f,g)=\inf \log K_h,$$
where $h$ runs over all hybrid conjugacies between $f$ and $g$,
and $K_h$ denotes the qc dilatation of $h$. 
The Teichm\"{u}ller distance from $f$  to the  quadratic 
$P_{c(f)}: z\mapsto z^2+c(f)$ in its hybrid class is controlled by the
modulus of $f$:

\proclaim Proposition 2.5. If $\mod (f)\geq\mu>0$ then 
 $\dist_T (f, P_{c(f)})\leq C$ with a $C=C(\mu)$ depending only on
$\mu$. Moreover, $C(\mu)\to 0$ as $\mu\to\infty.$

This is the reason why the control of moduli of polynomial-like maps 
is crucial  for the renormalization theory (see [S2], [McM]).   

Let us finish with the following remark:
Given a polynomial-like map with connected Julia set,
 we can define external rays and equipotentials
near the filled Julia set by conjugating it to a polynomial. This definition
is certainly not canonical but rather depends on the choice of conjugacy.
We will specify the particular choice whenever it matters.

\medskip\noindent{\bf 2.5. Douady-Hubbard renormalization. } 
The reverse
 procedure under
the name of {\sl tuning} is discussed in [DH2], [D] and [M3]. A more general
point of view (but which is equivalent to the tuning, after all)
is presented in [McM].

Let $f: U'\rightarrow U$ be a quadratic-like map. Let $\bar a$ be
a dividing repelling  cycle, so that there are
 at least two rays landing at each point of $\bar a$.
Let $\RR\equiv \RR(\bar a)$ denote the configuration of rays landing at
$\bar a$, and let $\RR'=-\RR$ be the symmetic configuration.
Let us aslo consider an arbitrary equipotential $E$.
Let now 
$\Omega$ be the component of 
$\C\backslash (E\cup \RR\cup \RR')$
containing the critical point 0.   By Proposition 2.2,
it is bounded by four arcs $\gamma_i$ of
external rays ant two pieces of the equipotential $E$.

\midinsert
\centerline{\psfig{figure=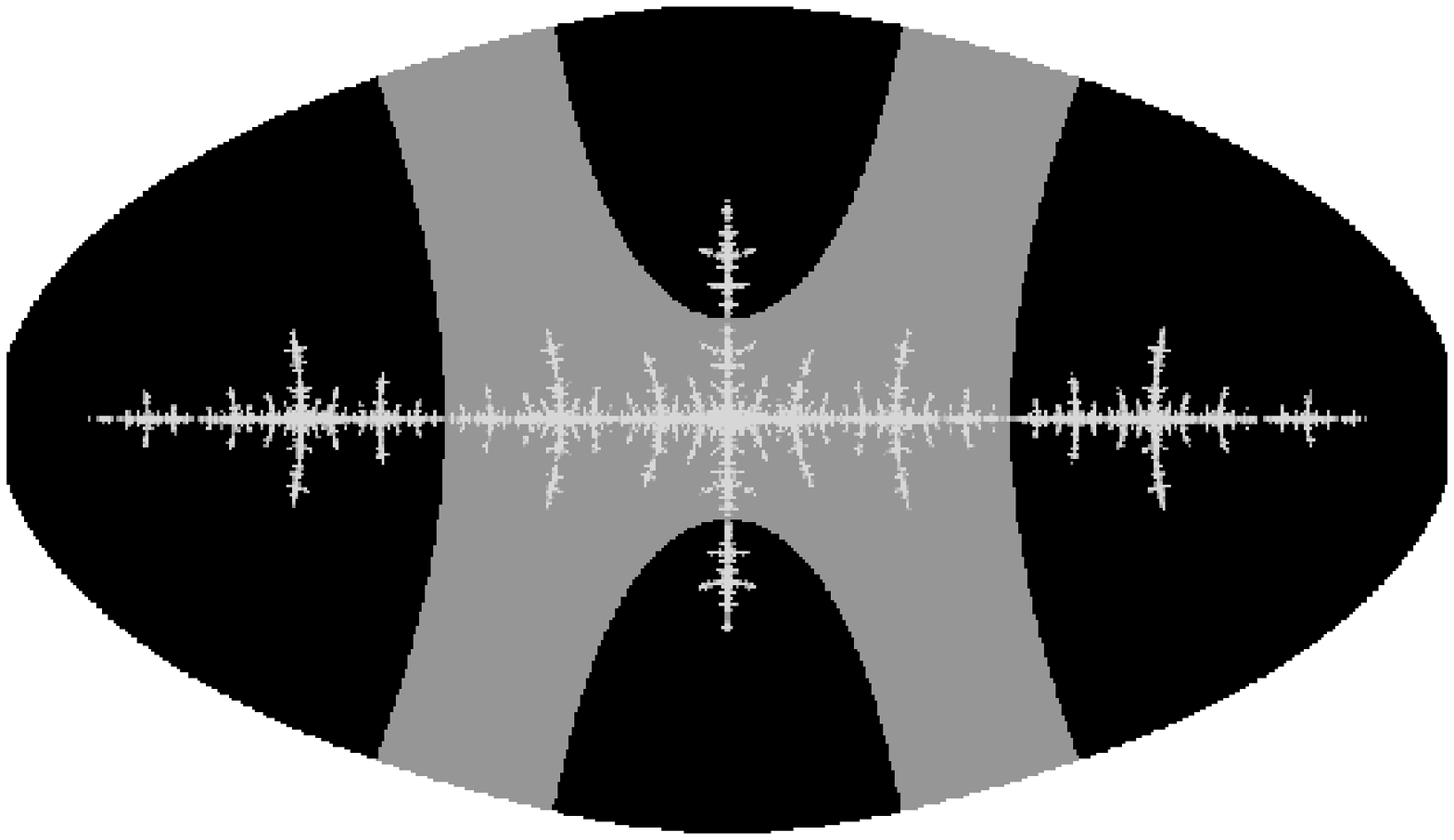,width=.8\hsize}}
\bigskip\centerline{Figure 2: A renormalization domain for the
    Feigenbaum polynomial.}
\medskip\endinsert 

Let $p$ be the period of the above rays, and $a\in \partial \Omega$ be the
periodic point of cycle $\bar a$. Let us consider a domain 
$\Omega'\subset \Omega$, 
the component of $f^{-p} \Omega$ attached to $a$ (see Figure 2).
If $\Omega'\ni 0$ then $f^p: \Omega'\rightarrow \Omega$ is a double covering
map. 
A quadratic-like map $f$ is called {\sl DH-renormalizable}
 if there is a repelling
cycle  $\bar a$ as above such that $\Omega'\ni 0$, 
and 0 does not escape $\Omega'$ under iterates of 
$f^p$. We will also say
that this renormalization is associated to the periodic point $a$.
We call $f$ 
{\sl immediately DH renormalizable} if $a$ can be selected as the
dividing fixed point $\alpha$ of $f$.

Note that the disks $\Omega', f\Omega',\ldots, f^{p-1}\Omega'$ have
disjoint interiors. Indeed, otherwise  $ f^k\Omega'$ would
 be inside
$\Omega$ for some $k<p$. But  this is impossible since the external
rays which bound $f^k \Omega'$ are outside of $\Omega$.

In the DH-renormalizable case by means of a ``thickening procedure"
(see [DH1] or [M2]) one can extract a polynomial-like map
$f^p: V'\rightarrow V$. Namely,  let us consider
a little bit bigger domain $V\supset \Omega$ 
bounded by arcs of four external rays
close to $\gamma_i$, two arcs of  circles going around the points $a$ and 
$a'=-a$,
and two arcs of $E$.  Pulling $V$ back by $f^p$, we 
obtain a domain $V'\subset\subset V$ such that 
the map $f^p: V'\rightarrow V$ is quadratic-like. This map considered up to
conformal conjugacy is called the {\it DH renormalization} of $f$.

Let now  $f: z\mapsto z^2+c_0$ be a quadratic polynomial, $c_0\in M$. 
If it is renormalizable then there is a homeomorphic copy $M_0\ni c_0$ 
of the Mandelbrot set with the following properties (see [DH2, [D]). For 
$z\in M_0'= M_0\backslash\{one\; point\}$ the polynomial
$P_c: z\mapsto z^2+c$ is renormalizable. Moreover, 
there is the parameter analytic
extension $a_c$ of the periodic point $a$ to a neighborhood of
$M_0'$ such that the above renormalization of $P_c$ is associated to $a_c$.
At the parameter value $b$ removed from $M_0$ the periodic point $a_c$
is becoming parabolic with multiplier one.  This parameter value is called
the {\sl root} of $M_0$.  The component $H$ of $M_0$ corresponding to the
component of $M$ enclosed by the main cardioid ``gives  origin" for the
copy $M_0$. 
Vice versa, any hyperbolic component $H$ of the Mandelbrot set gives
origin to a copy of $M$.  In particular, the copies  corresponding
to the immediate renormalization are attached to the main cardioid.

We will see below that among all renormalizations there is the {\sl first}
 one, which we denote $Rf$ (see \S 3.4). This renormalization corresponds
to a  {\sl maximal} copy  of the Mandelbrot (that is a copy, which
 is not contained in
any bigger copies except $M$ itself).
 Let $\MM$ denote the family of maximal Mandelbrot copies.

Given any sequence $M_0, M_1,\ldots$ of maximal copies of $M$,
there is an infinitely renormalizable quadratic polynomial $P_b$ such that
$c(R^m P_b)\in M_m$, $m=0,1,\ldots$. Indeed, the sets 
$$C_N=\{b: c(R^m P_b)\in M_m,\;  m=0,1,\ldots, N\}$$
form a nest of copies  of $M$ whose intersection  consists of the desired
parameter values. 

We say that these infinitely renormalizable quadratics have the same
combinatorics determined by the sequence $M_0, M_1,\ldots$. 
The MLC problem for these parameter values is equivalent to the
assertion that there is only one quadratic with a given combinatorics.
In other words, the above copies $C_N$ shrink to a point as $N\to \infty$.

\medskip\noindent {\bf 2.6. Yoccoz puzzle.} Let $f: U'\rightarrow U$ be a 
quadratic-like map with both fixed points $\alpha$ and $\beta$
repelling. As usual, $\alpha$ denotes the dividing fixed point with
rotation number $\rho(\alpha)=q/p$, $p>1$.
Let $E$ be an equipotential
sufficiently close to $K(f)$ (so that both $E$ and $fE$ are closed curves).  
   Let $\RR_{\alpha}$ denote the union of external rays landing
at $\alpha$. 
These rays  cut the domain bounded by $E$ into
$p$ closed topological disks $Y^{(0)}_i$, $i=0,\ldots , p-1$, called {\sl puzzle
pieces of zero depth} (Figure 3). The main property of this partition is that
$f\partial (\cup Y_i^{(0)})$ is outside of $\inter (\cup Y_i^{(0)})$.

Let us now define {\sl puzzle-pieces $Y^{(n)}_i$ of depth $n$} as the connected 
components of $f^{-n} Y^{(0)}_k$. They form a finite partition of the
the neighborhood of $K(f)$ bounded by $f^{-n} E$.
If the critical orbit does not land at $\alpha$, then
 for every depth there is
a single puzzle-piece  containing the critical point. It is called
{\sl critical} and is labeled as $Y^{(n)}\equiv Y^{(n)}_0$.

Let $\MM(f)$ denote 
the family of all puzzle pieces  of $f$ of all levels. 
It is {\sl Markov} in the following sense:

\smallskip\noindent (i) Any two puzzle pieces are either nested or have
  disjoint interiors. In the former case the puzzle-piece of bigger
  depth is contained in the one of smaller depth. 

\smallskip\noindent (ii) The image of any puzzle-piece $Y^{(n)}_i$ of depth
  $n>0$ is a puzzle piece $Y^{(n-1)}_k$ of the previous depth.  Moreover,
   $f: Y^{(n)}_i\rightarrow Y^{(n-1)}_k$ is a two-to-one branched covering
   or a  conformal isomorphism
   depending on whether $Y^{(n)}_i$ is critical or not.

We say that $f^k Y^{(n)}_i$ $l$-{\sl to-one covers} a union of pieces
$\cup_{m,j} Y^{(m)}_j$ if $f^k|\inter Y^{(n)}_i$ is $l$-to-one covering map
onto  its image, and
$$f^k| Y^{(n)}_i\cap J(f)=\bigcup_{m,j} Y^{(m)}_j\cap J(f).$$
In this case $\cup Y^{(m)}_j$ is obtained from $f^k| Y^{(n)}_i$ by cutting with
appropriate equipotentials.

On depth 1 we have $2p-1$ puzzle pieces: one central $Y^{(1)}$,
$p-1$ non-central $Y_i^{(1)}$ attached to the fixed point $\alpha$ 
(cuts of $Y^{(0)}_i$ by the equipotential $f^{-1} E$), and $p-1$ symmetric ones
$Z^1_i$ attached to $\alpha'$.  Moreover, $fY^{(1)}$ two-to-one covers 
$Y^{(1)}_1$, $Y^{(1)}_i$ univalently covers $Y^{(1)}_{i+1}$, $i=1,\ldots , p-2$,
and $f Y^{(1)}_{p-1}$ univalently covers $Y^{(1)}\cup_i Z^{(1)}_i$.
Thus $f^p Y^{(1)}$ cut by $f^{-1} E$
is the union of $Y^{(1)}$ and $Z^{(1)}_i$ (Figure 3). 

\proclaim Theorem 2.6. (Yoccoz, 1990). Assume that both fixed points of 
a polynomial-like map $f$ are repelling,
 and that $f$ is DH non-renormalizable.
Then the following divergence property holds:
$$\sum_{n=0}^{\infty}\mod (Y^{(n)}\backslash Y^{(n+1)})=\infty.$$
Hence $\diam Y^{(n)}\to 0$ as $n\to\infty$.

\proclaim Corollary 2.7. Under the circumstances of the above theorem
   the Julia set $J(f)$ is locally connected.

The reader can consult [H], [M] or [L2]  for a proof 
(or read \S 5 of this paper).

The Yoccoz puzzle provides us with the Markov family of puzzle pieces to
play with. There are several different ways to do this: by means of
the Branner-Hubbard tableaux [BH], or by means of the Yoccoz $\tau$-function
(unpublished), or by means of the principal nest and
 generalized renormalization, as  will be
described below (compare [L2], [L3]).

\bigskip
\centerline{\bf \S 3. Principal nest and generalized renormalization.}\medskip 

{\sl In what follows we will always assume that both fixed points of the
DH quadratic-like maps under consideration are repelling. }

\medskip\noindent{\bf 3.1 Principal nest.} Given a set $W=\cl(\inter W)$ 
 and a point $z$  such that $f^l z\in \inter W$,
 let us define the {\sl pull-back}
of $W$ along the ${\rm orb}_l z$ as the chain of sets 
$W_0=W, W_{-1}\ni f^{n-1}z,\ldots, W_{-l}\ni z$ such that 
$W_{-k}$ is the closure of the component of $f^{-1} (\inter W_{-k+1})$ containing
$f^{l-k}z$. In particular if $z\in \inter W$ and $l>0$ is the moment of 
first return of $\orb z$ back to $\inter W$ we will refer to
 the {\sl pull-backs corresponding
to the first return of} $\orb z$ to $\inter W$.

Let us consider the puzzle pieces of depth 1 as described above: 
 $Y^{(1)}, Y^{(1)}_i$ and $Z^{(1)}_i$, $i=1,\ldots p-1$ (Figure 3). 
 If $z\in Y^{(1)}$ then
$f^p z$ is either in $Y^{(1)}$ or in one of $Z^{(1)}_i$. Hence either
$f^{pk}0\in Y^{(1)}$ for all $k=0,1, \ldots$,  
or there is the smallest $t>0$ and a $\nu$ such that
$f^{tp}0\in Z^{(1)}_{\nu}$. Thus either $f$ is immediately DH renormalizable, or
the critical point escapes through one of the non-critical pieces, attached
to $\alpha'$.

In the immediately renormalizable case the principal nest of puzzle pieces
consists of just single puzzle piece $Y^{0}$ (which is not too informative). 
In the  escaping  case we will construct {\sl the principal nest}
 $$Y^0\supset V^1\supset V^2\supset\ldots  \eqno (3-0)$$
 in the following way. 
Let  $t$ be the first moment when $f^{tp}0\in Z^{(1)}_i$. Then let $V^0\ni 0$
be the pull-back of $Z^{(1)}_i$ along the ${\rm orb}_{tp}0$.
Further, let us define 
$V^{n+1}$ as the pull-back of $V^n$ corresponding to the first return of the
critical point 0 back to $\inter V^n$. Of course it may happen that 
the critical point never returns back to $\inter V^n$.
Then we stop, and the principal nest turns out to be finite.
This case is called  {\it combinatorially non-recurrent}.  
 If the critical point is {\sl recurrent} in the usual sense, that is 
 $\omega(0)\ni 0$, it is also combinatorially recurrent, 
 and the principal nest is infinite.

Let $l=l(n)$ be the first return time of the critical point back 
to $\inter V^{n-1}$. Then the map $g_n=f^{l(n)}: V^n\rightarrow V^{n-1}$
is two-to-one branched covering. Indeed, by the Markov property of the 
puzzle, $f^k V^n\cap \inter V^{n-1}=\emptyset$ for $k=1,\ldots, l-1$, 
so that the maps $f: V^k\rightarrow V^{k-1}$ are univalent for those $k$'s. 

 Let us call
 return to level $n-1$  {\sl central} if $g_n 0\in V^n$. In other words
$l(n)=l(n+1)$.
Let us say that a sequence $n, n+1,   n+N-1$ (with $N\geq 1$) 
 of levels (or corresponding 
puzzle pieces) of the principal nest form a
{\sl (central) cascade} if the returns to all levels $n, \ldots, n+N-2$ are
central, while the return to level   $n+N-1$ is non-central.
In this case 
$$g_{n+k}|V^{n+k}=g_{n+1}|V^{n+k},\;\; k=1,\ldots , N.$$
and $g_{n+1} 0\in V^{n+N-1}\backslash V^{n+N}$. 
Thus all maps $g_{n+1},\ldots, g_{n+N}$ are the same as quadratic-like. 
We call the number $N$ of levels
in the cascade its {\it length}. Note that a cascade of length 1 consists
of a single non-central level. Let us call the cascade {\it maximal} 
if the return to level $n-1$ is non-central.
 Clearly the whole principal nest is the 
union of disjoint maximal  cascades.
The number of such cascades is called the {\it height } $\chi(f)$  of $f$.
In other words, $\chi(f)$ is the number of different quadratic-like maps
among the $g_n$'s. (If $f$ is immediately renormalizable set $\chi(f)=-1$.)

\medskip\noindent {\bf Remark.} Given a quadratic polynomial $f: z\mapsto z^2+c$, 
the principal nest
determines a specific way to approximate $c$ by superattracting parameter values.
Namely, one should perturb $c$ in such a way that 
the critical point becomes  fixed under $g_n$, while the combinatorics on the preceeding
levels keeps unchanged. The number of points in this approximating sequence is equal
to the height $\chi(f)$. This resembles "internal addresses" of Dierk Schleicher [Sch]
but turns out to be different.

\medskip\noindent{\bf 3.2. Initial Markov partition.} 
Let $P_i$ be a finite or countable family of topological disc with
disjoint interiors,
and $g: \cup P_i\rightarrow \C$ be a map such that the restrictions $g|P_i$
  are branched coverings onto
their images. 
This map is called {\sl Markov} if $g P_i\supset P_j$ whenever
$\inter g P_i\cap \inter P_j\not=\emptyset$. Let us call it a
{\sl standard Markov map} if all restrictions $g|P_i$ are one-to-one onto
their images.

A  Markov map is called  
{ \sl Bernoulli (with range $D$)} if 
$gP_i\supset D\supset \cup P_j$ for all $i$.  Similarly we can define a
standard Bernoulli map.

We know that $f^p Y^{(1)}$
two-to-one covers $Y^{(1)}$ itself and the puzzle pieces 
$Z^{(1)}_i$ attached to $\alpha'$. If $f^p0\in Y^{(1)}$ (central return)
then  the pull back of $Y^{(1)}$ by this map is the critical piece
$Y^{(1+p)}$, while each $Z^{(1)}_i$ has two univalent pull-backs
$Z^{(1+p)}_j$ (we count them by $j$ in an arbitrary way) (see Figure 3). 

\midinsert
\centerline{\psfig{figure=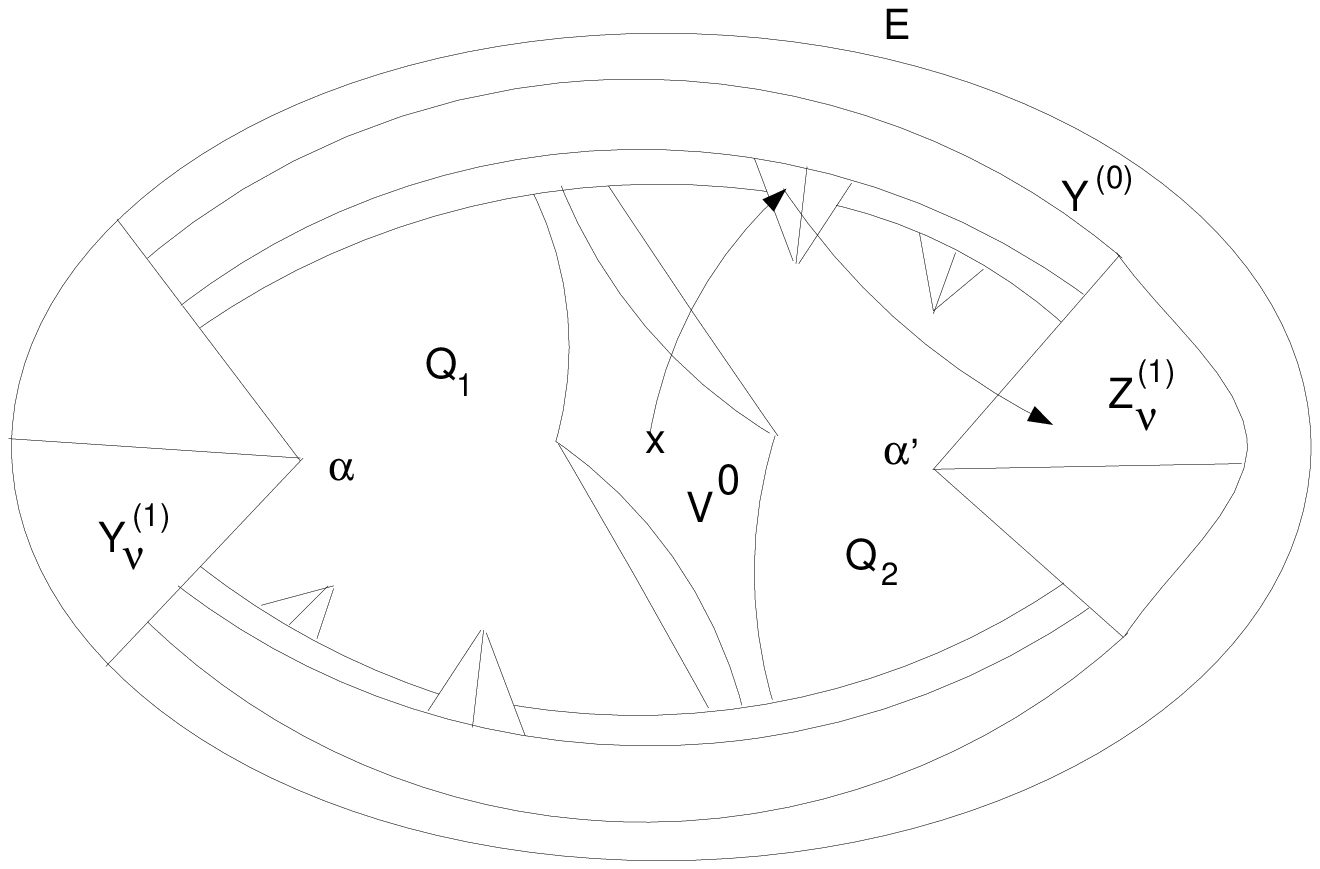,width=.8\hsize}}
\bigskip\centerline{Figure  3: Initial partition ($p=3$, $t=2$).}\medskip
\endinsert

   Now, $f^p Y^{(1+p)}$ two-to-one covers all these puzzle pieces.
If we again have a central return, that is $f^p 0\in Y^{(1+p)}$,
then the pull-back will give us one critical piece $Y^{(1+2p)}$,
and $4(p-1)$ off-critical $Z^{(1+2p)}_j$. 

Repeating this procedure $t$ times (where $f^{tp}0\in Z^{(1)}_i$),
we obtain the initial central nest
$$Y^{(1)}\supset Y^{(1+p)}\supset\ldots\supset Y^{(1+tp)},  \eqno (3-1)$$
and a family of non-critical puzzle pieces $Z^{(1+sp)}_j$,
$0\leq s\leq k-1$. Moreover
$$f^p 0\in Z^{(1+tp)}_{\nu}  \eqno (3-2)$$ for some $\nu$.

 Let us say that a set
$D$ is {\sl partitioned by pieces $W_i \rel F(f)$} if the
$\inter W_i$ are disjoint, and $D\cap J(f)=\cup W_i\cap J(f)$.

Thus we have partitioned $ Y^{(0)}\rel F(f)$ by
the pieces $Z^{(1+sp)}_i$, $0\leq s\leq t$,  and $Y^{(1+tp)}$.
 Let us look closer at this last piece. Its image under $f^p$ 
two-to-one covers all above puzzle pieces of depth $1+tp$. The pull-back of
$ Z^{(1+tp)}_{\nu}$ from (3-2) gives us exactly $V^0\ni 0$, the first
puzzle-piece in the principal nest. The pull-backs of other pieces
$Z^{(1+tp)}_j$ give some non-critical pieces $Z^{(1+(t+1)p)}_j $.
Finally, we have two univalent pull-backs $Q_1$ and $Q_2$ of $Y^{(1+tp)}$.
Altogether these pieces give a partition of $ Y^{(1+tp)}\rel F(f)$.

To understand how the critical point returns back to $V^0$ we need 
to partition further  $Q_1\cup Q_2$. To this end let us iterate
the standard Bernoulli map 
$f^p|Q_1\cup Q_2$ with range
$Q_1\cup Q_2\cup V^0\cup Z^{(1+(t+1)p)}_j$. So take a point
$z\in Q_1\cup Q_2$ and push it forward by iterated $f^p$ until it escapes
$Q_1\cup Q_2$ (or iterate forever if it does not escape). 
It can escape through the piece $V^0$ or
through a piece $Z^{(1+(t+1)p)}_j$. In any case pull the corresponding
piece back to this point. In such a way we will obtain the partition
$$Q_1\cup Q_2=\bigcup_{k>0}\bigcup_i X^{kp}_i 
  \bigcup_{k>t+1}\bigcup_j Z^{(1+kp)}_j\cup K
\;\rel F(f),$$
where $X^{k}_i$ denote the pull-backs of $V^0$ under $f^{kp}$
(for $k=0$ we have just one piece $X^0_0\equiv V^0$),
 $Z^{1+kp}_j$ denote the
pull-backs of  $Z^{(1+(t+1)p)}_j$ under $f^{(k-t-1)p}$,
 and $K$ denote the residual set of
non-escaping points.

Altogether we have constructed the {\sl initial Markov partition}:
$$Y^{(0)}\backslash K=V^0\bigcup_{k>0}\bigcup_i X^{kp}_i 
  \bigcup_{k>0}\bigcup_j Z^{(1+kp)}_j
 \rel F(f). \eqno (3-3)$$

It is convenient (to restrict the number of iterates in the following
considerations) to pass here to the following Markov map:
 $$G: \bigcup_{k,i} X^k_i\bigcup_{k, j} 
Z^{(1+kp)}_j\rightarrow \C \eqno(3-4).$$  
Observe that for any $j$ there is an $i$ such that
$f^{sp} Z^{(1+sp)}_j$ univalently  covers  $Z^{(1)}_i$. 
Moreover $f^{p-i} Z^{(1)}_i$ univalently covers 
$Y^{(0)}$. Let us set $G|Z^{1+sp}_j=f^{ps+(p-i)}$. 
 The image of each piece 
$Z^{(1+sp)}_j$ under this map univalently covers 
$Y^{(0)}$.
Similarly let us set $G|V^0=f^{tp+p-i}$, where $f^{tp }0\in Z^{(1)}_i$.
The image of this piece two-to-one covers $Y^{(0)}$. 
Finally  $G|W^k_i=f^{kp}$ for $k>0$. These pieces are univalently mapped 
onto $V^0$.

\medskip\noindent{\bf 3.3. Non-degenerate annulus.}  Yoccoz has shown that
if $f$ is non-renormalizable then in the nest 
$Y^{(1)}\supset Y^{(2)}\ldots$ there is a non-degenerate annulus
$Y^{(n)}\backslash Y^{(n+1)}$. However the modulus of this annulus is not under
control. We will construct a different non-degenerate annulus whose modulus
we will be able to control. 

Let $A^n=V^n\backslash V^{n+1}$.

\proclaim Proposition 3.1. Let $f$ be a DH  quadratic-like map 
 which is not immediately DH renormalizable. Then all annuli
 $A^n$ are non-degenerate.  

\medskip\noindent{\bf Proof.} Observe first that $V^0$ is strictly inside
 $Y^{(0)}$, that is, the annulus $Y^{(0)}\backslash V^0$ is non-degenerate.
 Indeed, $V^0$ is the pull-back of $Z^{(1)}_i$ which is strictly inside
 $Y^{(0)}$. As the iterates of $\partial Y^{(0)}$ stay outside 
 $\inter Y^{(0)}$, $V^0$ may not touch $\partial Y^{(0)}$.

For the same reason all other pieces $Z^{(1+kp)}_j$ and $X^k_i$
of the initial Markov partition  are strictly inside $Y^{(0)}$ as well. 

Let us consider the orbit of the critical point
  0 under iterates of the map $G$ from (3-4) until it returns back to
$V^0$.
It first goes through the 
  $Z$-pieces of the initial Markov partition, then at some moment $l\geq 1$
  lands at a $X^s_i$, and the latest at the next moment returns back to $V^0$.  
 
Since the map $G: V^0\cup Z^{(1+kp)}_j\rightarrow \C$ is Bernoulli
with range $Y^{(0)}$, there is a piece $P\subset V^0$, such that
$G^l P$ two-to-one covers $Y^{(0)}$.  Clearly $V^1$ is the pull-back
of $X^s_i$  by $G^l$. Since $X^s_i\subset\subset Y^{(0)}$, we conclude that
$V^1\subset\subset V^0$.

Now it is easy to see that all annuli $A^n$ are non-degenerate as well.
Indeed, it follows that the orbit of $\partial V^1$ stays away from $V^1$.
Hence $V^2$ cannot touch $\partial V^1$, for otherwise there would be
a point on $\partial V^1$ which returns back to $V^1$. So $A^2$ is
no-degenerate. Now we can proceed inductively.   \QED

\medskip\noindent
\noindent {\bf 3.4. DH renormalization and central cascades.}
In this section quadratic-like maps and renormalization are understood in
the sense of Douady and Hubbard.

\proclaim Proposition 3.2. A quadratic-like map is renormalizable
   if and only if it is either immediately renormalizable, or
  the principal nest $V^0\supset V^1\supset\ldots$
  ends with an infinite cascade of central returns.
   Thus the height $\chi(f)$ is finite if and only if $f$ is renormalizable.

\noindent{\bf Proof.} Let the principle nest ends with an infinite
central cascade $V^{m-1}\supset V^m\supset\ldots$ Then $l_m=l_{m+1}=\ldots=l$,
and $g\equiv g_m=f^l|V^m$. Hence $\cap V^k$ consists of all points which never
escape $ V^m$ under iterates of $g$, that is, $\cap V^k=K(g)$.
Since $0\in \cap V^k$, $K(g)$ is connected, and
 $g: V^m\rightarrow V^{m-1}$ is quadratic-like (take into account that
$V^m\subset\subset V^{m-1}$ by Proposition 3.1).

Take now the non-dividing fixed point $b$ of $g$. Let us show that
$b$ is dividing for the big Julia set $K(f)$. To this end let us consider
the configuration of the {\it full} external rays whose segments
  bound $V\equiv V^m$. They
divide the plane into the central component containing $V$, and the
family $\SS=\{S_i\}_{i=1}^t$ of sectors. As $g$ maps $V$ onto a bigger piece,
$g(\partial S_i)\subset  S_{\sigma(i)}$ for any sector $S_i\in \SS$.

Let $i_1,\ldots, i_r$ be a set of indeces which are cyclically permuted
by the map $\sigma: \{1,\ldots, t\}\rightarrow \{1,\ldots, t\}$.
Then there is a sequence of sectors $S_{i_1}\equiv T_1,\ldots, T_r,\ldots$
with the following properties:

\smallskip\noindent$\bullet$ Each $T_k$ is bounded by two external rays,
 and does not contain 0;

\smallskip\noindent$\bullet$ The angle of each $T_k$ at infinity is
  smaller than $\pi$;

\smallskip\noindent$\bullet$ $T_k\supset S_{i_l}$ with $l\equiv k\; \mod r$;

\smallskip\noindent$\bullet$ $g( \partial T_{k+1})=\partial 
  T_{k}$, $k=1,2, \ldots$

\smallskip\noindent$\bullet$ $T_k\subset T_{k+r}$, while
 $T_{k+1}, \ldots, T_{k+r}$
  are pairwise disjoint.

\smallskip Then for any $k$ we have a nest 
of sectors, $T_k\subset T_{k+r}\subset T_{k+2r}\ldots$, converging
to a sector $\tilde T_k$.  These sectors also satisfy all the above properties,
and moreover $\tilde T_k=\tilde T_{k+r}$.
 Hence $g$ cyclically permut the vertices of these sectors, which therefore
form a cycle of period $r$. 
Moreover, the points of this cycle are clearly dividing.

Let us finally note that actually $r$ must be equal to 1, so that this
cycle actually coincides with the fixed point $b$. Indeed, the opposite
situation would contradict Proposition 2.2 (ii).

So we have a dividing periodic point $b$.
Let $\Omega'\subset\Omega$
be the  corresponding domains constructed in \S 2.5. 
Then $\Omega'\supset K(g)$ and hence
$\Omega'\ni 0$. It follows that $g: \Omega'\rightarrow \Omega$ is a double
covering. Moreover, $g^n 0\in K(g)\subset\Omega'$, $n=0,1,\ldots$. Thus
$f$ is renormalizable.

Assume now that $f$ is renormalizable. Let $a$ be the corresponding
periodic point, and $\Omega'$ be the region bounded by rays landing
at $a$ and arcs of an equipotential, such that
$Rf=f^l: \Omega'\rightarrow\Omega$ is a double covering.

Then the fixed point $\alpha$ may not lie in $\inter \Omega'$, for otherwise
$\inter f\Omega'$ would intersect $\inter\Omega'$. Hence $\alpha$ does not
cut the filled Julia set $K(Rf)$.
But then the preimages of $\alpha$ don't cut $K(Rf)$ either.
Hence given a puzzle-piece $Y^{(n)}_i$, either  $K(Rf)$ is contained
in $Y^{(n)}_i$, or $K(Rf)\cap \inter Y^{(n)}_i=\emptyset$. In particular
$V^m\supset K(Rf)$. But then $f^l 0\in V^m$ for all $m$, so that
the first return times to $V^m$ are uniformly bounded. But then 
this nest must end up with a central cascade.  \QED

The above discussion shows that there is a well-defined
 first  renormalization $Rf$
with the biggest Julia set, and it can be constructed
 in the following way. If $f$ is immediately renormalizable, then
$Rf$ is  obtained by thickening $Y^{(1)}\rightarrow Y^{(0)}$. Otherwise
the principal nest ends up with the infinite central cascade
$V^{m-1}\supset V^{m}\supset\ldots$. Then 
$Rf=g_m: V^m\rightarrow V^{m-1}$. 

The internal class $c(Rf)$ of the first renormalization 
 belongs to a maximal copy $M_0$ of the
Mandelbrot set, that is, a copy which is not contained in any other
one  except the whole $M$.

\medskip\noindent{\bf 3.5. Return maps and Koebe space.} 
Let $f$ be a DH quadratic-like map, and let $V\in \MM(f)$ be a puzzle piece.

\proclaim Lemma 3.3. 
 Let $z$ be a point whose orbit passes
  through  $\inter V$. Let $l$ be the first positive moment of time
  for which $f^l z\in \inter V$. Let $U\ni z$ be the puzzle piece mapped
   onto $V$ by $f^l$. Then $f^l: U\rightarrow V$ is either a
   univalent map or two-to-one branched covering depending on whether
  $ U$ is non-critical or otherwise.

\noindent{\bf Proof.}  Let $U_k=f^k U,\; k=0,1,\ldots l$. Since
$f^k z\not\in \inter V$ for $0<k<l$, by the Markov property of the puzzle,
$U_k\cap \inter V=\emptyset$ for those $k$'s. 
 Hence $f: U_k\rightarrow U_{k+1}$
is univalent for $k=1,\ldots, l-1$, and the conclusion follows.  \QED

Let $z\in \inter V$ be a point which returns back to $\inter V$, and
let $l>0$ be the first return time. 
Then there is a puzzle piece $V(z)\subset V$
containing $z$ such that $f^l V(z)=V$.
It follows that the first return map $A_V f$ to $\inter V$  
is defined on the union   of disjoint open
puzzle pieces $ \inter V_i$. Moreover, if
$$f^m \partial V\cap V=\emptyset,\; m=1,2, \ldots  \eqno (3-4)$$ 
then it is easy to see that 
the closed pieces $V_i$ are pairwise disjoint and are contained in $\inter V$.
Indeed, otherwise there would be a boundary  point 
$\zeta\in \partial V$ whose orbit would return back to $V$, despite (3-4).

 Somewhat non-rigorously,  we will call the map
$$A_V f: \bigcup V_i\rightarrow V \eqno (3-5)$$
 the first return map to $V$.
(Note: it perhaps may happen that a point $z\in \partial V$ returns back to $V$
but does not belong to $\cup V_i$.) 
Let $V_0$ denote the critical (``central") puzzle piece (provided
the critical point returns back to $V$). Now Lemma 3.3 immediately 
implies:

\proclaim Lemma 3.4. The first return map  $f_V$  univalently maps
 all non-critical pieces $V_i$ onto $V$, and maps the critical piece
 $V_0$ onto $V$  as a double branched covering.

Thus the first return map $A_V f$ is Bernoulli, and is  standard Bernoulli
on $\bigcup_{i\not=0} V_i$.

Let us now improve Lemma 3.3. 

\proclaim Lemma 3.5. Let $z$ be a point whose orbit passes through the
central domain $\inter V_0$ of the first retun map (3-5), 
and $l\geq 0$ be the first  moment when $f^l z\in V_0$. Then  there is
a puzzle piece $\Omega\ni z$ mapped univalently by $f^l$  onto $V$.

\noindent{\bf Proof.} Let $s$ be the first moment
when $f^s z\in V$. Then $f^l z= f_V^k(f^s z)$ for some $k\geq 0$.
Moreover, $(A_Vf)^r (f^s z)\not\in \inter V_0$ for $r<k$. 

Since the
return map is standard Bernoulli outside of the central piece,
there is a piece $X\subset V$ containing $f^s z$ which is univalently 
mapped by $(A_V f)^k$ onto $V$. On the  
other hand, by Lemma 3.3 there is a domain $D\ni z$ which is univalently mapped
by $f^s$ onto $V$. Hence the domain $D\cap f^{-s} X$ is univalently mapped
by $f^l$ onto $V$.             \QED

Let us  now consider the principal nest of $f$:
$$Y^0\supset V^0\supset V^1\supset\ldots.$$
Let $$g_n: \cup V^n_i\rightarrow V^{n-1}  \eqno  (3-6)$$ 
be the first return map
to $V^{n-1}$, where $V^n_0\equiv V^n\ni 0$. We will call them the 
{\it principal return maps}. We will also let $g_0\equiv f$.

Let $\Phi: z\mapsto z^2$ be the quadratic map. 
Since  $f: U'\rightarrow U$ is a double covering with the critical point
at  0, it can be decomposed
as $h\circ\Phi$ where $h$ is a univalent map with range $U$.

\proclaim Corollary 3.6. For $n\geq 1$ the pieces $V^n_i$ are pairwise
  disjoint, and the annuli $V^{n-1}\backslash V^n_i$ are non-degenerate.
For $n\geq 2$ the map $g_n|V^n_i$ can be decomposed as $h_{n,i}\circ \Phi$ where
$h_{n,i}$ is a univalent map with range $V^{n-2}$.

\noindent{\bf Proof.} As
$V^n\subset \inter V^{n-1} $ (Proposition 3.1), 
and $$f^m(\partial  V^n)\cap \inter V^{n-1}=\emptyset,$$
condition (3-4) is satisfied here, and the first statement follows.   

  Take a piece $V^n_i$, and let $g_n|V^n_i=f^l$.  
By Lemma 3.5 the map $f^{l-1}: (f V^n_i)\rightarrow V^{n-1}$
 can be extended to
a univalent map with range $V^{n-2}$, and the second statement  also follows.
\QED

\medskip\noindent{\bf 3.6. Bernoulli scheme associated to a central cascade.}
In the case of a central cascade we need a more precise analysis of the 
Koebe space. Let us consider a central cascade $\CC\equiv \CC^{m+N}$:

$$V^m\supset V^{m+1}\supset\ldots\supset V^{m+N-1}\supset V^{m+N},\eqno (3-7)$$
where $g_{m+1}0\in V^{m+N-1}\backslash V^{m+N}$. Set $g= g_{m+1}|V^{m+1}$. Then
$g: V^{k}\rightarrow V^{k-1}$ is a double branched covering, 
$k=m+1,\ldots m+N$.

Let us consider the first return map 
$g_{m+1}: V^{m+1}_i\rightarrow V^m$ as in (3-6).  Let us pull the pieces
$V^{m+1}_i$ back to the annuli $A^k=V^{k-1}\backslash V^k$ 
by iterates of $g$, $k=m+N,\ldots, m+1$.
We obtain a family $\WW(\CC)$ of pieces $W^k_j$. By construction,
$W^k_j\subset A^k$ and $g^{k-m+1}$ univalently maps $W^k_j$ onto 
some $V^{m+1}_i\equiv W^{m+1}_i$.

Let us define a standard Bernoulli map $G\equiv G^{m+N}$:
$$G: \bigcup_{\WW(\CC)} W^k_j\rightarrow  V^m  \eqno(3-8)$$
as follows: $G| W^k_j= f_{m+1}\circ g^{k-m+1}$ (see Figure 4).

\medskip\noindent{\bf Remark.} This Bernoulli map is  similar to the 
 initial Bernoulli map constructed in \S 3.5. Actually in the initial 
construction we deal with the central cascade (3-1)
with degenerate annuli.

\midinsert
\centerline{\psfig{figure=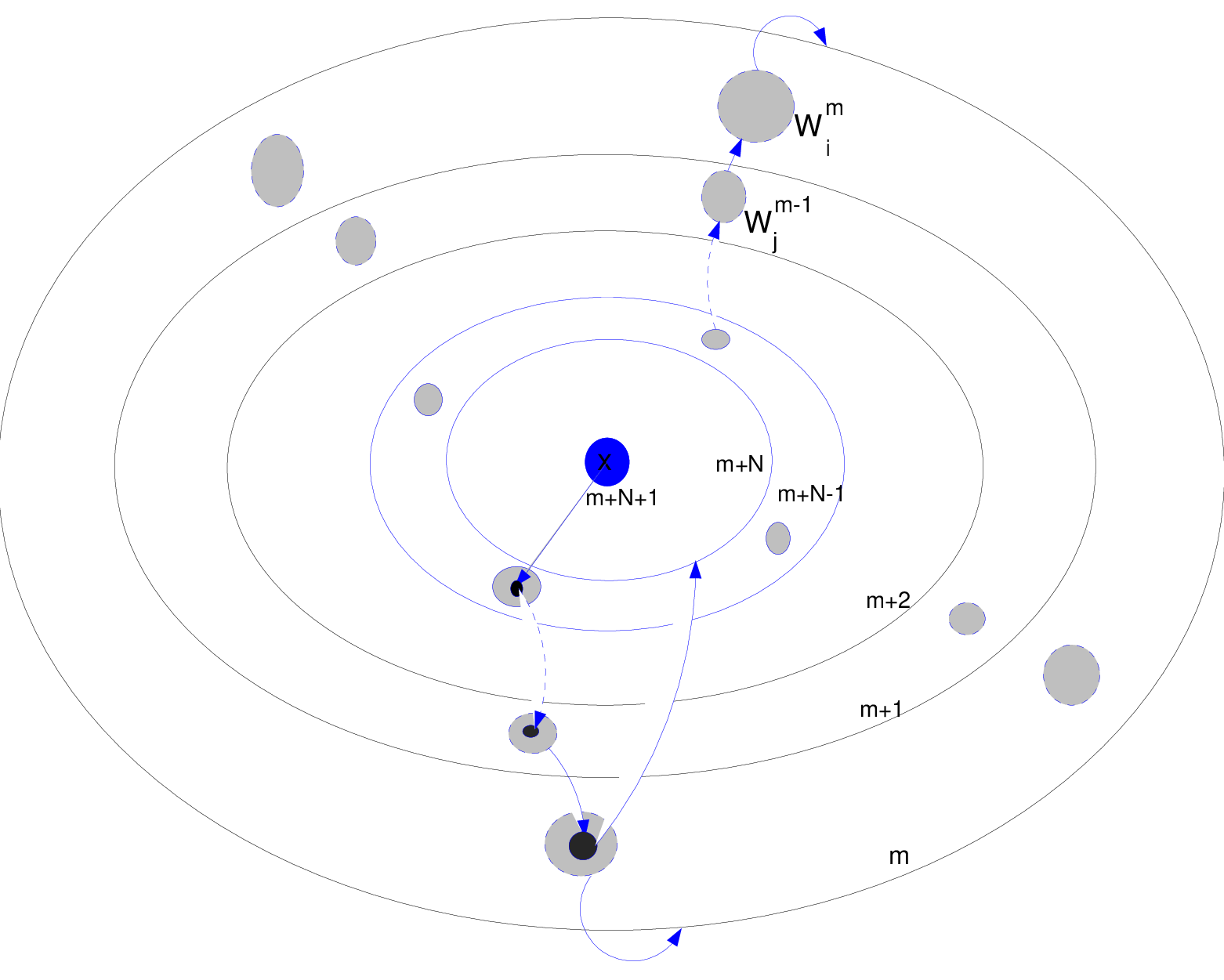,width=.8\hsize}}
\bigskip\centerline{Figure 4: Bernoulli scheme associated to a central cascade.}
\medskip\endinsert

\proclaim Lemma 3.7. Let us consider the central cascade (3-7).
  Let $z$ be a point whose orbit passes through
  $\inter V^{m+N}$, and $l$ be the first moment for which 
  $f^l z\in V^{m+N}$. Then there is a piece $\Omega\ni z$ which is
  univalently mapped by $f^l$ onto $V^m$. 

\noindent{\bf Proof.}  Let $s$ be the first moment
for which  $f^s z\in V^m$. Then $f^l z= G^k(f^s z)$ where $G$ is the
Bernoulli map (3-8). 
Now repeat the  argument of Lemma 3.5  just using $G$ 
instead of the first return map.  \QED

\proclaim Corollary 3.8. Let us consider the central cascade (3-7).
Then the map $g_{m+N+1}: V^{m+N+1}\rightarrow V^{m+N}$ can
be represented
 as $h_{m+N+1}(z^2)$ where $h_{m+N+1}$ is a univalent
map with range $V^{m}$.

\noindent{\bf Proof.}  Repeat the argument for Corollary 3.6
using  Lemma 3.7 instead of 3.5.  \QED

\medskip\noindent{\bf 3.7. Generalized polynomial-like maps and 
      renormalization.}
Let $\{U_i\}$ be a finite or countable 
family of topological discs with disjoint interiors
strictly contained in a topological disk $U$. We call a map
$g: \cup U_i\rightarrow U$ a {\sl (generalized ) polynomial-like map}
if $g:U_i\rightarrow U$ is a branched covering of 
   finite degree which is univalent on all but finitely many $U_i$.


Then we can define the filled Julia set $K(g)$ as the set of all
non-escaping points, and the Julia set $J(g)$ as its boundary.
The DH polynomial-like maps correspond to the case of a single disk $U_0$.

Let us say that a polynomial-like map $g$ is {\it of finite type} if its domain
consists of finitely many disks $U_i$.  

\proclaim Generalized Straightening Theorem.
 Any generalized polynomial-like map
of finite type is qc conjugate to a polynomial with the same number of
non-escaping critical points. 

\proclaim Lemma 3.9. A generalized polynomial-like map with 
non-escaping critical point has a connected
Julia set/(filled Julia set) if and only if it is DH polynomial-like.

Let us call a (generalized) polynomial-like map a {\sl (generalized)
quadratic-like} if it has a single (and non-degenerate) critical point.
 In such a case
we will always assume that 0 is the critical point, and count the
discs $U_i$ in such a way that $U_0\ni 0$. In what follows we will deal
exclusively with quadratic-like maps, namely with the 
 principal sequence $g_n$ of the first return maps (3-6).  

Given a $V_j^{n+1}$, let  $l$ be its first
return time back to $V^n$, that is, $g_{n+1}|V^{n+1}_j= g^l_n$. 
Then $$g_n^k V_j^{n+1}\subset V_{i(k)}^n,\; k=0,1,\ldots, l,$$
with $i(0)=i(s)=0$. Moreover, 
$g_n^k V_j\subset\subset V^n_{i(k)}$ for $k<s$. The sequence 
$0=i(0), i(1),\ldots, i(s)=0$ is called the {\sl itinerary }
of $V_j^{n+1}$ through the domains of previous level.

Philosophically the dynamical renormalization is the first return map
to an appropriate piece of the phase space considered
up to rescaling. 
In our setting
let us define the $n$-fold { \it generalized renormalization } $T^n f$
of $f$ as the first return map $g_n$ restricted to  the union of puzzle
pieces $V^n_i$ meeting the critical set $\omega(0)$, and considered up to
rescaling.  In the most interesting situations these maps are of finite type
(compare [L3]):

\proclaim Lemma 3.10. If $f$ is  a DH renormalizable quadratic-like map,
 then all maps $T^n f$ are of finite type.

{\bf Proof.} In the tail of the principal nest the maps $g_n$ are DH 
quadratic-like, and their domains consist just of one component.
So we should take care only of the initial piece of the cascade.

Let us take the  renormalization 
$Rf=f^l: V^{t+1}\rightarrow V^t$ with $t\geq n$. 
Since 0 is non-escaping under iterates of
$Rf$, we have the following property: the first return time of any point
$f^k 0$ back to $V^{t+1}$ is at most $l$. All the more, the return time to 
any bigger domain $V^n\supset V^t$ is bounded by $l$.
Hence the components of $f^{-t}V^n$, $t=0,1,\ldots, l-1$, cover
the whole postcritical set. For sure there is only finitely many such
components. But the domain of $T^n f$ is the union of some of
these components (which are inside $V^n$).   \QED

\medskip\noindent{\bf  3.8. Return graph.} Let $\II^n$ be the family of puzzle
  pieces $V^n_i$ intersecting $\omega(0)$, that is, the pieces in the domain of
  the generalized renormalization 
$$T^n f: \bigcup_{\II^n} V^n_i\rightarrow V^{n-1}.$$
Let us consider a graded graph $\Upsilon_f$ whose vertices of level $n$
are the pieces $V^n_j\in \II^n$, $n=0,1,\ldots$.
 Let us take a vertex  $V^n_j\in \II^n$, and
let $i(1),\ldots, i(t)=0$ be its itinerary through the pieces of the
previous level under the iterates of $g_n$ (see the previous section).
  Then for $n\geq 1$ join $V^{n+1}_j$
with $V^{n}_i$ by $k$ edges, provided the symbol $i$ appears in the above 
itinerary $k$ times. This means that the piece $V^{n+1}_j$ under iterates
 of $g_n$ passes through
$V^{n}_i$ $k$ times before the first return back to $V^{n-1}$.
As to the top level, let us join each $V^1_j$ with $V^0$ by the number
of edges equal to the first return time of $V^1_j$ back to $V^0$
under iterates of $f=g_0$.

\smallskip{\bf Remark.} A similar graph in the  real one-dimensional setting
 was introduced by Marco Martens [Ma]. The above graph is not 
exactly the same as Marten's graph, as the latter is related to the 
iterates of $f$ itself rather than the renormalized maps.\smallskip

Note that for any vertex $V^{n+1}_j$ with $n>0$
 there is exactly one edge joining it
to the critical vertex $V^n_0$ of the previous level.  Note also that 
by Lemma 3.10 in the DH renormalizable case the
number of vertices on a given level is finite. In any case there are clearly
only finitely many edges leading from a $V^{n+1}_j$ to the previous level.

By a  {\it path} in the graph $\Upsilon_f$ we mean a connected sequence
of edges such that no two of them join the same two levels, up to reversing the
order of the sequence. So we don't endow
the paths with orientation, and can go along them either strictly upwards
or strictly downwards.

Diverse combinatorial data can be easily read off this graph.
For example, given $n\geq m$, the number  
 of paths joining $V^{n+1}_j$ to $V^m_i$ is equal to the number of times  which
the $g_m$-orbit of $V^{n+1}_j$ passes through $V^m_i$ before the first return
back to $V^n$. 
Hence the return time of $V^{n+1}_j$ back to $V^{n}$ under
iterates of $g_m$ is equal to  the  total number
of paths in $\Upsilon_f$  leading from $V^{n+1}_j$  up to level
$m$. Let $g_m-time(V^{n+1}_j)$ denote this return time.  

Let now the map $f$ be DH renormalizable, 
and $t$ be a renormalization level in the
principal nest, that is, 
$g_{t+1}: V^{t+1}\rightarrow V^t$ is a quadratic-like map with non-escaping
critical point. 
 Then there is a single vertex $V^{t+1}$ at level $t+1$,
and below it the return graph  is just 
the ``vertical path" through the critical vertices. 
Let  $per(f)\equiv g_0-time(V^{t+1})$ 
denote {\it the renormalization period}, that is, the
return time of $V^{t+1}$ back to $V^{t}$ under iterates of the original map 
$f=g_0$.   


By the above discussion, $per(f)$ is equal to the 
total number of paths in the graph  $\Upsilon_f$
(joining the top vertex $V^0$ to the bottom vertex $V^{t+1}$).
It follows that the $per(f)$ is bounded if and only if the DH-level $t$ is
bounded, and all return times $g_m-time(V^{m+1}_i)$ are bounded for
$1\leq m\leq t$ and any $i$. For example, the ``if" statement means:
 If $t\leq T$ and $g_m-time(V^{m+1}_i)\leq R$, then
$per(f)\leq P(T, R)$. Indeed, in this case
the total number of paths in the graph is bounded by $R^T$.

Let us now take a closer look at the central cascades.
Let us consider  a central cascade (3-7) and a non-pre-critical
piece $V^{m+N+1}_j$, that is $g_{m+N} V^{m+N+1}_j\not=V^{m+N}_0$ . 
Then  let us denote by 
$G_{m+N}-time(V^{m+N+1}_j)$ the first return time 
of $g_{m+N} V^{m+N+1}_i$ 
back to $V^{m+N}$ under the iterates of the map $G_{m+N}$ from (3-8).

This time can be expressed in terms of the following {\it reduced graph}
$\Upsilon_f^r$. This graph is obtained from the $\Upsilon_f$ by removing
from all central cascades (3-7)
the edges leading from the non-critical pieces $V^{k+1}_j$, $j\not= 0$,
to  the critical piece $V^k_0$, $k=m+N-1,\ldots , m$.  
 Then  $G_{m+N}-time(V^{m+N+1}_j)$ is equal to the number of paths 
in the reduced graph leading from $V^{m+N+1}_j$ up to level $m$. 

Let us define the {\it reduced period }  $per_r (f)$ as the total number
of paths in the reduced graph. This means that we don't count the moments
of time when the orbit goes through the intermediate  levels of central
cascades. Clearly {\it the reduced period of $f$ is bounded if and only if
the height  $\chi(f)$ and all $G-times$ are bounded.}
 
Let us finally define one more combinatorial notion, the rank. 
Let $V^{n+1}_j$ and $V^n_i$ be two $\Upsilon_f$-adjacent  puzzle pieces,
and let $\gamma$ denote an edge  joining them. 
Let $D^n\supset V^{n}_i$
 be  a puzzle piece  contained in $V^{n-1}$. 
Consider the first moment $t$ for which
 $g_n^t V^{n+1}_j\subset V^n_i$.
The piece $D^{n+1}\supset V^{n+1}_j$ in  $V^n$ 
  such that $g_n^t D^{n+1}=D^n$ will be  called the {\it pull-back of
$D^{n}$ along the  edge} $\gamma$. 
 More generally, 
let us define the {\it pull-back of $D^n$ along a path} $\gamma$
leading from $V^n_i\subset D^n$ downwards by composing the pull-backs 
along the edges.

Given a  piece $V^n_k$, let $rank V^n_k$ be the
number of non-central levels in the shortest path leading from $V^n_k$
 down to a critical piece. For a piece $D^n\subset V^{n-1}$ 
 as above, let $rank(D^n)$ denote the minimum of ranks of puzzle pieces
$V^n_k$ contained in $D^n$.
Note that $rank(D^n)=0$ iff $D^n$ 
is critical.

By a path  through $D^n$ we will mean a path through a piece $V^n_i\subset D^n$.

\proclaim Lemma 3.11.  Let $\gamma$ be the shortest path leading
from $D^n$ down to a central piece $V^{n+t}_0$, and let $D^{n+t}$ be the
pull-back of $D^n$ along this path. Then 
$$V^{n+t}\subset D^{n+t}\subset V^{n+t-1},$$
and the map $D^{n+t}\rightarrow D^n$ is a double branched covering.

\noindent{\bf Proof.} Follows from the definitions.  \QED

\medskip\noindent{\bf 3.9. Full principal nest.} Let $f$ be a
DH renormalizable, but not immediately, quadratic-like map.
Then its principal nest 
$$Y^{(0,0)}\supset V^{0,0}\supset V^{0,1}\supset\ldots
   \supset V^{0,t(0)}\supset V^{0,t(0)+1}
\supset\ldots$$
ends up with an infinite cascade of central returns (we mark this nest with
two labels for the reason which will become clear in a moment). Let us select
a level $t(0)$ of this cascade, so that the return map
$Rf=g_{0,t(0)+1}:  V^{0,t(0)+1}\rightarrow  V^{0,t(0)}$
is  DH quadratic-like. We will call such a level DH.
(The particular choice of DH renormalizable levels in what follows
 will depend on the geometry).  

If $Rf$ is not immediately DH renormalizable,
 let us cut the puzzle piece $ V^{0,t(0)+1}$ by the external rays, 
and construct its short principal nest:
$$Y^{(1,0)}\supset V^{1,0}\supset V^{1,1}\supset\ldots\supset 
V^{1,t(1)}\supset 
V^{1,t(1)+1}\supset\ldots$$
If $Rf$ is DH renormalizable, then this nest also
 ends up with an infinite central cascade.
Then select a DH level $t(1)+1$, and pass to the next short nest.

If $f$ is infinitely DH renormalizable but  non of the renormalizations are 
immediate, then in such a way we construct the {\sl full principal nest}
$$\eqalign
{Y^{(0,0)}\supset 
V^{0,0}\supset V^{0,1}\supset\ldots\supset V^{0,t(0)}\supset V^{0,t(0)+1}
\supset\cr
Y^{(1,0)}\supset
V^{1,0}\supset V^{1,1}\supset\ldots\supset V^{1,t(1)}\supset 
V^{1,t(1)+1}\supset\ldots\cr
Y^{(m,0)}\supset
V^{m,0}\supset V^{m,1}\supset\ldots\supset V^{m,t(m)}\supset 
V^{m,t(m)+1}\supset\ldots\cr}. \eqno (3-9)$$
Here $Y^{(m,0)}$ is the first critical Yoccoz puzzle piece for  the
$m$-fold DH renormalization $R^m f$, while
the pieces $V^{m,n}$ form the  corresponding short principal nest.
Moreover,  for $m>1$, $Y^{(m,0)}$ is obtained by cutting 
$V^{m-1,t(m-1)+1}$ with the 
external  rays of 
 $R^m f: V^{m-1, t(m-1)+1}\rightarrow V^{m-1,t(m-1)}$.

The annuli $A^{m,n}=V^{m,n-1}\backslash V^{m,n}$ will be called the
{\it principal annuli}.

\medskip\noindent{\bf 3.10. Big type: special families of Mandelbrot copies.}
Assume that we associated to any quadratic-like map its
``combinatorial type" $\tau(f)$, which depends only on the hybrid class
$c(f)$ and is 
constant over any maximal copy of the Mandelbrot set.
(Keep in mind the height function $\chi(f)$ or the period $per(f)$.)  
Thus we can use the notation $\tau(M')$.

Let $\SS\subset\MM$ be a family of maximal copies of the Mandelbrot set.
 Let us call it {\sl $\tau$-special}
if it satisfies the following property: For any truncated secondary
limb $L$ there is a  $\tau_L$ such that  $\MM$ contains
all copies $M'\subset L$ 
of the Mandelbrot set with $\tau(M')\geq\tau_L$.

Let $f$ be an infinitely DH renormalizable quadratic-like map. Let us 
say that it is of $\SS$-type
if all internal 
classes $c(R^n f)$ belong to copies  $M'$
from  $\SS$.

\bigskip

\centerline{\bf \S 4. Initial geometry.}\bigskip

\medskip\noindent{\bf 4.1.  Geometry of rays.} 
Let $L\equiv L_b$ be a limb of the Mandelbrot set with root at $b$.
By $L^{tr}=L^{tr}_b$ 
we denote a truncated  limb, that is,
$L$ with a neighborhood of the root removed.

We say that a pre-periodic point is {\sl repelling} if the corresponding 
periodic point is.  
Given a repelling periodic or pre-periodic point $a$  of $P_c$,  
let $\RR(a)$ be the union of the segments of rays landing at $a$ up to
the equipotential level 1 together with point $a$.
We will use the Hausdorff topology 
on the space of configurations of curves. We will often suppress the label
 $c$ unless it may cause confusion. Let $B_c$ denote the B\"{o}ttcher 
function (2-1) of $P_c$.

\proclaim Lemma 4.1 (see [GM]). Let $a_c$ be a periodic point of $P_c$
continuously depending on $c\in L_b$
which stays repelling on the unrooted limb $L_b\backslash \{b\}$.
Then the rays configuration $\RR(a_c)$ depends continuously on 
$c\in L_b\backslash \{b\}$. 

\noindent{\bf Proof.}  By [GM] the external angles of the rays landing at
$a_c$ are the same through the limb. So $\RR(a_c)$ is the $B^{-1}$-image
of the fixed configuration of segments in the annulus $\{z: 1<|z|<e\}$
with the point $a_c$ added.

Let us go back to the construction of the 
external  rays in \S 2.2. It is easy to see from the explicit formulas for 
the Rimann map and the linearizing coordinate that the map (2-3)
$ h_i^{-1}: \H\rightarrow \DD_i$ depends continuously 
(in the compact-open topology)  on $c$ ranging over the unrooted limb.
 Let $I=\{iy: 1\leq y\leq 2^p\}$ be a fundamental
segment of the vertical geodesic in $\H$.  It follows that the curves
$\hat J=h_i^{-1} I$  depend smoothly on $c$. Hence the curves
 $$\hat R_i=\bigcup_{n=0}^{\infty} \lambda^{-n} \hat J$$
depend continuously on $c$ in the Hausdorff topology. Then the rays
$R_i=\psi(\hat R_i)$ depend continuously on $c$ as well.
   \QED

Given  a configuration $\CC_0$ of finitely many parametrized curves and point
in   ${\bf C}$, let us consider the space Teich$(\CC_0)$ of all 
configurations qc equivalent to $\CC_0$ modulo conformal equivalence.
 There is a natural Teicm\"{u}ller distance on this space:
$$\dist_T(\CC_1, \CC_2)=\inf \log K_h,$$
where $h$ runs over all qc equivalences between $C_1$ and $C_2$. 

We say that  configurations  of some family 
 have  {\sl bounded  geometry} if they stay  bounded Teichm\"uller 
distance from a reference configuration $\CC_0$ whose curves are smooth
and  intersect transversally.

Let us consider a family of two topological nested disks $D_1\subset D_2$ with 
$\Gamma_i=\partial D_i$, and $A=D_2\backslash D_1$.
The statement that $\mod(A)>\epsilon$ with an $\epsilon>0$
 uniform over the family
will be freely expressed in the following ways:
``The annulus $A$ has a definite modulus", or ``$D_1$ is well inside $D_2$",
or ``There is a definite space in between $\Gamma_1$ and $\Gamma_2$."

A statement like ``If $f$ has a definite modulus then  a certain 
quantity  is bounded"
means: ``If $\mod(f)>\epsilon$ then there is a bound on that quantity
depending only on $\epsilon$".  

\proclaim Lemma 4.2. Under the circumstances of Lemma 4.1 
 the configuration $\RR(a_c)$ has  a bounded geometry
  when $c$ ranges over the truncated limb  $L_b^{tr}$ of the Mandelbrot set.


\noindent{\bf Proof.} Let us consider the  configuration $\II$
of intervals obtained by rotating  the interval $[0,1]$ by the cyclic
group $Z_p$ of order $p$. We will show that  the configurations
$\RR(\alpha_c)$ of rays stay bounded   Teichm\"{u}ller distance
from $\II$.

Let us use the notations from the proof of the previous lemma.
Since the configuration $\hat\RR$ of rays
 continuously depends on $c$ ranging over the truncated limb $L^{tr}$,
 we can find a smooth Jordan curve $\gamma$ enclosing 0 and
smoothly depending  on $c$, 
which intersects each infinite ray $\hat\RR_i^{\infty}$ at only one point,
namely the endpoint of $\hat\RR_i$,  and this 
intersection is transversal. The curves $\gamma_{-N}=\lambda^{-N}\gamma$ 
clearly satisfy the same property. But for sufficiently big $N$ (uniform
over the truncated limb) $\gamma_{-N}$ lies strictly inside $\gamma$ with
a definite space in between. 

Let us consider the annulus bounded by $\gamma$ and $\gamma_{-N}$
with the segments of the rays $R_i$ in between. It follows from the
above that this  configuration $\CC$ stays bounded Teichm\"{u}ller distance
from the standard one $\CC_0$: the round annulus
$\{re^{i\theta}: 1/2\leq r\leq 1\}$ with $p$ equally spaced straight intervals
inside.

So there a $K$-qc map $h: \CC\rightarrow \CC_0$ with dilatation $K$
depending only on the truncated limb, which conjugates $z\mapsto \lambda^Nz$
to $z\mapsto 2z$ on the inner boundaries of the configurations.
 Pulling this map back by linear 
dynamics, we obtain a desired $K$-qc equivalence
between the configuration $\hat\RR$ and $\II$. 

Finally, it is easy to see that the linearizing map $\psi$ is univalent
in a neighborhood of 0 whose size is uniform over the truncated limb.
Hence there is a uniform $l$ such that $\psi$ is  univalent on
$\RR_{-l}\equiv \lambda^{-l}\hat\RR$. Hence the  configuration 
$\RR=f^l\psi(\RR_{-l})$ has bounded geometry.
\QED

\proclaim Corollary 4.3. Let us consider 
 a limb $L_b$ of $M$ and a finite set $A_c$
 of periodic
or pre-periodic  points which stay repelling
through  the unrooted limb $L_b^{tr}$. Let 
$\CC_c $ be the union of configuration $\RR(A_c)$  of rays 
 landing at points of $A_c$  cut at level 1, together with $A_c$
and several equipotentials (whose levels don't depend on $c$).
Then the configuration $\CC_c$ has bounded geometry when $c$ ranges over 
  a truncated limb $L_b^{tr}$

\medskip\noindent{\bf Proof.}
It follows from Lemma 4.2 that near $A_c$ 
 the configuration $\CC_c$ has bounded geometry (where ``near" is
uniform over the truncated limb). On the other hand, it clearly has
bounded geometry outside a uniform neighborhood of  $A_c$,
since  $\phi_c^{-1}$ is a normal family of maps.  \QED

Let $\YY^{(n)}_c\equiv \YY^{(n)}_c(f)$ denote the configuration of 
cutting  curves for the Yoccoz puzzle of depth
$n$ for a map $f$. Now we immediately conclude:

\proclaim Corollary 4.4. For any given $n$ 
 the Yoccoz configuration $\YY^{(n)}(P_c)$ has bounded geometry when $c$ ranges
 over a truncated primary limb $L_b$.

If we consider a DH quadratic-like map $f$  with $\mod(f)>\epsilon$,
 then we can $K(\epsilon)$-qc conjugate it to a quadratic polynomial,
and transfer the net of rays and equipotentials from this quadratic.
In what follows we always assume that the choice of the net is made 
in this way.

\proclaim Corollary 4.5. If $f$ is a quadratic-like map with a definite 
 modulus and internal class $c(f)$ ranging over a primary truncated limb
  $L_b$ of $M$, then Yoccoz
  configuration  $\YY^{(n)}(f)$ has a bounded geometry for any given $n$.

\medskip\noindent{ \bf  4.2. Fundamental domain near the fixed point.}
The goal of this subsection is to construct a combinatorially defined
fundamental domain with bounded geometry  near the fixed point $\alpha$.
It is where the secondary limbs condition comes into the scene.
 
Let $f$ be a DH quadratic-like map whose $\alpha$-fixed point has
rotation number $q/p$.
This map has a single  periodic point $\gamma\in{\rm int}\, Y^{(1)}$ 
of period $p$. 
Let $\gamma'$ be the symmetric point. Consider the family  $\RR(\gamma')$
of rays landing at $\gamma'$. Let $D=D(f)$ be the component of 
$Y^{(1)}\backslash \RR(\gamma')$ attached to the fixed point $\alpha$
(see Figure 5). Then $f^p D$ univalently covers the component
of $Y^{(0)}\backslash \RR(\gamma)$ attached to $\alpha$.
Note that $\partial D\cap \partial (fD)$ is contained in the union of two
rays landing at $\alpha$.

\midinsert
\centerline{\psfig{figure=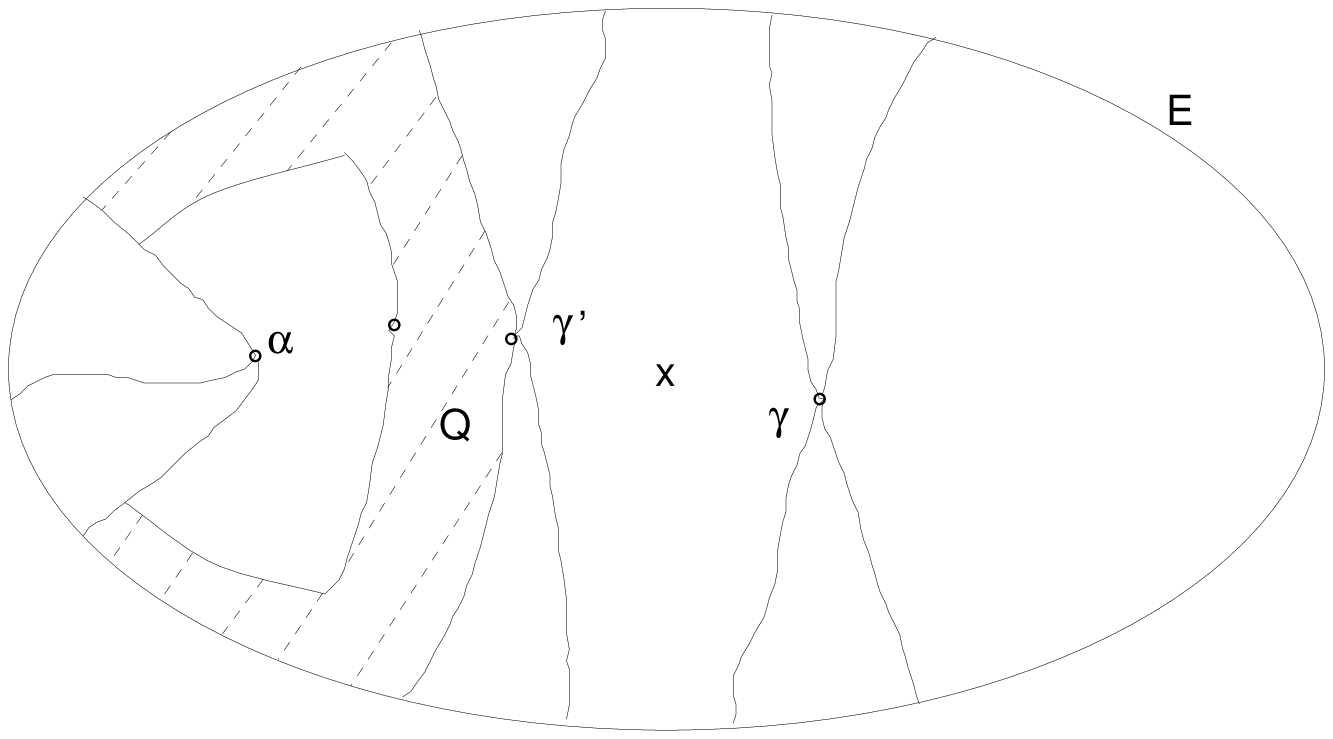,width=.8\hsize}}
\bigskip\centerline{Figure 5: Fundamental domain near the fixed point.}\medskip
\endinsert

 Hence there is a univalent branch of $f^{-p}$
which fixes $\alpha$ and maps $D$ inside itself. It  in now easy to see
 that $f^{-pn} D$ shrink  to $\alpha$ as $n\to\infty$.
So we can select $Q=Q(f)=D\backslash f^{-p}D$ as a fundamental domain for $f^p$
near $\alpha$: any trajectory which starts near $\alpha$ must pass through
$Q=Q(f)$. Now Corollary 4.3 yields:

\proclaim Lemma 4.6. Geometry of the fundamental domain $Q(f)$ is bounded
  if $c(f)$ ranges over a truncated secondary limb and $f$ has a definite
  modulus.  


\medskip\noindent{\bf 4.5. Modulus of the first annulus.}

\proclaim Lemma 4.7. Let $P_c$ be a quadratic polynomial with $c$ outside
  the main cardioid but not immediately renormalizable. If $c$ ranges 
  over a truncated secondary limb $L_b$, then all pieces $X$ of the  initial
  Markov partition (3-3) are well inside 
  $Y^{(0)}$:  $\mod (Y^{(0)}\backslash X)>\nu(L_b)>0.$ 

\noindent{\bf Proof.}  Take a little $\epsilon>0$. Then find an $N$ and 
 $\delta\in (0,\epsilon)$
  such that the equipotential $E_{1/2^N}$ does not intersect the
  $\delta$-neighborhood of 
 $\partial Y^{(0)}\backslash B(\alpha, \epsilon)$
  (for all $c$ in the truncated limb). 

It follows from  Proposition 3.1 and Corollary 4.4 that the statement 
 is true for all pieces of depth  $\leq N$. 

Any other piece $X$ is enclosed
  by the equipotential $E_{1/2^N}$ (where $E\equiv E_1$ is the outer-most
equipotential of external radius 1).
  Hence if $\dist (X, \{\alpha\})>\epsilon$ then
$\dist (X,\partial Y^{(0)})\geq \delta$. As $\diam X$ is uniformly bounded,
we conclude that $X$ is well inside $Y^{(0)}$.

Assume now that $\dist (X, \{\alpha\})<\epsilon$. Then $X$ intersects the domain
$D=D(f)$. Since $\partial D\cap \partial X=\emptyset$, $X\subset D$.
Let us consider the iterates $f^{pk} X$, $k=0,1,\ldots$ until the last 
moment $l$ such that $f^{pl} X\subset D$. At this moment $f^{pl}X$ must
intersect the fundamental domain $Q$. Since their boundaries don't intersect,
we conclude that $f^{pl}X\subset Q$. 

Let us consider domain $Q'\subset Q$ obtained by truncating $Q$ with the
equipotential $f^{-p} E$. This domain has a bounded geometry since 
the fundamental domain $Q$ does (Lemma 4.6). Hence $Q'$ is well inside $f^p D$.
Moreover, $f^{pl}X\subset Q'$ since all puzzle pieces inside $Y^{1}$ are
enclosed by the equipotential $f^{-p} E$ (see Figure 3). 
Hence $f^{pl}X $ is well inside $fD$ as well.

On the other hand, if part of the puzzle-piece $f^{pl}X$ lies in $Q\backslash Q'$,  
it belongs to one of the pieces of depth $n\leq 2p$, of the
 initial partition (3-2). As we have noted above,
 this piece is well inside $fD$.
Hence $f^{pl}X $ is also well inside $fD$.

We conclude that there is always  a definite space around $f^{pl}X$ in
$f^p D$. Pulling this space back by iterates of the univalent
branch $f^{-p} f^p D\rightarrow D$, we obtain a definite space around $X$
in $D$.   \QED

Remember that $A^n=V^{n-1}\backslash V^n$  are the the principal annuli.

\proclaim Theorem I. Let $f$  be  a quadratic-like map with 
internal class $c(f)$ ranging over a truncated secondary limb $L^{tr}_b$.
 If mod$(f)\geq R>0$ 
 then  $$\mod(A^1)\geq  C(R)\,\nu(L^{tr}_b)>0$$
 where $C(R)>0$ monotonically depends on $R>0$ and
$C(R)\nearrow 1$ as $\mu\nearrow\infty$.


\noindent{\bf Proof.} Let us go through the proof of Proposition 3.1.
 We found an $l$ and a puzzle piece $P\subset V^0$ such that $G^l P$
  two-to-one covers $Y^{(0)}$, where $G$ is the Markov map (3-4). Let
  $G^l 0\in X_i^s$ where $X_i^s$ 
  is a puzzle-piece of the initial partition (3-3).
  Then $V^1$ is the pull-back of $X_i^s$ by $G^l|P$. But by  Lemma  4.7
  $X_i^s$ is well inside $Y^{(0)}$. Hence $V^1$ is well inside $V^0$.  \QED

\bigskip

\centerline{\bf \S 5  Increasing of moduli.}\medskip

\noindent{\bf 5.1. Statement of the results.}
Let $f$ be a Douady-Hubbard quadratic-like map.
Let us consider its principal nest
 $Y^0\supset V^0\supset V^1\supset \ldots$,
 and the corresponding nest of annuli
$A^n=V^n\backslash V^{n-1}$. Let us call their moduli
$\mu_n=\mod(A^n)$ the {\it principal moduli} of $g_1$.

In this section we will prove the central result of the paper:

\proclaim Theorem III.  Let $n(k)$ counts the non-central
 levels in the
principal nest  $\{V^n\}$.  Then 
$$\mod A^{n(k)+1}\geq B k, $$
where the constant $B$ depends only on the first modulus $\mu_1=\mod A_1$.

On the way to this result we prove a priori moduli and
distortion  bounds along the principal nest (Theorem II). 
Note that already this result yields the divergence property
of the Yoccoz Theorem (see \S 2). 

Theorem III will also imply a priori bounds for infinitely renormalizable
quadratics of sufficiently big height:

\proclaim Theorem IV. There is a $\chi$-special family $\SS$ of the Mandelbrot 
  copies with the following property. Let  $f$ be an infinitely renormalizable
  quadratic of $\SS$-type then, and $A^{n,m}$ be its principal annuli.
  Then  there is an  $Q=Q(\chi)\to\infty$ as $\chi\to \infty$, such that 
   $\mod(R^n f)\geq \mod A^{n,1}\geq Q$, $n=0,1,\ldots$.

At the end of this section  we will describe other combinatorial
factors which yield big space. This is summarized in Theorem IV$'$
which loosely says  that if the periods of $R^n f$ are sufficiently
big and there are no ``parabolic" or ``Siegel cascades"  in the principal 
nests then there are a priori bounds.  





Recall that $g_n: \cup V^n_i\rightarrow V^{n-1}$
denotes the principal sequence of return maps (3-6),
and $\MM\equiv \MM(f)$ denotes the full family of all puzzle pieces (\S 2.6).
Let $\VV^n\subset \MM$ denote the family of all pieces $V^n_i$ of 
level $n$.

\medskip\noindent{\bf 5.2. First  estimates.} 
Let  us start with a lemma which partly explains the importance of the
principal nest: the principal moduli control the distortion of the
first return maps (see the Appendix for the definition of the 
distortion). Let us consider the decomposition:
$$g_n|V^n = h_n\circ \Phi,  \eqno (5-0)$$
where $\Phi$ is a purely quadratic map and $h_n$ is a diffeomorphism
of $\Phi V^n$ onto $V^{n-1}$.

\proclaim Lemma 5.1. Let $D\in \MM$ be a puzzle piece such that
$f^l D= V^n$, while $f^k D\cap V^n=\emptyset,\, k=0,\ldots, l-1$.
If $\mu_n\leq\bar\mu$ 
then the distortion of $f^l$ on $D$ is $O(\exp(-\mu_{n-1}))$ with
a constant depending on $\bar\mu$. 
Hence the distortion of $h_n$ is  $O(\exp(-\mu_{n-2}))$.

\medskip\noindent{\bf Proof.} This follows from  
Lemma 3.5, Corollary 3.6
 and the Koebe Theorem (see the Appendix).   \QED

Let us fix a level $n> 0$, denote $V^{n-1}=\Delta,\; V_i=V^n_i$,
$g=g_n$,  $A=A^n=\Delta\backslash V_0$, $\mu=\mu_n$, and mark the
objects of the next level $n+1$ with prime. Thus $\Delta'\equiv V\equiv V_0$,
and $g': \cup V_i'\rightarrow \Delta'$.
(We restore the index $n$ whenever we need it).

\proclaim Lemma 5.2.  Let $D'\subset \Delta'$ be a puzzle piece such that
$g^{i(k)} D'\subset V_{i(k)},\, k=0,1,\ldots, l$ with $i(k)\not=0$ for
$0<k<l$.   Then
 $${\rm mod} (\Delta'\backslash D')\geq {1\over 2} 
     \sum_{k=1}^l{\rm mod} (\Delta\backslash V_{i(k)}) .$$

\noindent{\bf Proof.} Let us consider the following nest of topological disks:
   $$\Delta'=W_1\supset  \ldots \supset W_{l}\supset W_{l+1}\equiv D',$$
where $W_k$   is the pullback of $\Delta$ under  $g^{k}$,  $k=1,\ldots l$
(which has itinerary $0=i(0), i(1), \ldots, i(k-1)$ through the pieces
of level $n$).
Then   $g^{k}$ is a two-to-one branched covering of the annulus 
 $W_{k}\backslash W_{k+1}$ over the annulus $\Delta\backslash V_{i(k)}$. Hence
$$\mod  (W_{k}\backslash W_{k+1})={1\over 2} \mod (\Delta\backslash V_{i(k)}),\quad
  (1\leq k\leq l). $$
But by the Gr\"{o}tzsch inequality 
$$\mod (\Delta'\backslash D')\geq\sum_{k=1}^l\mod  (W_{k}\backslash W_{k+1}),$$
and the desired estimate follows.
\QED

This lemma immediately yields:

\proclaim Corollary 5.3.
Given a puzzle piece $V_j'$, we have
 $${\rm mod} (\Delta'\backslash V_j') \geq {1\over 2}\mu.$$
Moreover, if the return to level $n$ is non-central, that is
$g 0\in V_i$ with an $i\not=0$, then
$$\mod (\Delta'\backslash V_j') \geq {1\over 2}(\mu+\mod(\Delta\backslash V_i)).$$

So, a definite principal modulus on some level  produces a definite space
on the next level.

\medskip\noindent{\bf 5.3. Isles and asymmetric moduli.}
 Let $\{V_i\}_{i\in \II}\subset \VV^n$  be a finite family of disjoint 
 puzzle pieces  consisting of at least two pieces (that is $|\II|\geq 2$)
and containing a critical puzzle piece $V_0$.
Let us call such a family {\it admissible}.
We will freely identify the label set $\II$ with the family itself.

Given a puzzle piece $D$, let $\II|D$ denote the family of puzzle pieces
of $\II$ contained in $D$.
Let $D$ be a puzzle piece containing at least two pieces of family $\II$.
For $V_i\subset D$ set 
$$R_i\equiv R_i(\II|D)\subset D\backslash\bigcup _{j\in \II|D}V_j$$ be
an annulus of maximal modulus enclosing $W_i$  but not enclosing other
pieces of the family $\II$. Such an annulus exists by the Montel
Theorem (see Figure 6). We will briefly call it the 
{\sl maximal annulus} enclosing
$V_i$ in $D$ (rel the family $\II$). 

\midinsert
\centerline{\psfig{figure=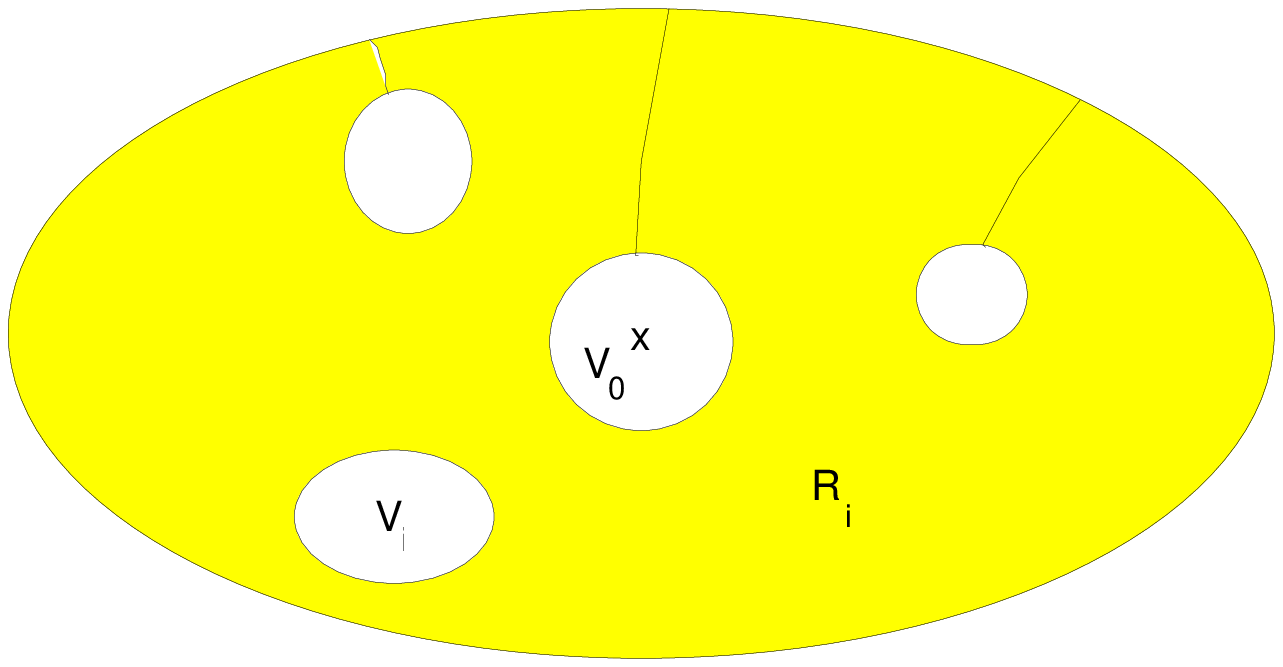,width=.8\hsize}}
\bigskip\centerline{Figure 6: Annulus $R_i$.}
\medskip\endinsert 

Let us define
{\sl the asymmetric modulus of the family $\II$ in $D$ } as
$$\sigma(\II|D)=\sum_{i\in \II} {1\over2^{1-\delta_{i0}}}
  {\rm mod} R_i(R(\II|D), \eqno(5-1)$$
where $\delta_{ji}$ is the Kronecker symbol. So  the critical
modulus is supplied with weight 1, 
while the non-critical moduli are supplied with weights 1/2
(if $D$ is a non-critical island then all weights are
actually 1/2). 

\smallskip
{\bf Remark.} A real analogue of this parameter, ``the asymmetic Poincar\'{e}
  length", appeared in [L3]. Its complex counterpart for the Fibonacci
  map was suggested by Jeremy Kahn. A general notion involving admissible
  familes and isles is given by the author. \smallskip

Take a level $n-1$ {\it which is not DH renormalizable}, that is, 
the family $\VV^n$ consists of more than one piece.
Let  $D=V^{n-1}$, and  let $\{V^n_i\}_{i\in \II}$  be an admissible subfamily
of $\VV^n$.
Then set $\sigma_n (\II)\equiv\sigma(\II|V^{n-1})$ and
$$\sigma_n=\min_\II\sigma_n(\II), \eqno(5-2)$$
where $\II$ runs over all admissible subfamilies of $\VV^n$.

The  principal moduli $\mu_n$ and the asymmetric moduli
 $\sigma_n$ are the main geometric parameters 
of the renormalized maps $g_n$. Again, in what follows the label $n$
will be suppressed as long as the level is not changed.

Let $\{V'_i\}_{i\in \II'}$ be an admissible subfamily of $\VV'$.
Let us organize the pieces of this family in {\sl isles}
in the following way. 
A puzzle piece $D'\subset \Delta'$ is called an {\it island} (for family
$I'$) if

\smallskip\noindent $\bullet$ $D'$ contains at least two puzzle pieces of
  family $\II'$;

\smallskip\noindent $\bullet$ There is a $t\geq 1$ such that $g^k D'\subset
  V_{i(k)}, \; k=1,\ldots t-1,$ with $i(k)\not=0$, while $g^t D=\Delta.$ 

Given an island $D'$, let 
$ \phi_{D'} = g^t : D'\rightarrow\Delta$. This map
is either a double covering or a biholomorphic isomorphism
 depending on whether $D'$ is critical or not.  In the former case,
 $D'\supset V'_0$ (for otherwise $D'\subset V'_0$ contradicting the 
 first part of the definition of isles). 

We call a  puzzle piece $V_j'\subset D'$  $\phi_{D'}$-{\sl pre-critical} 
if $\phi_{D'}(V_j')=V_0$. 
There are at most two pre-critical pieces in any $D'$.
If there are actually two of them, then they are non-critical and
symmetric with respect to the critical point 0. Thus in this case
$D'$ contains also the critical puzzle piece $V_0'$.


Let ${\cal D}'={\cal D}(\II')$ be the family of isles associated with $\II'$.
Let us  consider the asymmetric moduli $\sigma(\II'|D')$ as a function on
this family. This function is clearly monotone:
$$\sigma(\II'|D')\geq\sigma(\II'|D_1')\quad {\rm if} \quad D'\supset D_1',
    \eqno (5-3)$$
 and superadditive:
$$\sigma(\II'|D')\geq\sigma(\II'|D_1')+\sigma(\II'|D_2'), $$
provided $D_i'$ are disjoint subisles in $D'$.

Let us call an island $D'$ {\it innermost} if it does not contain any
other isles of the family $\DD(\II')$. As this family is finite,
innermost isles exist.

\medskip\noindent{\bf 5.4. Non-decreasing of moduli.}

\proclaim Lemma 5.4. Let $\II'$ be an admissible family of puzzle pieces. 
  Let $D'$ be an innermost island associated to the family $\II'$, and
let  $\JJ'=\II'|D$.
Let $i(j)$ is defined for $j\in \JJ'$ by the property
 $\phi_{D'}(V_j')\subset V_{i(j)}$, and let
 $\II=\{i(j) : j\in \JJ'\}\cup\{0\}$.
 Then $\{V_i\}_{i\in \II}$ is an admissible family of puzzle pieces, and
$$ \sigma(\II'|D')\geq 
{1\over 2}\left( (|\JJ'|-s)\mu +  s\; {\rm mod} R_0+
 \sum_{j\in \JJ',\, i(j)\not=0}{\rm mod} R_{i(j)}\right), 
\eqno (5-4)$$ 
where $s=\#\{j: i(j)=0\}$ is the number of $\phi_{D'}$-pre-critical pieces,
and $R_i$ are the maximal annuli enclosing $V_i$ in $\Delta$ rel family $\II$.

\noindent{\bf Proof.}   Let $\phi\equiv \phi_{D'}$.
Let us show first that the family $\II$ is admissible. This family is 
finite since $\JJ'\subset \II'$ is finite. The critical
puzzle piece belongs to $\II$ by definition. So the only 
property to check is that $|\II|\geq 2$. But otherwise $\JJ'$ would consist
of  two pre-critical puzzle pieces.  But then $D'$ would be critical,
and  thus should have   also contained the
critical piece $V_0'$, which is a contradiction.

Let us observe next that 
$$\mod(V_{i(j)}\backslash \phi V'_j)\geq\mu \quad {\rm if}\quad 
       i(j)\not=0.\eqno(5-5)$$
   Indeed, in this case $g^m (\phi V'_j)=V_0$ for
   some $m>0$. Let $W\subset V_{i(j)}$ be the pull-back of $\Delta$ under
   $g^m$.     Then the annulus $W\backslash \phi V'_j$ is univalently mapped by
   $g^m$ onto the annulus $\Delta\backslash V_0$. Hence
   $\mod(W\backslash \phi V'_j)=\mod(\Delta\backslash V_0)=\mu,$ and  (5-5)
follows. 

Given a $i\in I$, let us consider  a topological
     disk $Q_i = R_i\cup V_i\subset \Delta$ (``filled annulus $R_i$").
By the Gr\"{o}tzsch inequality and (5-5),
$$\mod(Q_{i(j)}\backslash \phi V_j) \geq \mod R_{i(j)}+(1-\delta_{0,i(j)})\mu.
  \eqno (5-6) $$

For a $j\in J'$, let us consider an annulus $B_j\subset D'$, 
   the component of $\phi^{-1} R_{i(j)}$ enclosing $V_j'$. 
This annulus does not enclose
any other pieces $V_k'\in \JJ'$, $k\not=j$. Indeed, otherwise the
inner component of $\C\backslash B_j'$ would be an island contained
in $D'$, despite the assumption that $D'$ is innermost. 

Let us now consider a topological disk $P_j'$
 obtained by filling the annulus $B_j'$. 
As it contains a single puzzle piece $V_j'$
of family $J'$, the annulus $P_j'\backslash V_j'$ does not go around any
other  puzzle piece
 $V_k'\in \JJ'$, $k\not=j$. Hence

$$\mod R_j'\geq \mod(P_j'\backslash V_j'), \eqno (5-7)$$
where $R_j'\subset D'$ is the maximal annulus enclosing $V_j'$ rel $\JJ'$.
Moreover $\phi: P_j'\rightarrow Q_{i(j)}$ is univalent or double covering
depending on whether $j\not=0$ or $j=0$. Hence

$$\mod(P_j'\backslash V_j')\geq {1\over 2^{\delta_{j0}}}
  \mod  (Q_{i(j)}\backslash \phi V_j).                           \eqno (5-8)$$
Inequalities  (5-6)-(5-8) yield
$$\mod R_j'\geq {1\over 2^{\delta_{j0}}}( \mod R_{i(j)}+(1-\delta_{0,i(j)})\mu).
                             \eqno(5-9)$$
Summing up estimastes (5-9) over $\JJ'$ with weights ${1/ 2^{1-\delta_{j0}}}$,
we obtain the desired inequality (5-5).  \QED
 

\proclaim Corollary 5.5. For any island $D'$ of the family ${\cal D}'$
 the following estimates hold:
$$\sigma(\II'|D')\geq {1\over 2}\mu\quad {\rm and}
\quad \sigma(\II'|D')\geq\sigma(\II)\geq\sigma.$$

\noindent{\bf Proof.} By  monotonicity (5-3), it is enough to check the case
of an innermost island $D'$. Let us use the notations of the previous lemma.
Since the family $\II$ is admissible, it contains a non-critical piece.
 Hence $|\JJ'|$ is always strictly greater than the number $s$ of
pre-critical pieces in $D'$, and (5-4) implies the first of the above
inequality. 

Furthermore, as 
 $\mu\geq {\rm mod}(R_0)$ and $|\JJ'|\geq 2$,
 the right-hand side in (5-4) is bounded from
 below by
$$ {1\over 2}\left( |\JJ'|\; {\rm mod}(R_0)+\sum_{i\in \II, i\not=0}
 {\rm mod}(R_i)\right)\geq \sigma(\II). $$ 
(Note that $\sigma(\II)$ makes sense since $\II$ 
 is admissible).
Finally $\sigma(\II)\geq\sigma$, and the second inequality
follows.
  \QED 


Let us fix a ``big" integer quantifier $N_*>0$.
 We say that a level $n$ is in the
 ``tail of a  cascade" if all levels $ n-1,   n-N_*$  belong to
a cascade (note that level $n-1$ itself may be non-central).
 Cascades of length at least $N_*$ we call ``long".

\proclaim Theorem II. Given a generalized quadratic-like
map $g_1$, we have the following bounds of the   geometric parameters
within its principal nest:

{\it$\bullet$ The asymmetric moduli $\sigma_n$ grow monotonically
  and hence stay away
  from 0 on all levels (until the first DH renormalizable level):
  $\sigma_n\geq\bar\sigma>0$.

  $\bullet$ The principal moduli $\mu_n$ stays away
   from 0  (that is, $\mu_n\geq \bar\mu>0$)
   everywhere except for  the case when $n-1$ is in the
  tail of a long  cascade (the bound $\bar \mu$ depends on the choice of
    $N_*$). 

$\bullet$ The non-critical puzzle-pieces $V^n_i$
   are well inside $V^{n-1}$ (that is, $\mod(V^{n-1}\backslash V^n_i\geq\bar\mu>0$)
   except for the case when $V^n_i$ 
  is pre-critical and $n-2$ is the last level of
   a long cascade.

$\bullet$ The distortion of $h_n$ from (5-0)
  is uniformly bounded on all levels by
  a constant $\bar K$.

\noindent All bounds depend only on the  first principal modulus $\mu_1$
and (as $\bar\mu$ is concerned) on the choice of $N_*$).}

\noindent{\bf Proof.} The first assertion follows from the second inequality
of Corollary 5.5.  Together with Corollary 5.3 it implies the second one
(note that the second inequality of this corollary implies that
$\mu'\geq \sigma/2$ in the non-central case).
One more application of Corollary 5.3 yields the next assertion.

Let us check the last statement. If $n-2$ is not in the tail of a
central cascade, then $\mu_{n-1}\geq \bar\mu$ by the second 
statement, and the desired follows from Lemma 5.1.

Let $n-2$ be in the tail of a central cascade 
$V^m\supset\ldots\supset V^{n-2}\supset\ldots$. If this is not the
last level of this cascade then $g_n|V^n= g_{m+2}|V^n$, so that
 $h_n$ is just a restriction of the map $h_{m+2}$ with bounded distortion.

Finally, if $n-2$ is the {\it last} level of a central cascade, then by 
Corollary 3.8 $h_n$ can be extended to a univalent map with range
$V^m$, and the Koebe Theorem implies the distortion bound.
  \QED

\medskip\noindent{\bf  5.5. Linear growth of moduli.}
Our goal is to prove that 
$\sigma'\geq\sigma+a $  with
a definite $a>0$ (that is, dependent only on $\mod A_0$)
at least on every other level, except for the tails
of long cascades and a couple of the following levels. 
(Theorem II shows the reason why these tails play a special
role: In the tails the principal moduli become tiny  which slows down
the  growth rate of asymmetric moduli.) 
 
Clearly it is enough to show that for any innermost island $D'$
$$\sigma(\II'|D')\geq\sigma+a \eqno (5-10)$$
with a definite $a>0$.
The analysis will be split into a tree of cases.

\medskip\noindent{\bf Case I. An island with  at least three puzzle pieces.}

\proclaim Proposition 5.6. If an innermost island $D'$ contains
  at least three puzzle-pieces $V'_j,\; j\in \JJ'$,  then
  $$\sigma(\JJ'|D')\geq\sigma(I)+{1\over 2}\mu.$$

\noindent{\bf Proof.}  Let us  split off
$(1/2)\mu$  in (5-4) and estimate all other $\mu$'s by mod$(R_0)$. 
This estimates the right-hand side by
$${1\over 2}\mu+{|\JJ|-1\over 2}{\rm mod}(R_0)+
{1\over 2}\sum_{i\in \II, i\not=0}{\rm mod}(R_i), $$
which immediately yields what is claimed.     \QED

Hence under the circumstances of Proposition 5.6 we observe 
a definite growth  of the asymmetric modulus
 provided level $n-1$ is  not in the tail of a long cascade. 
Indeed then by Theorem II $\mu$ is bounded away from 0,
and (5-10) follows.

\medskip
\noindent{\bf Case II. An island with two puzzle pieces.}  
The further analysis needs
some preparation in the geometric function theory summarized in the
Appendix. 

Suppose  we have an innermost island $D'$ containing
two puzzle-pieces $V_j',\; j\in \JJ'$. Let $\phi\equiv \phi_{D'}$ and let
$\phi V'_j\subset V_i$ with $i=i(j)$. Fix a quantifier $L_*>0$.
When we say that something is ``big", this means that it is at
least  $C(L_*)$ where $C(L_*)\to\infty$ as $L_*\to\infty$. Similarly ``small"
means an upper bound by $\epsilon(L_*)\to 0$ as $L_*\to\infty$.
The sign $\approx$ will mean an equality  up to an small (in the above sense)
error,
 while the sign $\succ$ will mean the inequality up to
a small error.
 
\medskip {\bf  Subcase (i)}.
 {\it Assume that  there is a non-critical puzzle-piece $V_{i(j)}$
   whose  Poincar\'{e} distance in $\Delta$  from the critical point
is less than $L_*$.} Then by Lemma A.1
$$\mu\geq {\rm mod}(R_0) + \alpha  \eqno (5-11)$$ with a definite
$\alpha=\alpha(L_*)>0$. But observe that when we passed from Lemma 5.4
to Corollary 5.5 we estimated $\mu$ by mod$(R_0)$. Using the better
estimate (5.10), we obtain a definite increase of $\sigma$.

\medskip {\bf Subcase (ii).}
  {\it Assume now that the hyperbolic distance in $\Delta$ from any
 non-critical puzzle piece $V_{i(j)}$, 
 from the critical point is at least $L_*$}. Let also $n-2$ don't belong to 
the tail of a
long cascade (for the sake of linear growth it is enough to prove 
definite growth on such levels).
Then $V_0$ may not  belong
to any non-trivial island together with some non-critical piece
  $V_{i(j)}$. Indeed, 
 by Theorem II all puzzle pieces of level $n-1$ are well
 inside $V^{n-2}$. But then by Lemma 5.2 all non-trivial isles of level $n$  
are well inside of $V^{n-1}\equiv \Delta$. (The quantifier $L_*$ should be
chosen bigger than the a priori bound on the hyperbolic diameters of the isles).

\medskip {\bf Subcase (ii-a).}
 {\it Assume that both $V_{i(j)}$ are non-critical}. Then 
by Corollary 5.5 $\sigma(\JJ'|D')$
is estimated by $\sigma_n( \II)$  
where the family $\II$  consists of three puzzle pieces:
two pieces $V_{i(j)}$ and the central puzzle piece $V_0$. 

If puzzle pieces $V_{i(j)},\; j\in \JJ'$, don't 
belong to the same non-trivial island, then by Proposition 5.6
$\sigma(\II)\geq \sigma_{n-1}+a$ with a definite $a>0$, and we are done.

 Otherwise the puzzle pieces $V_{i(j)}$ belong to an island $W$. 
Since by Lemma 5.2
 $W$ is well inside of $\Delta$, it stays on the big Poincar\'{e}
distance from the critical point (namely, on distance $L_*-O(1)$). 
Hence mod$(R_0)\approx\mu$, and
$$\sigma( \II)\geq \sigma(\II|W)+{\rm mod}(R_0)\succ
   \sigma_{n-1}+\mu$$
where $\mu\equiv \mu_{n}$ 
is bounded away from 0, since level $n-1$ is not in the tail
of a long cascade.
So we have  gained some extra growth, and can pass to the next case.

Below we will restore labels $n$ and $n+1$ since many levels
will be involved in the  consideration.

\medskip
 {\bf Subcase (ii-b).}
  {\sl Let one of the puzzle pieces $V_i^n$ be critical}.
  So we have the family $\II^n$ of two puzzle-pieces $V_0^n$ and $V_1^n$.
  Remember that we also assume that {\sl the hyperbolic distance
  between these pieces is at least $L$}. Hence, $V^{n-1}$ is the
  only island containing both of them, so that $g_{n-1} V_0^n$ and 
 $g_{n-1} V_1^n$
  belong to different puzzle-pieces of level $n-1$.  For the same reason 
   we can assume that one of these puzzle-pieces is critical. Denote
  them by $V^{n-1}_0$ and $V^{n-1}_1$. 
  Then one of the following  two possibilities on level $n-2$ can occur:

\smallskip\noindent
1) {\it Fibonacci return} 
   when $g_{n-1} V^n_0\subset V^{n-1}_1$ and
   $g_{n-1}V^n_1=V^{n-1}_0$ (see Figure 7);

\smallskip\noindent 2) {\it Central return}  when
    $g_{n-1} V^n_0= V^{n-1}_0$ and
   $g_{n-1}V^n_1\subset V^{n-1}_1$. 

\midinsert
\centerline{\psfig{figure=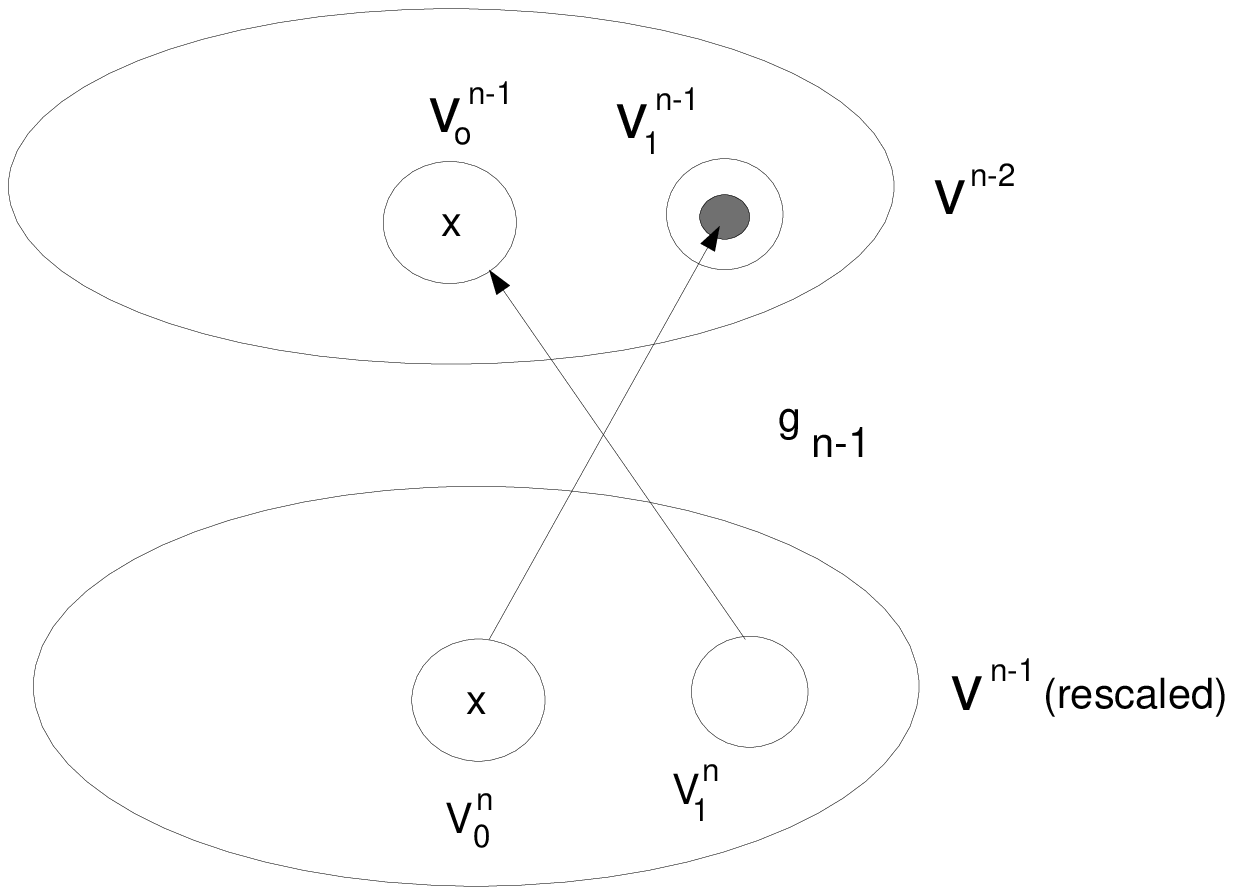,width=.8\hsize}}
\bigskip\centerline{Figure  7: Fibonacci scheme.}\medskip
\endinsert

  We can assume that one of these schemes occur on several previous levels
    $n-3, n-4,...$ as well (otherwise we gain an extra growth by the previous 
   considerations). To fix the idea, let us first consider
  the following particular case, which plays the key role for the whole 
  theorem.

\medskip\noindent{\bf Fibonacci cascade.} 
Assume that  on both
   levels $n-1$ and $n-2$ the Fibonacci returns  occur. 
  Let us look more carefully at the estimates of Lemma 4. 
 In the Fibonacci case we just have:

$${\rm mod}(R^n_1)\geq {\rm mod} (R^{n-1}_0),  \eqno (5-12)$$
$${\rm mod} (R^n_0)\geq {1\over 2} {\rm mod} 
   (Q^{n-1}_1\backslash g_{n-1} V^n_0  ), \eqno (5-13)$$
where $Q^n_i=V^n_i\cup R^n_i$. Applying $g_{n-2}$ we see that
$${\rm mod} ( Q^{n-1}_1\backslash g_{n-1} V^n_0 )\geq
 {\rm mod} (Q^{n-2}_0\backslash V_0^{n-1}). \eqno (5-14)$$
But since $V^{n-2}_1$ is hyperbolically far away from the critical point,
$${\rm mod} (Q^{n-2}_0\backslash V_0^{n-1})\approx
 {\rm mod}(V^{n-3}_0\backslash V_0^{n-1}). \eqno (5-15) $$
By the Gr\"{o}tzsch Inequality there is an $a\geq 0$ such that
$${\rm mod}(V^{n-3}_0\backslash V_0^{n-1})=\mu_{n-1}+\mu_{n-2}+a. \eqno (5-16)$$
Clearly
$$\mu_{n-1}\geq {\rm mod} (R_0^{n-1}). \eqno (5-17) $$
Furthermore, let $P_1^{n-1}\subset V^{n-2}$
 be the pull-back of $Q_0^{n-2}$ by $g_{n-2}$.
Since $\partial 
P_1^{n-1}$ is hyperbolically far away from $V_1^{n-1}$,
we have:
$$\mu_{n-2}\geq {\rm mod}(R_0^{n-2})={\rm mod} (P_1^{n-1}\backslash V_1^{n-1})
      \approx \mod(V^{n-2}\backslash V_1^{n-1})\geq {\rm mod} ( R_1^{n-1}).
    \eqno(1-18) $$
Combining estimates (5-13) through (5-18) we get
$${\rm mod}(R_0^n)\succ {1\over 2} ({\rm mod}(R_0^{n-1})+
   {\rm mod}(R_1^{n-1})+a).  \eqno (5-19)$$

We see from (5-12) and (5-19) that we need to check  that
the constant  $a$ in (5-16) is definitely positive.
Assume that this is not the case, that is, for any $\delta>0$ we can find
a level $n$ in the Fibonacci cascade as above  such that $a<\delta$.
 Set $\Gamma_n=\partial V^n$.
Then by the Definite Gr\"{o}tzsch Inequality (see the Appendix), 
the width$(\Gamma_{n-2})$
in the annulus $T=V^{n-3}\backslash V^{n-1}$  is at most $\xi(\delta)$ with
$\xi(\delta)\to 0$ as $\delta\to 0$. 
Since $\Gamma_{n-2}$ is well inside of $T$, we conclude by the Koebe
Distortion Theorem that $\Gamma_{n-2}$ is contained in a narrow
neighborhood of a curve $\gamma$ with a bounded geometry. 
Hence there is a $k=k(\delta)\to 0$ as $\delta\to 0$ and and
 $\epsilon=\epsilon(\delta, k)>0$ such that
the curve  $\Gamma_{n-2}$ is not $(k,\epsilon)$-pinched.

On the other hand, the 
 hyperbolic distance from the puzzle piece 
$V^{n-1}_1$ to the critical point 0 in $V^{n-2}$ is at least $L_*$. 
Hence  by Lemma A.4 it must be located  Euclideanly  very close 
 to $\Gamma_{n-2}$ relatively the Euclidean 
distance to the critical point (that is, the relative distance is at most
$\beta(\delta)\to 0$ as $\delta\to 0$). Hence the critical value $g_{n-1} 0$
is also very close to $\Gamma_{n-2}$ relatively the distance to the
critical point, that is 
$${ \dist (g_n 0, \Gamma_{n-2})\over\dist (g_n 0, 0)}\leq \epsilon (L_*),$$
where $\epsilon(L_*)\to 0$ as $L_*\to \infty$.

By the last statement of Theorem II,
$g_{n-1}$ is  a quadratic map up to a bounded distortion. Hence the
curve $\Gamma_{n-1}$ which is the pull-back of $\Gamma_{n-2}$ by $g_{n-1}$
must have a huge eccentricity around the critical point. 
But then by Lemma 7.2 
 the width of $\Gamma_{n-1}$ in $V^{n-2}\backslash V^{n}$ is also big,
which by the above considerations gives a definite 
linear growth on the next level.

\smallskip\noindent {\bf Remark.} The actual shape of a deep level puzzle-piece
  for the Fibonacci cascade is shown on Figure 8. There is a good reason 
 why it resembles the filled Julia set for $z\mapsto z^2-1$
  (see [L5]). As the geodesic in $V^{n-1}_0$ joining 
   the puzzle-pieces $V^n_0$ and $V^n_1$
   goes  through the  pinched region,
 the Poincar\'{e} distance
  between  these puzzle-pieces is, in fact, big.

\midinsert
\centerline{\psfig{figure=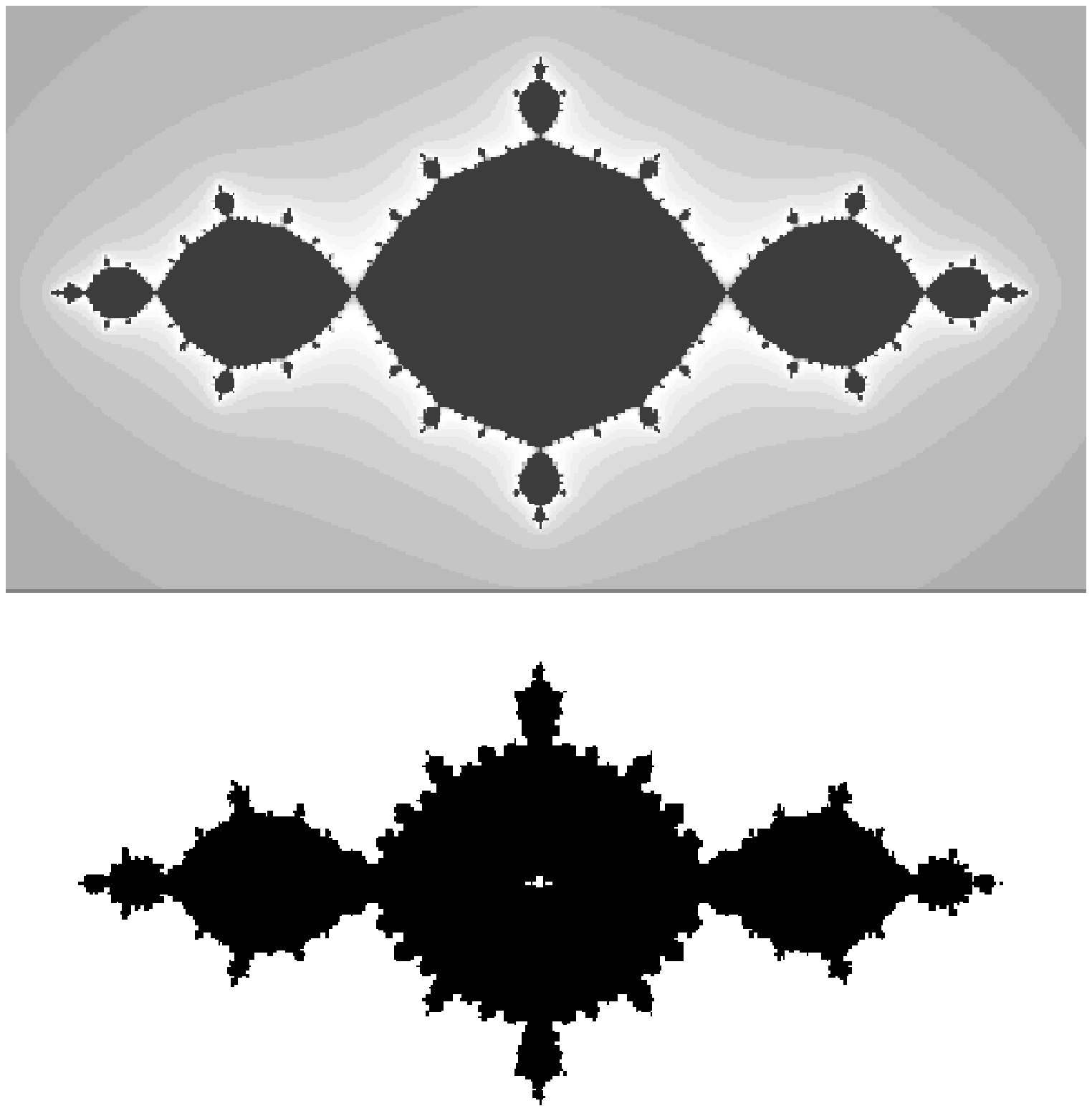,height=.9\hsize}}
\bigskip\centerline{Figure  8: Fibonacci puzzle piece (below)
 vs the Julia set of 
$z\mapsto z^2-1$ (above) .}\medskip
\endinsert

Now it is the time to look closer at central cascades. 

\medskip\noindent{\bf Central cascades.} Let $N\geq 2$, $n=m+N$, and
let us consider a nest $\CC^{m+N}$ of puzzle pieces   
$$V^{m}\supset V^{m+1}\supset\ldots\supset V^{m+N-1}\supset V^{m+N}
  \supset D^{m+N}
\eqno (5-20)$$  satisfying the following
properties (see Figure 8):

\noindent$\bullet$ The return on level $m-1$ is non-central:
   $g_m 0\not\in V_0^m$;

\smallskip\noindent $\bullet$ Central returns occur on levels $m, m+1,
 \ldots m+N-2$, 
 that is $g_{m+1} 0\in V^{m+N-1}$;

\smallskip\noindent $\bullet$ $D^{m+N}$ is an island 
  with a family $\II^{m+N+1}$  of
 two  puzzle pieces inside. Let $\phi_{m+N}\equiv \phi_{D^{m+N}}$
denote the corresponding double covering  $D^{m+N}\rightarrow V^{m+N-1}$;

\smallskip\noindent $\bullet$ One of the puzzle pieces $\phi_{m+N} V^{m+1}_0$,
   $\phi_{m+N} V^{m+1}_1$ is critical.

\midinsert
\centerline{\psfig{figure=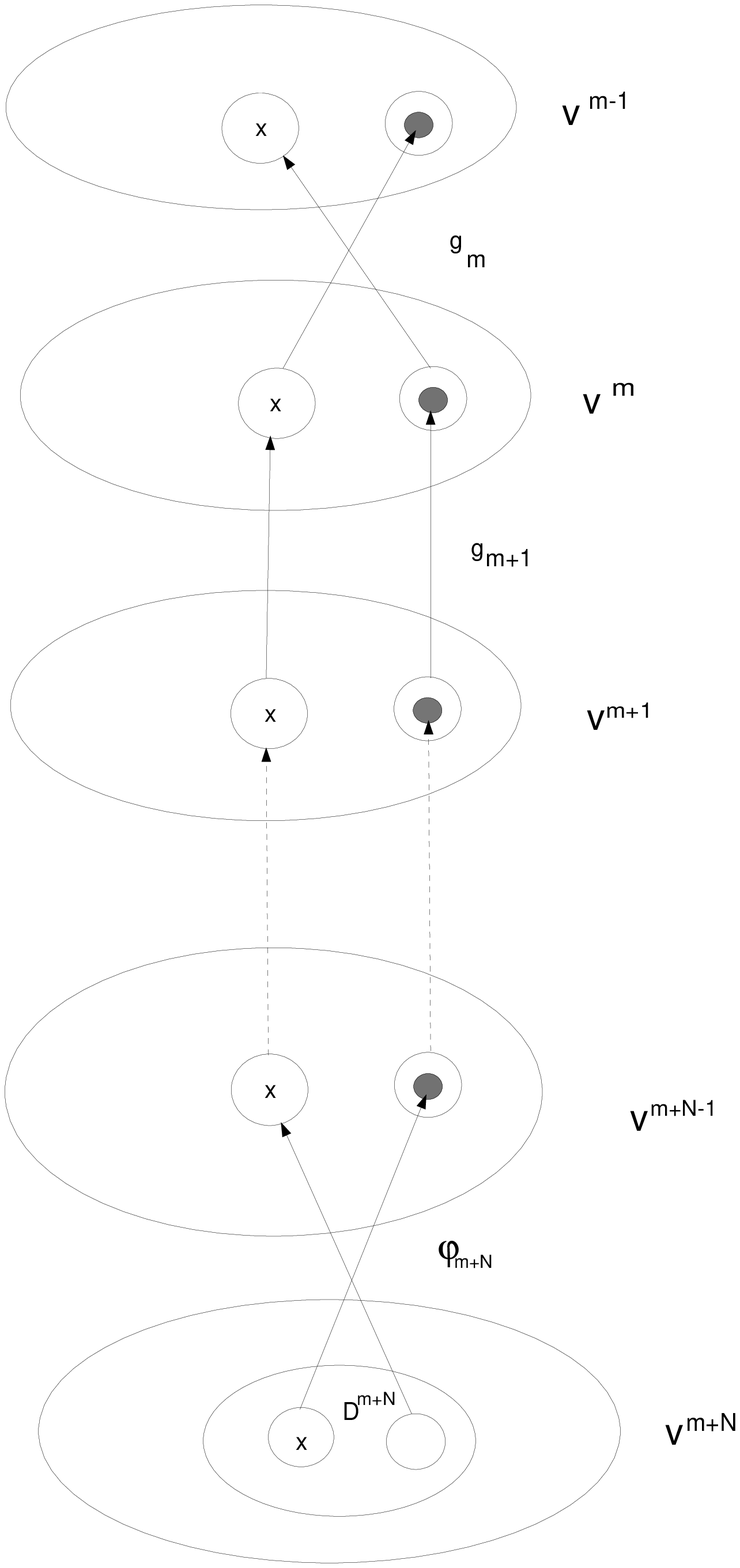,height=.8\vsize}}
\medskip\centerline{Figure  9:  Central cascade}
\centerline{ (with Fibonacci returns on top and bottom)  }\bigskip
\endinsert

Though the logic of our argument so far allowed to assume that $N\leq N_*$,
the following argument will require  consideration of arbitrary big $N$.
So let $N$  be arbitrary integer $\geq 2$.

We would like to analyze when
$$\sigma(I^{m+N+1}|D^{m+N})\geq \sigma_{m+1} + a \eqno (5-21)$$
with a definite $a>0$. To this end
 we need to pass from level $m+N$ all way up
to level $m$.

Let us consider the
 family $ \WW(\CC^{m+N})$
of puzzle pieces $W^k_i,\; k=m+1,...,m+N,$ $i\not=0$,
 the pull-backs
  of the $V^{m+1}_i\equiv W^{m+1}_i$ to the annuli $A^k$ (compare
with the initial Markov partition in \S 3.2). Given a $W^k_i\subset A^k$,
there is a $V^{m+1}_j\subset A^{m+1}$, $j\not=0$,
 such that $g_{m+1}^{m+1-k} W^k_i =V^{m+1}_j$.
Hence $g_m\circ g_{m+1}^{m+1-k} W^k_i =V^m$.
 So we can define a Bernoulli  map
  $$G: \cup W^k_i\rightarrow  V^m \eqno (5-22)$$
by letting $G| W^k_i=g_m\circ g_{m+1}^{m+1-k}| W^k_i.$

Let $V^{m+N+1}_*\subset D^{m+N}$ be a non-pre-critical piece of the family
$\II^{m+N+1}$, and 
$$\phi  V^{m+N+1}_*\subset V^{m+N}_1\subset W_i^{m+N}$$
for some $i\not=0$. Then the return map 
$g_{m+N+1}: V^{m+N+1}_*\rightarrow V^{m+N}$ can be decomposed as
$G^l\circ \phi$ for an appropriate $l\geq 1$. Since the map $G$
is Bernoulli with range $V^m$, 
$$\mod(W_i^{m+N}\backslash \phi V^{m+N+1}_*)\geq \mod(V^m\backslash V^{m+N}) \eqno (5-23).$$
Let $\Gamma^k=\partial V^k$, and
$$w_k=\width (\Gamma^k| V^{k-1}\backslash V^{k+1}). \eqno(5-24)$$

For $k\in [m+1, m+N]$ let $V^k_1 $ denote the puzzle piece of level $k$
containing
 $g_{m+1}^{\circ(k-m-N)}\phi V^{m+N+1}_*$, and $\II^k$ denote the family
of two puzzle pieces: $V_0^k$ and $V_1^k$. Moreover, let $R_i^k\subset V^{k-1}$
 denote an annulus
of maximal modulus going around $V^k_i$ but not going around  the other 
piece of family $\II^k$, $i=0,1$. 

By the Definite Gr\"{o}tzsch Inequality and the second part of Theorem II,
there is an $a=a(w_{m+1})$ such that
$$\eqalign{\mod(V^m\backslash V^{m+N})\geq \sum_{k=m+N}^{m+1} \mod A^k +a=\cr
  \sum_{k=0}^{N-1}{1\over 2^k}{\rm mod}(A^{m+1})+a\geq
   (2-{1\over 2^{N-1} })\mod R_0^{m+1}+a.\cr }\eqno (5-25)$$

Let $S^{m+N}_0$ and $S^{m+N}_1$ denote the  pull-backs of 
 the annuli $R^{m+1}_0$ and $R^{m+1}_1$ by the map
$g_{m+1}^{\circ(N-1)}: V^{m+N-1}\rightarrow V^m$. Then

$$\mod S^{m+N}_0\geq {1\over 2^{N-1}} \mod R^{m+1}_0 \quad {\rm and}\quad
  \mod S^{m+N}_1\geq \mod R^{m+1}_1. \eqno (5-26)$$

Note that the inner boundary of $S^{m+N}_1$ coincides with the outer
boundary of \par\noindent
  $W_i^{m+N}\backslash \phi V^{m+N+1}_*$.  Let $Q_1^{m+N}$
denote the union of these two annuli. This annulus goes around 
$\phi V^{m+N+1}_*$ but not around $V^{m+N}_0$. Now estimates
(5-23), (5-25), (5-26) yield

$$\mod S_0^{m+N}+ \mod Q_1^{m+N}\geq 2\mod R_0^{m+1}+ \mod R_1^{m+1}
                 +a\geq 2\sigma(\II^{m+1})+a.$$
Finally, pulling $S_0^{m+N}$ and $Q_1^{m+N}$
 back by $\phi_{m+N}$ to the island $D^{m+N}$ 
we obtain:

$$\sigma(\II^{m+N+1}|D^{m+N})\geq {1\over 2} (\mod S_0^{m+N}+\mod Q_1^{m+N})\geq
                   \sigma(\II^{m+1})+a/2. \eqno (5-27)$$ 
So we come up with the following statement:

\proclaim  Statement 5.7. There is an increasing function 
$a: \R_+\rightarrow \R_+$, $a(0)=0$, such that for the cascade $\CC^{m+N}$
estimate (5-21) holds with $a=a(w_{m+1})$, where 
$w_{m+1}=\width(\Gamma^{m+1}|V^m\backslash V^{m+2})$.

Let us fix a quantifier $ w_*$ which distinguishes  ``small width" $w$ from a
 ``definite" one. For  further analysis
let us go several levels up. 
Let $m-1-l$ be the highest non-central level preceding
$m-1$, $l\geq 1$. We are going to  study when
$$\sigma(\II^{m+N+1}|D^{m+N})\geq \sigma_{m-l}+a \eqno (5-28)$$
with a definite $a>0$.
We cannot now assume that $l$ is bounded, so we face 
a possibility of a long cascade $\CC^{m-1}$: 
$V^{m-l}\supset\ldots\supset V^{m-1}$. 
Set $g=g_{m-l+1}$; then $g0\in V^{m-1}$.

{\it Assume  first that $m-2$ is not the last piece of a long cascade}
(in particular, this is the case when central return occurs on level
$m-2$, that is, $l\geq 2$).
Then by the third part of
Theorem II all non-central pieces of level
$m$ are well inside $V^{m-1}$: $\mod (V^{m-1}\backslash V^m_j)\geq \bar\mu$.
 Hence  $g_m V_0^{m+1}$ and $g_m V_1^{m+1}$ belong  to different
pieces of level $m$. Indeed, otherwise the hyperbolic distance between
$V_0^{m+1}$ and $V_1^{m+1}$ in $V^m$ 
would be bounded by a constant $L(\bar\mu)$.
But according to our assumption this distance is  at least $L_*$. 
So this situation
is impossible if $L_*$ was a priori selected bigger than $L(\bar\mu)$.

For the same reason the pieces $g^k\circ g_m V_i^{m+1}$, $i=0,1$,
also belong to different pieces  $V^{m-k}_j$ for $0\leq k\leq l-3$.
Indeed, assume they belong to the same piece $V^{m-k}_j$. Clearly this piece
is non-central, that is $j\not=0$. Then it is contained in a piece
$W^{m-k}_j$ of the Bernoulli family  $\WW(\CC^{m-1})$
associated to the central cascade $\CC^{m-1}$. Hence $g_m V_0^{m+1}$ and
$g_m V_1^{m+1}$ belong to $W^{m}_j$, the pull-back of $W^{m-k}_j$ by
$g^k$.  As
 $\mod(W^{m}_j\backslash g_m V^{m+1}_i)\geq\bar\mu$, the hyperbolic distance
between $ V_0^{m+1}$ and $ V_1^{m+1}$ in
$V^{m}$ is at most $L(\bar\mu)$ contradicting our assumptions.

Let us show now that (5-28) holds if both $g_m V_i^{m+1}$ are
non-central. Indeed let us then consider the family $\II^{m}$
of three pieces: two pieces of level $m$ containing   $g_m V_i^{m+1}$
and the central piece $V^m$.  Let $\II^{m-k}$ denote the family
of puzzle pieces of level $m-k$ containing  the pieces of 
$g^k \II^{m}$. By the previous two paragraphs, $\II^{m-k}$ consists of
three puzzle pieces. Then by Corollaries 5.5 and 5.7,
$$\sigma(\II^{m+1})\geq \sigma(\II^m)\geq\ldots\geq  \sigma(\II^{m-l+2})
    \geq \sigma(\II^{m-l+1})+{1\over 2}\bar\mu, $$ 
and we are done.  

Thus let us assume that the Fibonacci return occurs on level $m-1$.
In this case let $\II^{m-k}$ denote the family of two puzzle pieces
$V_0^{m-k}$ and $V_1^{m-k}$ containing $g^k\circ g_m V^{m+1}_i$,
$i=1,2$, $k\leq l-1$. 

Note that in order to have (5-28) it is enough to to have a definite
increase of the $\sigma(\II^{m-k})$ in the beginning of the cascade
$\CC^{m-1}$. By  Statement 5.7 applied to this cascade this is the case
if $\width(\Gamma^{m-l+1}|V^{m-l}\backslash V^{m-l+2})\geq w_*.$ So assume
that the opposite inequality holds. Similarly,  because of Subcase (i),
we can assume that
the hyperbolic distance from $V^{m-l+2}_1$ to 0 in $V^{m-l+1}$ is at least
$L_*$. 

It follows from Lemma A.4  from the Appendix that
 the piece $V^{m-l+2}_1$ stays  Euclidean distance at most
$\epsilon\, \diam \Gamma^{m-l+1}$
 from $\Gamma^{m-l+1}$ where 
$\epsilon=\epsilon_{\bar\mu}(w_*, L_*)\to 0$ as $w_*\to 0$, $L_*\to\infty$ 
(for a fixed
 $\bar\mu>0$).  It follows that the Euclidean distance from $V^{m-l+2}_1$
to $\Gamma^{m-l+1}$ is relatively small as compared with its distance from
$\Gamma^{m-l}$ and $\Gamma^{m-l+2}$. More precisely,  there is a 
$\delta=\delta_{\bar\mu}(w_*, L_*)$ with the same properties as 
$\epsilon $ above such that  for any $z\in V^{m-l+2}_1$,
$$\dist(z, \Gamma^{m-l+1})\leq \delta\, 
 \dist (z,\partial (V^{m-l}\backslash V^{m-l+2}) \eqno(5-29)$$

Take $z_0\in V^{m-l}_1$, and let 
 $r=\dist (z_0,\partial (V^{m-l}\backslash V^{m-l+2})$. Note that
the disk $B(z_0, r)$ can be univalently pulled by $g^{l-3}$
to the annulus $V^{m-3}\backslash V^{m-1}$.  By the Koebe Distortion Theorem
and (5-29), for any $\zeta\in V^{m-3}_1$
$$\dist (\zeta, \Gamma^{m-2})\leq C \delta \dist 
(\zeta, \partial(V^{m-3}\backslash V^{m-1}))\leq C\delta\, \dist(\zeta, 0) $$
with an absolute constant $C$.
All the more,
$$\dist (\zeta, \Gamma^{m-2})\leq C \delta\, \diam V^{m-2},$$
so that $\Gamma^{m-2}$ has a big eccentricity about $V^{m-2}_1$
(that is, this eccentricity is at least $e(w_*, L_*)$, where
$e(w_*, L_*)\to \infty$ as $w_*\to 0$, $L_*\to\infty$).

Pulling  $\Gamma^{m-2}$ back by $g_{m+1}\circ g_m\circ g_{m-1}$,
we conclude that $\Gamma^{m+1}$ has a big eccentricity about 0.
Hence it has big width in the annulus $V^m\backslash V^{m+2}$, and 
 Statement 5.7 yields the desired.

Let us summarize the information which will be useful in what follows:

\proclaim  Statement 5.8. If the width $w_{m-l+1}$ is at most $w_*$ and
  the Poincar\'{e} distance from $V_1^{m-l+2}$ to 0 in $V^{m-l+1}$ is 
  at least $L_*$, then the eccentricity $\Gamma^{m}$ about the origin
  is at least $e(w_*, L_*)$, where   $e(w_*, L_*)\to\infty$ as
 $w_*\to 0$ and $L_*\to\infty$.

{\it Let us assume now  that $m-2$ is the last piece of a long cascade}
   $$\CC^{m-2}: V^{m-2-t}\supset\ldots V^{m-2},\quad  t\geq N_*.$$
Then non-central return occurs on level $m-2$. We will show that
$$\sigma_(\II^{m+N+1}|D^{m+N})\geq \sigma_{m-2-t}+a \eqno(5-30)$$
with a definite $a>0$.

Let $D^m\subset V^m$ be the island containing $V^m_0$ and $V^m_1$,
and $\phi_m: D_m\rightarrow V^{m-1}$ be the corresponding two-to-one map.
Note that in the case under consideration 
this island may be non-trivial and still the Poincar\'{e}
distance between $V^{m+1}_0$ and $V^{m+1}_1$ be big (since the pre-critical
puzzle pieces in $V^{m-1}$ are not well inside $V^{m-1}$). Moreover the map
$\phi_m$ is not necessarily bounded perturbation of the
quadratic map. These are the circumstances which make this case special.

As $m-2$ is a non-central level,
$\mu_{m+1}\leq \bar\mu$, and by the previous considerations we are
done unless

\smallskip\noindent$\bullet$ The return on level $m-1$ is Fibonacci, that is
  $\phi_m V^{m+1}_1=V^m_0$ and $\phi_m V^{m+1}_0\subset V^m_1$ for some
  puzzle piece $V^m_1$; 

\smallskip\noindent$\bullet$ The hyperbolic distance between the puzzle
  pieces $V^m_0$ and $V^m_1$ is at least $L_*$. 

\smallskip\noindent$\bullet$ The return on level $m-2$ is also Fibonacci:
  $g_{m-1} V^{m}_1=V^{m-1}_0$ and $g_{m-1} V^m_0= V^{m-1}_1$ for some
  puzzle piece $V^{m-1}_1$.

Let $V^k_0$ and $V^k_1$ be the pieces containing the corresponding push
forwads of $V^{m-1}_0$ and $V^{m-1}_1$ along the cascade $\CC^{m-1}$,
$m-1\leq k\leq m-t-1$. Then (5-30) follows unless

\smallskip\noindent$\bullet$ The width $w_{m-t-3}$ is at most  $w_*$,
  and the distance between  $V^{m-t-4}_0$ and $V^{m-t-4}_1$ in $V^{m-t-3}$
   is at least $L_*$. 

But then by  Statement 5.8  applied to the cascade $\CC^{m-2}$
the eccentricity of $\Gamma^{m-1}$ about 0 is at least 
$e=e(w_*, L_*)$. As $g_m$ is a bounded perturbation of the quadratic map,
 by Lemma A.5 the curve $\Gamma^m$ is
($0.1, \epsilon )$-pinched, where $\epsilon=\epsilon_{\bar\mu}(e)\to 0$
as $e\to \infty$. (Note that the pinched region is not necessarily around
$V^{m}_1$, since $\phi_m$ may differ from $g_m$).  Applying Lemma A.5 again,
we conclude that the curve $\Gamma^{m+1}$ is
 $((10 C)^{-1},C\sqrt{\epsilon})$-pinched.
By Lemma A.3 $\Gamma^{m+1}$ has a definite width inside
$V^m\backslash V^{m+2}$. Now  Statement 5.7 yields (5-30).




\medskip\noindent {\bf  5.6. Other factors yielding big space.}
Theorem III ensures that after  many central cascades we will observe a big
principal modulus.
 However, there are other
combinatorial factors which yield the same effect. Altogether they are 
quite close to a ``big renormalization period", except that ``parabolic or 
Siegel cascades" may interfere.

\medskip\noindent{\bf Big return time implies big modulus.}
This section will  rely on combinatorial considerations of \S 3.8. 
Recall that $\II^n$ denotes the family of puzzle pieces $V^n_i$
intersecting $\omega(0)$. 
Given two $\Upsilon_f$-adjacent puzzle pieces 
$V^{n+1}_j\in \II^{n-1}$ and
$V^n_i\in \II^n$, let $t$ be the first return time of  $V^{n+1}_j$ back
to $V^n_i$ under iterates of $g_n$.  Then we will use the
notation $\mod(V^{n+1}_j\rightarrow V^n_i)$ for 
$\mod(V^n_i\backslash g_n^t V^{n+1}_j).$ If $i\not=0$ then
$$\mod(V^{n+1}_j\rightarrow V^n_i)\geq \mu_n.  \eqno (5-31)$$

\proclaim Lemma 5.9.  Let $D^n\supset V^{n-1}$ be a puzzle piece
containing at least one piece of $\II^n$.
Let $\gamma$ be the shortest path leading
from $D^n$ down to some central piece $V^{n+t}_0$, and let $D^{n+t}$ be the
pull-back of $D^n$ along this path. Then
$$\mod(D^{n+t}\backslash V^{n+t})\geq {\bar\mu\over 2} rank(D^n) .$$ 

\noindent{\bf Proof.} Let 
$D^{n+k}\subset V^{n+k}_{j(k)}  ,\; k=0,1,\ldots, t,$ 
be the pull-back of $D^n$ along the path $\gamma$. By Lemma 3.11, 
all edges of this pull-back except the last one are univalent, and the
last one is a double covering. Hence
$$\mod(D^{n+t}\backslash  V^{n+t})\geq {1\over 2}\sum_{k=1}^{t}
       \mod(V^{n+k}_{j(k)}\rightarrow V^{n+k-1}_{j(k-1)}),$$
and by (5-31) and Theorem II the right-hand side of this 
inequality is estimated
from below by the right-hand side of the desired inequality.  \QED

\proclaim Lemma 5.10. Let $n-1$ be not in the tail of a central cascade.
 Assume that for a puzzle piece $V^{n+1}_j\in \II^{n+1}$,
  $g_n-time(V^{n+1}_j)\geq r$. Then there is a level $m$ such that
  $\mu_m\geq L(r)$, where $L(r)\to \infty $ as $r\to \infty$.

\noindent{\bf Proof.}  Let $M>0$.  We need to find a level $m$ with
 $\mu_m\geq M$, provided
$r$ is sufficiently big. If $rank(V^{n+1}_j)> M/\bar\mu $, then
 Lemma 5.9 yields the desired. So let us assume that
$$rank(V^{n+1}_j)\leq N\equiv M/\bar\mu +1.   \eqno (5-32)$$ 
Let $0=i(0), i(1),\ldots , i(r)=0$ be the itinerary of $V^{n+1}_j$
through the pieces of the previous level.
Let us consider the  nest of puzzle pieces 
$$ D_{r-1}\supset\ldots\supset  D_0\equiv V^n_i,    \eqno(5-33)$$
where $g_n^{r-k} D_{k}=V^{n}_{i(r-k)}$. 
Then 
$$D_{k+1}\backslash D_{k}\geq \bar\mu/2,  \quad   k=1,\ldots r-1.  \eqno(5-34)$$
Let us now pull the pieces $D_k$ down along the shortest path $\gamma$ 
  joining $V^{n+1}_j$ with a critical vertex $V^{n+t}$. 
Denote the corresponding pull-backs by $D^{n+l}_k$.
 If this pull-back
turns out to be univalent then by (5-34) 
$$\mu_{n+t}\geq \mod(D_r\backslash V^{n+1}_j)\geq   r\bar\mu/2,$$ 
which is greater than $M$ for sufficiently big $r$.
Otherwise let us consider the first level $n+s$ where $D^{n+s}_{r-1}$
 hits the critical
point. Let us find such an $l$ that
 $0\in D^{n+s}_{l+1}\backslash D^{n+s}_l$.

If $l> A \equiv 2M/\bar\mu$ then it follows from (5-33) that 
$\mu_{n+s}\geq M$, and we are done.
Otherwise  by (5-34) $\mod( D_l\backslash V^{n+1}_j)\geq  (r-A)\bar\mu$.
Let us now repeat the same procedure with $D_l^{n+s}$ instead of $D_{r-1}$. 
Note that $rank(D_l^{n+s})\leq rank (D_{r-1})-1$, since the pull-back
of $D_{r-1}$ 
along the first central cascade is univalent. Hence this procedure
can be repeated at most $N$ times, and the
principal modulus at the end will  be at least $M$, provided
$(r-NA)\bar\mu\geq 2 M$.   \QED

%

Let us improve this lemma in the case of central cascades.  Below 
we will use notions from \S 3.8.

\proclaim Lemma 5.11. Let us consider a central cascade (3-7)
 such that $m-1$ is not in the tail of a central cascade.
  Assume that $G_{m+N}-time(V^{m+N+1}_j)\geq r$ for some
  non-pre-critical puzzle piece  $V^{m+N+1}_j$. Then 
  there is a level $m$ such that
  $\mu_m\geq M(r)$, where $M(r)\to \infty $ as $r\to \infty$.

\noindent{\bf Proof.} The argument is the same as for the previous lemma
  except one modification:  To construct a  nest (5-33), use
the Bernoulli map $G_{m+N}$ instead of $g_{m+N}$.  \QED

Note that if $V^{m-1}\supset V^m\ldots\supset V^{m+N}$ is a central 
  cascade then the  first condition is satisfied for the sub-cascade
  $V^m\ldots\supset V^{m+N}$.

\medskip\noindent{\bf Parabolic and Siegel cascades.}
We will show that we usually will observe a big principal modulus 
after just one long central cascade. 
Let us consider a central cascade:
$$ V^{m}\supset...\supset V^{m+N-1}\supset V^{m+N},  \eqno (5-35)$$ 
where $ g_{m+1} 0\in V^{m+N-1}\backslash V^{m+N}$. 
The double covering $ g_{m+1}: V^{m+1}\rightarrow V^m$ can be viewed  as
a small perturbation of a quadratic-like map $g_*$ with a definite modulus
and with non-escaping critical point. 

To make this precise, let us consider the space $\QQ$ of double coverings
$g: U'\rightarrow U$, with $0\in U'\subset\subset U$ and $g'(0)=0$,
modulo affine conjugacy.
Let us supply it with
the {\it Catath\'{e}odory topology} (see  [McM]).
Convergence in this topology  means  Carath\'{e}odory convergence of the
domains and the ranges, and uniform convergence of the 
maps  on compact subsets. 


Given a $\mu>0$, let $\QQ(\mu)$ denote the set of double coverings
 $g\in \QQ$ with mod$(g)\geq \mu$. 
By Theorem II, the return maps $g_{m+1}: V^{m+1}\rightarrow V^m$ 
of the principal nest belong to $\QQ(\bar\mu)$.

\proclaim Compactness Lemma (see [McM]). The  set $\QQ(\mu)$ is
  Carath\'{e}odory compact. 

Let $\QQ_N(\mu)$ denote the space of double coverings $g: U'\rightarrow U$ 
from  $\QQ(\mu)$
such that $g^n 0\in U, n=0,1,\ldots N$.
Note that $\QQ_{\infty} (\mu)$ is the space of DH quadratic-like maps
with modulus at least $\mu$.

 As $\bigcap_N \QQ_N(\mu)=\QQ_{\infty}(\mu)$, for any 
neighborhood $\UU\supset \QQ_{\infty}(\mu)$, there is
an $N$ such that $\QQ_N(\mu)\subset \UU$. In this sense
any double map $g\in \QQ_N(\mu)$ is close to some quadratic-like map
$g_*$.
In particular, this concerns the above return map $g_{m+1}$
generating the cascade (5-35)
of big length $N$. Moreover, since $g_{m+1}$ has an escaping fixed point,
the neighborhood of $g_*$ containing $g_{m+1}$ also contains a quadratic-like
map with hybrid class $c(g_*)\in \partial M$. 

If we have a sequence of maps $f_n\in \QQ$ converging to a map 
$g_*\in \partial M$, we also say that the  $f_n$-{\it central cascades
converge to the $g_*$-cascade}.   

Let us say that the principal nest is minor modified 
if a piece $V^m$ is replaced by a piece $\tilde V^n\subset V^n$ such that
$\cl V^{n+1}_i\subset \tilde V^n$ for all pieces $V^{n+1}_i\in \II^{n+1}$.

\proclaim Lemma 5.12. Let $g_*$ be a DH quadratic-like map with
  $c(g_*)\in \partial M$ which does not have neither
  parabolic points, nor Siegel disks. Let $g_{m+1}$  be the return map
 of the principal nest generating cascade (5-35). Take an
arbitrary big $M>0$.
 If $g_{m+1}$ is sufficiently close to $g_*$ 
(depending on a priori bound $\bar\mu$ from Theorem II)
then  the principal nest can be minor modified in such a way that 
 $\tilde\mu_n\geq  M$ for some $n>m+N$.

\noindent{\bf Proof.} Take a big number $e>0$.

By the above assumptions, the Julia set
$J(g_*)$ has empty interior.  If $g_{m+1}$ is sufficiently close to $g_*$ 
then  $\Gamma^{m+N-1}=\partial V^{m+N-1}$
is close in the Hausdorff metric to the Julia set $J(g_*)$.
Hence  $\Gamma^{m+N-1}$ has an eccentricity at least $e$
 with respect to any point $z\in V^{m+N-1}$.

As the $g_m$ are quadratic maps up to bounded distortion (Theorem II),
 the curves
$\Gamma_{m+N}$, $\Gamma_{m+N+1}$ and $\Gamma_{m+N+2}$ also have big eccentricity
with respect to any enclosed point. Moreover, by the same theorem,
there is a definite
space in between these two curves. Hence  by
 Lemma A.2,
 mod$(V^{m+N+1}\backslash V^{m+N+3})$
is at least $M(e)$ where $M(e)\to  \infty$ as $e\to\infty$.

Let us  assume that non-central return occurs on level 
$m+N+1$: $g_{m+N+2} 0\in V^{m+N+2}_i$ with $i\not=0$. 
As the map $g_{m+N+2}:  V^{m+N+2}_i\rightarrow V^{m+N+1}$ is 
quadratic up to bounded distortion, the curve $\Gamma^{m+N+2}_i=
\partial   V^{m+N+2}_i$ has a big eccentricity $e'$ about any enclosed point
(that is, $e'$ can be made arbitrary big by a sufficiently big choice 
of $e$, depending on a priori bound $\bar\mu$). By  Lemma A.2 
$$\mod(V^{m+N+1}\backslash g_{m+N+2} V^{m+N+3}_0)\geq M(e), $$
where $M(e)\to \infty$ as $e\to\infty$. Hence
$\mod A^{m+N+3}\geq M(e)/2,$
and we are done.  

Let  the central return occurs on level   $m+N+1$   but this is not yet 
a DH-renormalizable level. Then the corresponding central cascade
is finite. Let $m+N+T$ be the last level of this cascade.
  Then by Statement 5.8 and Lemma A.2, $\mu_{m+N+T+2}\geq M(e)$, where
 $M(e)\to \infty$ as $e\to\infty$.
   
Assume  finally  that    $m+N+1$   is a DH-renormalizable level.
Then let us take a horizontal curve $\Gamma\subset A^{m+N+2}$
which divides this annulus into two subannuli of moduli at least
$\bar\mu/2$. Let $\Gamma'\subset A^{m+N+3}$ be its pull-back
by $g_{m+N+3}$, and $\tilde A$ be the annulus bounded by 
$\Gamma$ and $\Gamma'$. Then by Lemma A.2 
$ {\rm mod}(\tilde A) \geq M(e)$ with $M(e)$ as above.
As this is a minor modification of the nest, we are done. \QED


\medskip\noindent {\bf 5.7. Proof of Theorem IV.} Let us fix a $Q>0$.
Take a truncated secondary limb $L\equiv L_b^{tr}$, and  find
$q=C(Q)\, \nu(L)$ from  Theorem I. 
 Note that $q(L,Q)\geq \nu(L)/2$ for  sufficiently big $Q$ 
(independently of $L$). Let us now select all
copies $M'$ of the Mandelbrot set with the height 
$\chi(M')\geq  Q/B$, where   $B=B(q)$ is the constant from 
Theorem III. Taking the union of all these copies over all
truncated limbs, we obtain a desired special family $\MM$.

Let us now consider an infinitely renormalizable DH quadraticlike map $f$
of $\SS$-type with $\mod(f)\geq Q$ (to start with, take a quadratic
polynomial). Then by Theorem I, $\mod (A^1)\geq q$. Hence by Theorem III,
$\mod (Rf)\geq B\,\chi(f)\geq Q.$ 

By induction, 
$\mod (R^n f)\geq R$ for all $n$.    \QED

\medskip\noindent{\bf 5.8. Variations.}
Let us now improve Theorem IV by taking into account not only the height
but also the other factors yielding big space. \medskip

\noindent{\bf Theorem IV}$'$. {\it  Let $f$ be an infinitely renormalizable 
quadratic  polynomial, and let $P_m: z\mapsto z^2+c_m$
  be the straightened $R^m f$. Assume that

\noindent $\bullet$ All $c_m$ are selected from
  a finite number of truncated secondary limbs $L_i$, $i=i,\ldots s$;
 
\noindent $\bullet$ The  set $\AA\subset \QQ$ of accumulation points 
 of the central cascades of
  $P_m$ (of lengths growing to $\infty$)
  does  not contain parabolic or Siegel maps;

\noindent $\bullet$ $per (R^m f)\geq p$.

Then $\lim\inf_{n\to\infty}\mod(R^n f)\geq Q(p)$,  where the function 
 $Q(p)$ depends on the choice of the limbs and the accumulation set $\AA$,
and $Q(p)\to \infty$ as $p\to\infty$.}


\noindent{\bf Proof.}  By Theorem II the top modulus of the central
cascades of $P_m$ is bounded from below by some $\bar\mu$. Hence the set
 $\AA\subset \QQ(\bar\mu)$ is compact.
By Lemma 5.12,  for any $Q$ there is a neighborhood  $\UU\supset\AA$
such that: If $f\in \UU$ is renormalizable then $\mod Rf>Q$.

 As $\AA$ is
the accumulation set for the central cascades of the $P_m$,
there is an $N$ such that all but finitely many of these cascades
 of length  $\geq N$ belong to $\UU$. Hence
if the principal nest of $P_m$ contains a cascade of length
 $\geq N$ then $\mod(R (P_m))\geq Q$ (for sufficiently big $m$). 

Further, by Theorems I and III, there is a $\chi$ such that
if the height $\chi(P_m)\geq \chi$ then $\mod (R (P_m))\geq Q$. 
Let us also find a $T$ such that if 
for some cascade the return time from Lemma 5.11
is at least $T$, then $\mod (R (P_m))\geq Q$.

It is easy to see that there 
 is a $p$ such that:  If $per(P_m)\geq p$ then either $P_m$
has a central cascade of length at least $N$, or $\chi(P_m)\geq \chi$,
or one of the above return times is at least $T$.
In any case $\mod(R (P_m))\geq Q$.

Now the same argument as for Theorem IV  yields a priori bounds.  \QED

\bigskip
\centerline{ \bf \S 6. Local connectivity of the Julia sets.}\bigskip
  
In this section we will show that the Julia sets of quadratic polynomials
from Theorems IV and  IV$'$  are locally connected. This follows 
from  the moduli bounds in the full principal nest, which make puzzle pieces
shrink to points
(compare Yoccoz (see [H]), and  Hu-Jiang [HJ], [J]). 
 I thank J.~Kahn and C. McMullen for useful discussions of this issue.

\proclaim Theorem V. Let $\SS$ be a special family of Mandelbrot copies
  from Theorem V. Let $f$ be an infinitely renormalizable quadratic
  of $\SS$-type. Then the Julia set $J(f)$ is locally connected.

\noindent{\bf Proof.}
Let us consider the full principal nest (3-9).
Let $f_m \equiv R^m f$ and $J_m\equiv J(f_m)$. 
 It follows from  Theorem IV and the Gr\"{o}tcsz inequality that 
 $$\mod(Y^{0,0}\backslash J_m)\geq \sum_{k=0}^{m-1}\mod A^{m,0}\geq \epsilon m
\to\infty \quad as \quad m\to\infty.$$
 Hence  the ``little" Julia sets $J_m$ shrink down to the critical point.
Let us take an $\delta>0$, and find an $m$ such that $J_m$ is
contained in the $B(0,\delta)$.

Let us now inscribe into $B(0,\delta)$ a domain bounded by equipotentials and
external rays of the original map 
$f$ (compare Hu and Jiang [HJ], [J]).
Let $\alpha_m$ denote
the dividing fixed point of the Julia set $J_m$, and $\alpha_m'=-\alpha_m$
be the symmetric point. Let us consider a puzzle piece $P_m^0\ni 0$
bounded by any equipotential and four
external rays of the {\it original map $f$}
 landing at $\alpha_m$ and $\alpha_{m}'$. This is 
a ``degenerate" domain of the renormalized map $f_m$ (see \S 2.5). By
definition of the renormalized Julia set, the preimages 
$P_m^k\equiv f_m^{-k} P_m^0$  shrink down to  $J_m$.
Hence there is a puzzle piece $P_{m}^l$
 contained in the $\delta$-neighborhood
of the critical point. As $J(f)\cap P_{m}^l$ is clearly connected,
the Julia set $J(f)$ is locally connected at the critical point.

Let us  now prove local connectivity at any other point $z\in J(f)$.
This is done by a standard spreading  of the local information
near the critical point  around the whole dynamical plane. 
Let $f_m:  U_m'\rightarrow U_m$, where 
$U_m\equiv V^{m,t(m)}$, $U_m'\equiv V^{m,t(m)+1}$ are domains
 from the principal nest (3-9).
Let us consider two cases.\smallskip

{\sl Case (i)}. Let the orbit of $z$ accumulates on all Julia sets
  $J_m$. Find an $l=l(m)$ such that $P_m^{l}\subset U_m'$, and then
take the first moment $k=k(m)\geq 0$ such that $f^kz\in P_m^l$.
Let us consider the pull-back $V_m\supset Q_m^l\ni z$ of $U_m\supset P_m^l$
 along the orbit $\orb_k(z)=\{z,...,f^k z\}$. 
By Lemma 3.3, the pull-back of  the
puzzle piece $P_m^l$ is univalent. Moreover, by construction of the
principal nest, $U_m\backslash J_m$ does not intersect the critical set $\omega(0)$.
Hence the pull-back of  $U_m$ along the $\orb_k(z)$ is also univalent.

Let $\Gamma_m\subset U_m\backslash U_m'$ be a horizontal curve in the
annulus $U_m\backslash U_m'$ which divides it into two sub-annuli of
modulus at least $\epsilon/2$.  By the Koebe Theorem, it has a 
bounded eccentricity about 0 (with a bound depending on $\epsilon$).
Applying  Koebe again, we conclude that its pull-back $\gamma_m$ along
$\orb_k(z)$ has a bounded eccentricity about $z$. Since the inner
radius of this curve about $z$ tends to 0 as $m\to\infty$
(follows from the fact that the sufficiently high iterates
 of any disk in $J(f)$ cover the whole $J(f)$), the $\diam\gamma_m\to 0$
as well. All the more, the $\diam Q_m^l\to 0$ as $m\to \infty$.
As $Q_m^l\cap J(f)$ are connected, the Julia set is locally connected
at $z$.

\smallskip{\sl Case (ii)}. Assume now that the orbit of $z$
 does not accumulate on  some $J_m$. 
Hence it accumulates on some point $a\not\in \omega(0)$.
  Let us consider the puzzle associated with the periodic point $\alpha_m$
(so that the initial
 configuration consists of a certain  equipotential and the
 external rays landing at $\alpha_m$). 
 Since the critical puzzle pieces shrink to 
$J_m$,  the puzzle pieces  $Y^{l}_i$ of sufficiently big depth $l$
containing $a$ are
 disjoint from $\omega(0)$ (there are several such pieces if  $a$
is a preimage of $\alpha_m$). Take such an $l$, and 
let $X$ be the union of these
puzzle pieces. It is a closed topological disk disjoint from $\omega(0)$ whose
interior contains  $a$.  
Hence there is a simply connected
neighborhood $V\supset X$ still disjoint from
$\omega(0)$. Consider now the moments $k_i\to\infty$ when the orbit of
$z$ lands at $\inter X$, and pull $V\supset X$ back to $z$. As the pull-backs 
of $V$ are univalent, the pull-backs of $X$ shrink to $z$ (by the
Koebe argument as above). 
It follows that $J(f)$ is locally
connected at $z$.           \QED

\bigskip
\centerline{\bf \S 7. Forthcoming notes.}\medskip

Let us briefly outline the content of the forthcoming notes. They will be
mostly based on already existing preprints:

\medskip\noindent$\bullet$ Rigidity of quadratics of sufficiently big
  height (\S 4 of [L4]).

\medskip\noindent$\bullet$ Parapuzzle geometry: the moduli
in the principal  parameter nest grow  at the same rate as in the dynamical
  nest. This also yields the 
  rigidity of the corresponding quadratics. 
The Fibonacci case has been worked out by LeRoy Wenstrom (in preparation). 

\medskip\noindent$\bullet$ Geometry of real quadratics ([L3] and [L4], \S 5).
  We give a criterion when the geometry of a real quadratic
  polynomial is ``essentially bounded": It happens if and only if its "essential
period" is bounded.  
 On each level with sufficiently high
   essential period  the renormalized map has a big modulus.

\medskip\noindent$\bullet$ We prove local connectivity of the Mandelbrot set
at all real infinitely renormalizable points with sufficiently high essential 
period on all levels (following \S 5 of [L4]).
Rigidity for all  real quadratics follows 
  (compare Swiatek [Sw]).

\medskip\noindent$\bullet$ Extension of Sullivan's complex a priori bounds
  onto infinitely renormalizable 
 quadratics of essentially bounded type (joint work with
  Michael Yampolsky [LY]). Together with [L4] this
  yields complex bounds (and hence local connectivity of the Julia sets)
  for all real quadratics.

\medskip\noindent$\bullet$ Teichm\"uller metric on the space of
  quasi-quadratic maps following [L5], [L6].

\medskip\noindent$\bullet$  Applications to complex and real
    measurable dynamics.

\bigskip

\centerline{\bf Appendix: Conformal maps and geometry of curves.}\medskip

\noindent{\bf A.1. Poincar\'{e} metric and distortion.} 
A domain $D\subset\C$ is called {\it hyperbolic} if its universal
covering space is conformally equivalent to the unit disk. This happens
if and only if $\C\backslash D$ consists of at least two point. Hyperbolic 
domains possess the {\it hyperbolic (or Poincar\'{e})  metric} $\rho_D$
 of constant negative curveture.
This metric is obtained by pushing down the Poincar\'{e} metric
$d\rho_{\D}=| dz/(1-z^2)|$ from the unit disk $\D$.

 In the case of a simply connected
hyperbolic domain $D$ (``conformal disk"),
 $d\rho_D=p_D(z)|dz|$ is the pull-back of the $\rho_{\D}$
by the Riemann mapping $D\rightarrow \D$. In this case its density $p_D(z)$
is comparable with $1/\dist(z, \partial D)$:
$$ (1/4) \dist(z, \partial D)^{-1}\leq   p(z)\leq \dist(z, \partial D)^{-1}
                                        \eqno (A-1).$$

In the simply connected case,
  a set $K\subset D$ has a bounded hyperbolic diameter $\diam_D K$
if and only if there is an annulus $A\subset D$
of  definite modulus surrounding $K$. More precisely,
let $\mu_{min}(R)$ and $\mu_{max}(R)$ denote the minimal and maximal 
possible modulus of an annulus $A\subset D$ surrounding $K$, where
$K$ runs over all subsets of hyperbolic diameter $R$.
Then $0<\mu_{min}<\mu_{max}<\infty$ (all estimates are clearly independent. 
of $D$). 
This can be readily seen by passing to the disk
model and moving one point of $K$ to the origin.
The extremal moduli correspond to the cases of a pair of points and
hyperbolic disk of radius $R$. 
 Moreover, both minimal an maximal moduli behave  as 
 $\log (1/R)+0(1)$
(see [A], Ch. III).

Given
a univalent holomorphic function  $f: D\rightarrow \C$, the
{\it distortion} of $f$ on $ K$ is defined as
$$\sup_{z,\zeta\in K} \log \left|  f'(z)\over f'(\zeta)\right|.$$

 \proclaim Koebe Distortion Theorem. Let $D$ be a conformal disk,
$K\subset D$, $r=\diam_D K$ be the Poincar\'{e} diameter
of $K$ in $D$. Then the distortion of any univalent function  $f$ on $K$
is bounded by a constant $C_D(r)$ independent of a particular choice of $K$.
Moreover $C_D(r)=O(r)$ as $r\to 0$.

\noindent{\bf A.2. Moduli defect  and capacity.} 
Let $D$ be a topological disk, $\Gamma=\partial D$,
 $a\in D$, and $\psi: (D,a)\rightarrow
({\bf D}_r,0)$ be the Riemann map onto a round disk of radius $r$
with $\psi'(a)=1$.
Then $r\equiv r_a(\Gamma)$ is called the {\sl conformal radius} of $\Gamma$ 
about $a$. The {\sl capacity} of $\Gamma$ rel $a$ is defined as
$${\rm cap}_a(\Gamma)=\log r_a(\Gamma).$$

\proclaim Lemma A.1. Let $D_0\supset D_1\supset K$, where $D_i$ are topological
disks and $K$ is a connected compact. Assume that 
the hyperbolic diameter of $K$ in $D_0$ and the hyperbolic 
dist$(K,\partial D_1)$ are both bounded by a $L$.
  Then there is an $\alpha(L)>0$
such that
$${\rm mod }(D_1\backslash K)\leq {\rm mod} (D_0\backslash K)-\alpha(L).  $$

\noindent{\bf Proof.}  Let us take a point $z\in \partial D_1$ whose hyperbolic
distance to $K$ is at most $L$. Then there is an annulus  of a definite
modulus contained in $D_0$ and enclosing both $K$ and $z$. 

Let us uniformize $D_0\backslash K$ by a round annulus
 $A_r=\{\zeta : r<|\zeta|<1\},$
and let $\tilde z$ correspond to $z$ under this uniformization. Then 
$\tilde z$ stays a definite Euclidian distance $d$ from the unit circle.

If $R\subset A_r$  is any annulus enclosing the inner boundary of $A_r$
but not enclosing $\tilde z$ then by the normality argument
mod$(R)<{\rm mod}(A_r)-	\alpha_r(d)$ with an $\alpha_r(d)>0$.
(Actually, the extremal annulus is just $A_r$ slit along the radius from
$\tilde z$ to the unit circle).

We have to check that $\alpha_r(d)$ is not vanishing as $r\to 0$.
Let us fix an outer boundary $\Gamma$ of $B$ (the unit circle + the slit 
in the extremal case). We may certainly
assume that the inner boundary coincides with the $r$-circle.
Then the defect mod$(R)-\log(1/r)$ monotonically increases to the
 cap$_0(\Gamma)$. By normality this capacity is   bounded above
 by an $-\alpha(d)<0$, and we are done.      \QED

Let $A$ be a standard cylinder of finite modulus, $K\subset A$. Define
the  $\width(K)\equiv\width(K|A)$) 
as the modulus of the smallest concentric sub-cylinder 
$A'\subset A$ containing $K$.

\proclaim Definite Gr\"{o}tzsch Inequality. Let $A_1$ and $A_2$
be homotopically non-trivial disjoint topological annuli in $A$.
Let $K$ be the set of points in their complement which are separated 
by $A_1\cup A_2$ from
the boundary of $A$. Then there is a function $\beta(x)>0$
($x>0$) such that
$${\rm mod }(A)\geq {\rm mod}(A_1)+{\rm mod}(A_2)+\beta({\rm width}(K)).$$

\noindent {\bf Proof.} For a given cylinder this follows from the
usual Gr\"{o}tzsch Inequality and the normality argument. 
Let us fix a $K$, and let  mod$(A)\to\infty$. We can assume that $A_i$
are lower and upper components of $A\backslash K$ correspondingly.
Then the modulus defect
$${\rm mod }(A)-{\rm mod}(A_1)-{\rm mod}(A_2)$$
 decreases by the usual Gr\"{o}tzsch inequality. At the limit the cylinder
becomes the punctured plane, and the modulus defect converges to
-(cap$_0(K)+$cap$_{\infty}(K)$). 

It follows from the area inequality that
this sum of capacities is negative, unless $K$ is a circle centered at the
origin. Moreover the estimate depends only on $\width(K)$. 
Indeed, let $D_0$ and $D_{\infty}$ be the components of $\C\setminus K$
containing 0 and $\infty$ correspondingly. Let $\phi_0: B(0, R_0)\rightarrow
D_0$ and $\phi_{\infty}: \C\backslash B(0, R_{\infty})\rightarrow D_{\infty}$
be the Riemann mappings normalized by: $\phi_0(z)\sim z$ as $z\to 0$, and
$\phi_{\infty} (z)\sim z$ as $z\to\infty$. 
Then $\capac_0 K=\log R_0$ and $\capac_{\infty} K= \log (1/R_{\infty})$.

As scaling does not change  $\capac_0(K)+\capac_{\infty} (K)$, we can assume that
$R_{\infty}(K)=1$. Let $$\phi_{\infty}(z)=z-\sum_{k=1}^{\infty}{a_k\over z^k}.$$
Then:
$$\area(\C\backslash D_{\infty})=
{i\over 2}\int_{|z|=1} \phi_{\infty}d\bar\phi_{\infty}=
\pi(1-\sum_{k=1}^{\infty} k|a_k|^2)\leq \pi,$$
with equality only in the case when $\phi_{\infty}=\id$.
Hence $\area(D_0)\leq \pi$ with equality only in the case when $K$ is the unit
circle $S^1$.

As $|\phi_0|^2$ is a subharmonic function, 
$$1=|\phi_0(0)|^2\leq {\area(D_0)\over \pi R_0^2}\leq {1\over R_0^2},$$
with equality only on the case when $K=S^1$. Hence $\capac_0(K)<0$
unless $K=S^1$. Moreover, by a normality argument 
$\capac_0(K)\leq c(\width(K))<0$. (Indeed, otherwise there would be a sequence
of domains $D_0^m$ as above converging to a domain $\Omega$
different from the unit disk, with $\capac_0(\Omega)=0$.) 
The lemma is proved.  \QED

\medskip\noindent{\bf A.3. Eccentricity and pinching.}

Let $\Gamma$ be a Jordan curve surrounding a point $a$.
Let $d_a(\Gamma)$ and $\rho_a(\Gamma)$ be the Euclidian radii of the 
inscribed and circumscribed circles about  $\Gamma$ 
centered at $a$. Then let us define the eccentricity of $\Gamma$ about $a$ as
$$e_a(\Gamma)=\log{\rho_a(\Gamma)\over d_a(\Gamma)}.$$

\proclaim Lemma A.2. Let $A\subset {\bf C}\backslash \{0\}$ be an annulus
homotopically non-trivially embedded
in the punctured plane,  $\Gamma\subset A$ be 
a homotopically non-trivial Jordan curve, and $A_i$ be the components
of $A\backslash \Gamma$. Assume that mod$(A_i)\geq\mu>0$. If $e_0(K)\geq e$ 
then $\width(\Gamma|A)\geq w(e)$, where $w(e)\to\infty $ as $e\to \infty$.

\noindent{\bf Proof.}  Assume that there is a sequence of annuli $A^n$
and curves $\Gamma^m\subset A^m$ satisfying the assumptions of the lemma,
such that  $\width(\Gamma^m|A^m)\leq w$,
while $e_0(\Gamma^m)\to\infty$.
  Let us  consider the uniformization $\phi_m:\bar A^m
\rightarrow A^m$ of the $A^m$  by round annuli centered at 0. 
Let $\bar \Gamma^m=\phi^{-1} \Gamma^m$.
 Then
$\bar \Gamma^m$ is contained in a round annulus  $\bar R^m$ of modulus $\leq w$
 cocentric with  $\bar A^m$. Let $R^m=\phi \bar R^m$. 

Let us normalize the annuli $A^m$ and $\bar A^m$ (by scaling and rotation)
so that the inner radii of $R^m$ and $\bar R^m$  are equal to 1,
 and $\phi_m(1)=1$.
Passing to a subsequence (without change of notations) we can find
cocentric annuli $\bar R\subset \bar A$ such that the inner radius of 
$\bar R$ is equal to 1,
$\mod \bar R\leq w$, the both components of $\bar A\backslash \bar R$
have moduli at least $\alpha(\mu,w)$, and
 $\bar R^m\subset\bar R$, $\bar A^m\supset \bar A$.

By the Koebe Theorem, the family of functions $\phi_m$ is normal in $A$.
Hence these functions are  uniformly bounded on $\bar R$ contradicting the
assumption that the eccentricities of $\Gamma^m$  about 0 go to $\infty$. \QED
 
``Pinching" of a Jordan curve means creating of a narrow region which
in limit makes the curve non-simple. Below we will quantify this 
process.

Let us take a number $0<k<1$  called the ``pinching parameter".
Let us define the $k$-{\it pinching} of a Jordan curve $\Gamma$ 
as $\xi_k(\Gamma)=\inf \dist(z_1,z_2)$, where the infimum  is taken
  over all pairs of points 
$ z_i\in \Gamma$ such that both components $\Gamma_1$ and $\Gamma_2$  
of $\Gamma\backslash\{z_1, z_2\}$ have diameter at least $k\, \diam\Gamma$.
 
We say that a curve $\Gamma$ is $(k,\epsilon)$-{\it pinched}  if
$\xi_k(\Gamma)<\epsilon$. Note that if the curve is symmetric about 0
and $e_0(\Gamma)\geq e$, then it is $(0.5-e^{-1}, e^{-1})$-pinched.

The following lemma shows that  a sufficiently 
pinched curve has a definite width:

\proclaim Lemma A.3. Let $\Gamma$, $ A$ and $A_i$ be the same objects
as  in Lemma A.2.  Let also  mod$(A_i)\geq\mu>0$. 
If $\Gamma$ is  $(k,\epsilon)$-pinched, then there exists a  $w=w(\mu, k)>0$
such that $\width(\Gamma | A)\geq w>0$ for all  sufficiently
 small $\epsilon>0$. 

\noindent{\bf Proof.} Otherwise we can find a sequence $A^m$ of annuli
as above with $\width(\Gamma^m | A^m)\to 0$ as $m\to\infty$, and
$\Gamma^m$ is $(k, 1/m)$-pinched.  As in the
previous lemma, let 
$\phi_m : \bar A^m\rightarrow A^m$ be the uniformizations by  round 
annuli normalized in such a way that 
$\Gamma$ and $\bar \Gamma=\phi^{-1}\Gamma$ pass through 1.

Then the curves $\bar \Gamma^m$ should converge to the unit circle in
the Hausdorff metric.  Moreover, the family $\phi_m$ is well-defined and
normal on a cocentric  annulus $\bar A$ of modulus, say, $\sqrt{\mu}$.
Hence any Hausdorff limit of the famlily of curves $\Gamma^m$ is an
analytic Jordan curve. On the other hand, these curves should be
$(k,0)$-pinched (that is, they are non-simple). Contradiction.          \QED 

\proclaim Lemma A.4. Let $\Gamma$, $ A$ and $A_i$ be the same objects
as  in Lemma A.2, and  mod$(A_i)\geq\mu>0$. Let $D$ be a topological
disk bounded by $\Gamma$, and $b\in D$. Then there is a function
$\delta(L)\to 0 $ as $L\to \infty$ such that  
$$\dist(b, \bd \Gamma)\leq \delta(L)\; \diam\Gamma,$$
provided $\rho_D(0,b)\geq L$.

(Thus, if $b$ is hyperbolically far away from 0 then it is 
Euclideanly close to the $\bd D$, in the scale of $D$.) \smallskip

\noindent{\bf Proof.}  Otherwise there is a sequence of the curves
$\Gamma^m=\bd D^m$ as above, and points $b^m\in D^m$ such that
$\diam D^m=1$,
$$\rho_{D^m}(0,b_m)\to \infty, \eqno(A-2),$$
  and 
$$\dist (b^m, \Gamma^m) \geq\delta.   \eqno (A-3)$$ Passing to a
Caratheodory limit along some subsequence (without change of notations),
 we have: $(D^m,0, b^m)\to (D,0, b)$, where $D$ is a topological disk
(note that $b\in D$ due to (A-3)). But then 
$\rho_{D^m} (0,b_m)\to \rho_D(0,b)<\infty$, contradicting (A-2).
                                             \QED



\proclaim Lemma A.5. Let $\Gamma$ be a Jordan curve  which does not pass
through 0. If it 
 is $(k,\epsilon)$-pinched then its 
 pull-back under  the quadratic map  $\Phi: z\mapsto z^2$
is $(C^{-1}k, C\sqrt{\epsilon})$-pinched, where $C>1$ is an absolute
constant. 

\noindent{\bf Proof.}
Let $z_1$ and $z_2$ be two points on $\Gamma$
such that $\dist(z_1, z_2)<\epsilon$, while $\diam \Gamma_i> k$, where
$\Gamma_i$ are complementary components of $\Gamma\backslash \{z_1, z_2\}$.
Let us mark the $\Phi$-preimages of the corresponding objects with twilde
(select the closest preimages of the points $z_i$).

We can assume that $\diam \Gamma=1/4$. If $\dist (\Gamma,0)>1/4$ then the
distortion of the quadratic map on $\Gamma$ is bounded by an absolute
constant, and the conclusion follows. Otherwise $\Gamma$ is contained
in the unit disk. Hence $1/2\geq \diam \tilde\Gamma\geq 1/8$ 
and $\diam\tilde \Gamma_i\geq(1/2)\diam \Gamma_i\geq k/8$.
Moreover, $\dist(\tilde z_1, \tilde z_2)\leq \sqrt{\epsilon}$,
and we are done.
                                     \QED

\bigskip\centerline{\bf References.}\bigskip

\noindent
[A] L. Ahlfors. Lectures on quasi-conformal maps. Van Nostrand Co, 1966.
      \smallskip

\noindent [At] P. Atela. Bifurcations of dynamical rays in complex 
  polynomials of degree two. Ergod. Th. \& Dynam. Sys., v. 12 (1991),
 401-423.

\noindent
[Be] A. Beardon. Iteration of rational functions. Graduate texts in mathematics.
    Springer-Verlag, 1991. \smallskip

\noindent
[B] B.Branner. The Mandelbrot set. In: ``Chaos and fractals",
 Proceedings of Symposia in Applied Mathematics, v. 39, 75-106. \smallskip

\noindent
[BH] B.Branner \& J.H.Hubbard.   The iteration of cubic polynomials,
  Part II.  Acta Math. v. 169 (1992), 229-325.\smallskip

\noindent
[CG] L. Carleson \& W. Gamelin. Complex dynamics.  Springer-Verlag, 1993.
  \smallskip

\noindent
 [D1] A. Douady. Shirurgie sur les applications holomorphes. In:
  ``Proc. Internat. Congress  Math. Berkeley", 1986, v. 1, pp. 724-738.
    \smallskip

\noindent
[D2] A. Douady. Description of compact sets in ${\bf C}$.
  In:
    ``Topological Methods in Modern Mathematics, A Symposium in
     Honor of John Milnor's 60th Birthday", Publish or Perish, 1993.\smallskip
 
\noindent
[DH1] A.Douady \& J.H.Hubbard. \'{E}tude dynamique des polyn\^{o}mes
   complexes. Publication Mathematiques d'Orsay, 84-02 and 85-04.\smallskip

\noindent
[DH2] A. Douady \& J.H. Hubbard. On the dynamics of polynomial-like maps.
    Ann. Sc. \'{E}c. Norm. Sup., v. 18 (1985), 287-343.\smallskip

\noindent
[GS] J. Graczyk \& G. Swiatek. Induced expansion for quadratic polynomials.
  Stony Brook IMS Preprint 1993/8.  \smallskip

\noindent
[GM]  L. Goldberg \& J. Milnor. Part II: Fixed point portraits. 
     Ann. Sc. \'{E}c. Norm. Sup.,
\smallskip

\noindent 
[H] J.H. Hubbard.  Local connectivity of Julia sets and bifurcation
     loci: three theorems of J.-C. Yoccoz. In:
    ``Topological Methods in Modern Mathematics, A Symposium in
     Honor of John Milnor's 60th Birthday", Publish or Perish, 1993.\smallskip

\noindent
[HJ] J. Hu \& Y. Jiang. The Julia set of the Feigenbaum quadratic
  polynomial is locally connected. Preprint 1993.\smallskip

\noindent
[J] Y. Jiang. Infinitely renormalizable quadratic Julia sets.
  Preprint, 1993.\smallskip

\noindent
[LS] G. Levin \& S. van Strien. Local connectivity of the Julia sets of 
  real polynomials. Preprint 1995.

\noindent
[L1] M. Lyubich. The dynamics of rational transforms: the topological picture.
     Russian Math. Surveys, v. 41 (1986), \# 4, 43-117.

\noindent
[L2] M. Lyubich.  Milnor's attractors, persistent recurrence and
  renormalization, 
  ``Topological Methods in Modern Mathematics, A Symposium in
     Honor of John Milnor's 60th Birthday", Publish or Perish, 1993.\smallskip

\noindent
[L3] M. Lyubich. Combinatorics, geometry and attractors of
     quasi-quadratic maps. \break   
     Preprint IMS  at Stony Brook \#1992/18,  Annals of Math., v. 140 (1994),
    347-404.
  \smallskip

\noindent
[L4] M. Lyubich. Geometry of quadratic polynomials: moduli, rigidity and
   local connectivity.   
     Preprint IMS  at Stony Brook \#1993/9.\smallskip

\noindent
[L5] M. Lyubich. Teichm\"{u}ller space of Fibonacci maps.
   Preprint IMS at Stony Brook, \#1993/12. \smallskip

\noindent
[L6] M. Lyubich. Teichm\"{u}ller space of quasi-quadratic maps.
   Manuscript in preparation. \smallskip

\noindent 
[LY] M. Lyubich \& M. Yampolsky. Dynamics of quadratic polynomials:
  Complex bounds for real maps, Manuscript 1995.
  \smallskip

\noindent
[LM] M. Lyubich \& J.Milnor. The unimodal Fibonacci map. 
         Journal of AMS, v. 6 (1993), 425-457. \smallskip

\noindent
[M1] J. Milnor. Dynamics in one complex variable: Introductory lectures,
          Stony Brook  IMS Preprint \# 1990/5. \smallskip

\noindent
[M2] J. Milnor. Local connectivity of Julia sets: expository
   lectures. Preprint IMS Stony Brook, \#1992/11.\smallskip

\noindent
[M3] J. Milnor. Self-similarity and hairiness in the Mandelbrot set,
          pp. 211-257 of ``Computers in geometry and topology",
          Lect. Notes in Pure Appl Math, v. 114, Dekker 1989.
         \smallskip

\noindent
[M4] J. Milnor. Periodic orbits, external rays and the Mandelbrot set:
    An expository account.
   Manuscript. \smallskip

\noindent
[McM] C. McMullen. Complex dynamics and renormalization.
  Preprint, 1993.\smallskip

\noindent
[Ma] M. Martens. Rotation numbers for minimal Cantor sets.
          Manuscript in preparation.\smallskip

\noindent
[MS] W. de Melo \& S. van Strien. One dimensional dynamics.
   Springer-Verlag, 1993. \smallskip

\noindent
[P] C. L. Petersen. Local connectivity of some Julia sets containing a circle
   with an irrational rotation, Preprint IHES/M/94/26.

\noindent
[R] M. Rees. A possible approach to a complex renormalization problem.
   In: ``Linear and Complex Analysis Problem Book 3, part II", p. 437-440. 
  Lecture Notes in Math., v. 1574. \smallskip

\noindent
[Sch] D. Schleicher. Internal addresses in the Mandelbrot set and
  irreducibility of polynomials. Thesis, 1994. 

\noindent
[Sw] G. Swiatek. Hyperbolicity is dense in the real quadratic family.
   Preprint IMS at Stony Brook, \#1992/10. \smallskip

\noindent
[S1] D.Sullivan. Quasiconformal homeomorphisms in dynamics, topology
    and geometry.   Proceedings of the ICM, Berkeley, 1986, v. 2,
    1216. \smallskip

\noindent
[S2] D.Sullivan. Bounds, quadratic differentials, and renormalization
conjectures.  AMS Centennial Publications.
{\bf 2}: Mathematics into Twenty-first Century (1992). \smallskip

\end